\RequirePackage{fix-cm}
\documentclass[smallextended]{svjour3}       
\smartqed  
\usepackage{graphicx}
\usepackage{cite}
\usepackage{subcaption}
\usepackage[export]{adjustbox} 

\usepackage{amsmath}
\usepackage{amssymb}

\usepackage{xspace} 
\newcommand{\etal}[0]{{\em et~al.\@}\xspace}
\newcommand{\eg}[0]{{e.g.\@}\xspace}
\newcommand{\ie}[0]{{i.e.\@}\xspace}

\newcommand{\norm}[1]{\left\lVert#1\right\rVert}
\newcommand{\abs}[1]{\left\lvert#1\right\rvert}
\newcommand{\T}[0]{\mathrm{T}} 
\usepackage{mathtools}  
\newcommand{\defn}[0]{\coloneqq} 
\newcommand{\order}[1]{\ensuremath{\mathop{}\!\mathcal{O}\!\left(h^{#1}\right)}}
\DeclareMathOperator{\sgn}{sgn}
\newcommand{\mean}[1]{\left\{\!\!\left\{#1\right\}\!\!\right\}}

\usepackage{xfrac} 
\usepackage{xcolor}
\usepackage[normalem]{ulem} 

\usepackage{hyperref} 
\hypersetup{
    colorlinks,
    linkcolor={black},
    citecolor={black},
    urlcolor={blue!80!black}
}

\spnewtheorem{assumption}{Assumption}{\bfseries}{\itshape}

\newcommand{\bsym}[1]{\ensuremath{\boldsymbol{#1}}}

\newcommand{\fnc}[1]{\ensuremath{\mathcal{#1}}}
\newcommand{\vecfnc}[1]{\ensuremath{\bsym{\mathcal{#1}}}} 

\renewcommand{\vec}[1]{\ensuremath{\bsym{#1}}}

\newcommand{\mat}[1]{\ensuremath{\mathsf{#1}}}

\newcommand{\bu}[0]{\ensuremath{\vec{u}}}
\newcommand{\bv}[0]{\ensuremath{\vec{v}}}
\newcommand{\tl}[0]{\vec{t}_{\mathrm{L}}}
\newcommand{\tr}[0]{\vec{t}_{\mathrm{R}}}
\newcommand{\tlT}[0]{\vec{t}_{\mathrm{L}}^\T}
\newcommand{\trT}[0]{\vec{t}_{\mathrm{R}}^\T}

\newcommand{\Idty}[0]{\mat{I}} 
\newcommand{\Hnrm}[0]{\mat{H}}
\newcommand{\HI}[0]{\mat{H}^{-1}}

\newcommand{\pder}[2][]{\dfrac{\partial #1}{\partial #2}} 
\newcommand{\der}[2][]{\dfrac{\mathrm{d} #1}{\mathrm{d} #2}} 

\newcommand{\diag}[1]{\mathrm{diag} \!\left( #1 \right)}

\begin{document}

\title{On the Practical Impact of Local Linear Instabilities in Entropy-Stable Schemes}


\author{Alex Bercik         \and
        David W. Zingg 
}


\institute{A.\ Bercik \at
              University of Toronto Institute for Aerospace Studies, Toronto, Canada \\
              \email{alex.bercik@mail.utoronto.ca}
           \and
           D.\ W.\ Zingg \at
              University of Toronto Institute for Aerospace Studies, Toronto, Canada \\
              \email{david.zingg@utoronto.ca}
}

\date{Received: date / Accepted: date}

\maketitle

\begin{abstract}

Local linear instability refers to the linearized discrete operator exhibiting perturbation growth exceeding that of the corresponding continuous linearized problem. In the context of nonlinear entropy-stable discretizations, we argue that local linear instabilities should be interpreted primarily as a mechanism for introducing numerical error whose practical impact is often negligible compared with other sources of discretization error. For split-form discretizations of the variable-coefficient linear advection equation, such as those resulting from linearizations of entropy-stable discretizations of the Burgers equation, perturbations can indeed exhibit unphysical modal growth. However, we demonstrate that this growth satisfies physically interpretable bounds and is typically small. Furthermore, through both modified-equation analysis and numerical experiments, we show that the growth is dominated by highly oscillatory and boundary-localized unphysical modes, and can therefore be readily controlled by employing small amounts of numerical dissipation. More generally, this modal perturbation growth does not extend directly to nonlinear two-point-flux discretizations of the type used in entropy-stable discretizations of the Euler equations. Floquet analysis is used to demonstrate that unstable spectra of frozen-baseflow Jacobians do not necessarily lead to unstable perturbation growth. Using the geometric flux for the variable-coefficient linear advection equation, we derive a sharp perturbation growth bound predicting negligible growth, then use numerical experiments to show analogous behaviour for the logarithmic flux. Finally, we argue that the robustness issues observed for entropy-stable discretizations in near-vacuum density-wave problems are better attributed to poor near-vacuum behaviour of the logarithmic mean than to local linear instability itself. Overall, our results suggest that local linear instabilities do not pose a practical obstacle to the use of high-order entropy-stable methods.

\end{abstract}

\section{Introduction}

Entropy-stable methods provide a systematic framework for constructing high-order discretizations with provable nonlinear stability properties for hyperbolic conservation laws and related balance laws~\cite{gassner_review}. Building on the seminal work of Tadmor~\cite{tadmor_2003}, modern constructions employ entropy-conservative two-point fluxes together with appropriate numerical dissipation to obtain a discrete entropy inequality~\cite{leFloch2002,fisher_2013,carpenter2014,chen2017,chan2018}. These ideas have been applied across a broad range of problems, including the compressible Euler and Navier--Stokes equations~\cite{crean,gassner_splitform_2016,chan_entropyprojection,carpenter2014,Fernandez2020}, the shallow water equations~\cite{gassner2016shallow,winters2017shallow}, magnetohydrodynamics and plasma models~\cite{liu2018mhd,bohm2020resistive,biswas2021tenmoment}, atmospheric flows with gravity~\cite{waruszewski2022gravity}, and multiphase or nonconservative systems~\cite{renac2019,ntoukas2022}. A central tool for enabling high-order constructions is the summation-by-parts (SBP) property, which together with simultaneous approximation terms (SATs) for enforcing boundary conditions, provides an algebraic framework for reproducing stability and conservation proofs at the semidiscrete level~\cite{Kreiss1974,strand1994,carpenter1999,DelReyFernandez2014_review}. This framework can be applied to finite-difference~\cite{Kreiss1974,strand1994,Fisher_2013_SBP}, finite-volume~\cite{nordstrom_finite_2001, nordstrom_finite_2003, Chandrashekar_finite_2016}, continuous Galerkin~\cite{Hicken2020, hicken_multidimensional_2016, Abgrall_cg1, Abgrall_cg2}, discontinuous Galerkin~\cite{carpenter2014, chan2018, gassner2013}, flux reconstruction~\cite{Ranocha_fr, Montoya_fr}, Galerkin difference~\cite{yan2023dgd}, and point-cloud or meshfree methods~\cite{hicken_pointcloud,Glaubitz2024}. For more details regarding entropy-stable and SBP methods, we refer the reader to the reviews of~\cite{gassner_review, chen_review, DelReyFernandez2014_review}. 

Despite the practical success of entropy-stable methods, recent work has raised concerns regarding the local linear instability of entropy-stable and split-form schemes~\cite{gassner_stability_2022,ranocha_2022_preventing_pressure}. Local linear instability refers to the possibility that perturbations about a smooth baseflow grow faster in the semidiscretization than in the corresponding continuous linearized problem. This behaviour can be diagnosed by linearizing the nonlinear semidiscretization about a given baseflow state and examining the spectrum of the resulting Jacobian. In model problems for which the continuous linearized problem admits no exponential growth (\eg a time-invariant operator with a purely imaginary spectrum), eigenvalues of the semidiscrete Jacobian with positive real parts indicate unphysical modal growth of perturbations. In~\cite{gassner_stability_2022}, this growth is attributed to anti-dissipative contributions of entropy-conservative two-point fluxes and their related split-form volume terms. These exponentially growing unphysical modes may eventually degrade the physical reliability of the simulation, trigger nonlinear mechanisms such as turbulence, or lead to positivity violations of thermodynamic variables. From this perspective, a global entropy estimate does not by itself guarantee that the scheme accurately reproduces the correct dynamics of the model problem.

The practical significance of this mechanism, however, remains unclear. All numerical discretizations introduce numerical error through the finite-dimensional approximation of continuous operators and solution spaces, aliasing, boundary closures, interface coupling, time integration, and roundoff. Therefore, the relevant question is not merely whether a frozen semidiscrete Jacobian has eigenvalues with positive real parts, but whether the associated perturbation growth is significant enough to affect the numerical solution in comparison with other sources of error. This distinction is especially important in the presence of nonlinear stability estimates: if the unphysical perturbation growth is eventually capped by a global entropy estimate, then local linear instabilities should be interpreted primarily as a mechanism for introducing numerical error, whose practical relevance must be assessed relative to other error sources.

In this work, we study the unphysical growth of perturbations by directly considering the variable-coefficient linear advection equation, which arises from the linearization of scalar hyperbolic conservation laws. The dynamics of this linear model problem, however, can change significantly depending on whether the coefficient varies in time or not. Spectral analysis directly characterizes modal growth in the time-invariant problem, but can be misleading for the time-varying problem~\cite{nonautonomous_systems,Chicone2006}. By contrast, Floquet analysis provides a natural framework for evaluating the stability of time-periodic linear systems~\cite{Chicone2006}. This distinction is also relevant to the local linear stability analysis of entropy-stable schemes, as even for a constant-coefficient model problem, nonlinear discretizations can produce time-varying linear operators when linearized about an evolving numerical solution. Floquet analysis long been used in conjunction with numerical methods to study the stability of differential equations, with applications for example to rotor dynamics~\cite{Peters1971,friedmann_1977}, fluid flows~\cite{Barkley_Henderson_1996}, and delay differential equations~\cite{Insperger_2002}. In these works, the numerical approximation is primarily used to approximate the Floquet multipliers of an underlying model problem. Here, we adapt these techniques to study the stability properties of the numerical discretizations themselves. A related perspective appears in~\cite{Higham1994}, where explicit Runge--Kutta time-marching methods were analyzed for a time-periodic ordinary differential model equation, and spectra of frozen Jacobians were shown to give misleading stability predictions.

The remainder of this paper is organized as follows. In \S\ref{sec_continuous_VCLCE}--\S\ref{sec_stability_proof}, we introduce the variable-coefficient linear advection model problem, formulate the flux-differencing SBP discretizations used throughout the paper, and establish basic stability, convergence, and perturbation bounds for the non-central schemes. In \S\ref{sec_modified_pde} and \S\ref{sec_frozen_coeff}, we revisit the local linear stability arguments of~\cite{gassner_stability_2022} using modified-equation and frozen-coefficient analysis to gain insight into the structure of locally linearly unstable modes. To confirm these observations, in \S\ref{sec:spectra} we analyze the spectra of time-invariant split-forms, and in \S\ref{sec:time_invariant_accuracy_experiments} perform numerical experiments to quantify the practical impact of perturbation growth relative to other sources of numerical error. In \S\ref{sec:time-varying}, we turn to nonlinear flux-differencing schemes whose linearizations lead to time-varying perturbation equations. Floquet analysis is introduced in \S\ref{sec:floquet} as a diagnostic for these time-periodic semidiscrete problems, and is then compared with numerical experiments in \S\ref{sec_logarithmic_LCE_experiments} to assess accumulated perturbation growth. The density-wave problem of \cite{gassner_stability_2022}, which is known to be challenging for entropy-stable schemes, is then studied in \S\ref{sec:nearvacuum} through our model problem, drawing attention to the numerical issues introduced by logarithmic fluxes in near-vacuum regions. Finally, in \S\ref{sec:nonlinear} we explore the broader implications of our results in the nonlinear context for the Burgers and compressible Euler equations.

\section{Time-Invariant Variable-Coefficient Linear Advection} \label{sec_lin_advec}

\subsection{The Continuous Problem} \label{sec_continuous_VCLCE}

\begin{figure}[t]
    \centering
    \begin{subfigure}[t]{0.44\textwidth}
        \includegraphics[height=3.8cm,
        trim={3mm 3mm 3mm 3mm},clip]{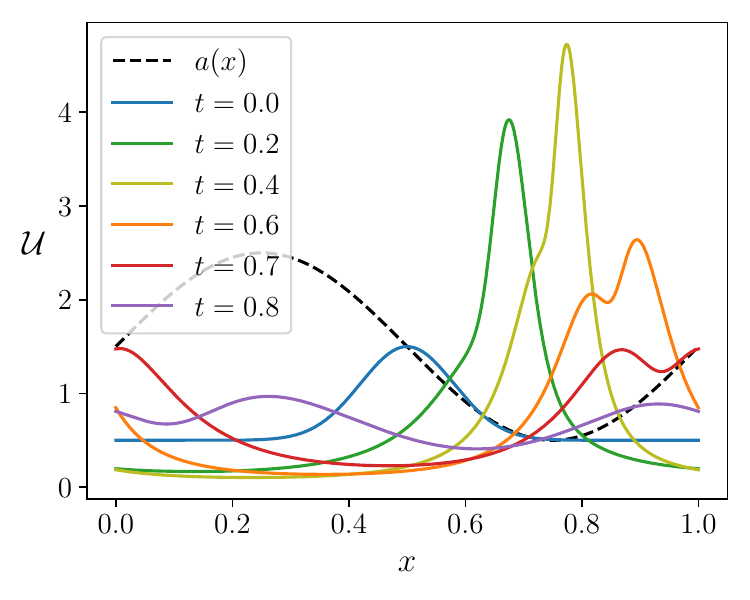}
    \end{subfigure}
    \hfill
    \begin{subfigure}[t]{0.55\textwidth}
        \includegraphics[height=3.8cm,
        trim={3mm 3mm 3mm 3mm},clip]{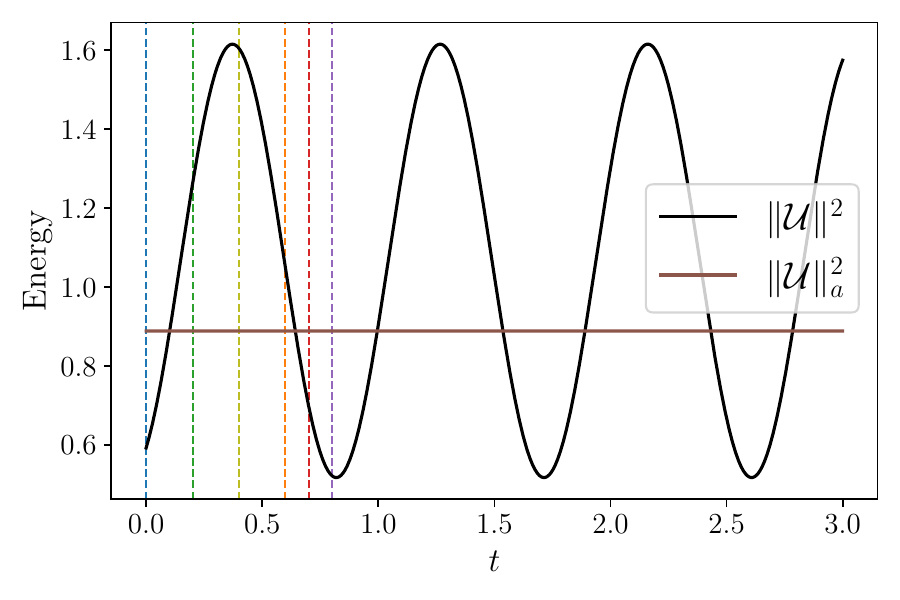}
    \end{subfigure}

    \caption{The exact solution and energy (in the $L^2$ and $a$-weighted norms) for the variable-coefficient linear advection equation \eqref{eq_lin_conv} with a Gaussian initial condition and sinusoidal coefficient $a(x)=\sin{(2 \pi x)}+\sfrac{3}{2}$. The times at which the solution snapshots are plotted in the left plot are indicated in the right plot by correspondingly coloured vertical dashed lines.}
    \label{fig:VCE_ex}
\end{figure}

The majority of this work considers the variable-coefficient linear advection equation,
\begin{equation} \label{eq_lin_conv}
    \pder[\fnc{U}]{t} + \pder[a(x) \fnc{U}]{x} = 0 , \quad (x,t) \in [0,L] \times [0,\infty), \quad a(x) > 0,
\end{equation}
with periodic boundary conditions and some initial condition $\fnc{U}(x,t=0)=\fnc{U}_0(x)$. As measured in the $L^2$ norm, the variable coefficient can induce constrained energy growth (or decay, depending on the sign of $\partial a / \partial x$), seen through the energy balance $\frac{\mathrm{d}}{\mathrm{d}t} \norm{\fnc{U}}^2 = - \int_0^L \frac{\partial a(x)}{\partial x} \fnc{U}^2 \, \mathrm{d}x$ \cite{kopriva_energy_2014, kreiss_lorenz}. However, this energy production is norm-dependent, and for periodic boundary conditions, remains dynamically balanced as the solution is advected through alternating regions of growth and decay. For example, the $a$-energy, $\norm{\fnc{U}}_a^2 \defn \int_0^L a(x) \fnc{U}^2 \, \mathrm{d}x $, is conserved for periodic boundary conditions~\cite{manzanero}. An illustrative example is shown in Figure~\ref{fig:VCE_ex}, where we take $L=1$, a sinusoidal coefficient $a(x)=\sin{(2 \pi x)}+\sfrac{3}{2}$, and a Gaussian initial condition of $\fnc{U}_0(x) = e^{-0.5 \left( \frac{x-0.5}{0.08} \right)^2} + \sfrac{1}{2}$. The solution is periodic in time with a period of $T=\sfrac{2}{\sqrt{5}}\approx0.9$, initially growing as the Gaussian is centred in a region of decreasing $a(x)$, before decaying as the wave convects through the region of increasing $a(x)$. Correspondingly, the standard energy norm oscillates in time, whereas the weighted $a$-energy illustrates the well-balanced nature of the growth.

By equivalence of norms, conservation of the $a$-energy implies that the standard $L^2$ energy is uniformly bounded in time~\cite{manzanero}. For the periodic problem~\eqref{eq_lin_conv}, we have
\[
\frac{\min_{x \in [0,L]}(a(x))}{\max_{x \in [0,L]}(a(x))} \norm{\fnc{U}_0}^2 \leq \norm{\fnc{U}(t)}^2 \leq \frac{\max_{x \in [0,L]}(a(x))}{\min_{x \in [0,L]}(a(x))} \norm{\fnc{U}_0}^2 .
\]
Since \eqref{eq_lin_conv} is a linear well-posed problem, uniform boundedness implies that the differential operator has a purely imaginary spectrum~\cite{engel_nagel_2000}. It is widely assumed that numerical schemes should respect the spectral bounds of the continuous operator to ensure both numerical stability and structural consistency~\cite{gassner_stability_2022, manzanero}. For linear numerical semidiscretizations of the form $\mathrm{d} \bu / \mathrm{d}t = \mat{J} \bu$ with $\mat{J}$ independent of $\bu$, the analogous condition is that the eigenvalues of the semidiscretization Jacobian $\mat{J}$ have zero real parts, \ie $\Re{(\lambda(\mat{J}))} = 0$. In practice, eigenvalues with negative real parts are permitted as they do not contribute to modal growth, resulting in a desired spectral bound of the form~${\Re{(\lambda(\mat{J}))}_{\max} \leq 0}$. In this work, we will argue that for nonlinear schemes possessing alternative stability estimates, such as entropy-stable schemes, satisfying such spectral bounds is often not necessary.

\subsection{Flux-Differencing SBP Schemes} \label{sec:flux_differencing}

Given a diagonal-norm SBP operator $\mat{D} \in \mathbb{R}^{N \times N}$, a typical flux-differencing volume term is (\eg \cite{fisher_2013, gassner_splitform_2016})
\begin{equation} \label{eq_flux_difference}
    \pder[\fnc{F}(\fnc{U})]{x} \approx 2 \sum_{j=1}^N \mat{D}_{ij} \fnc{F}^\star(u_i, u_j) , 
\end{equation}
where $\fnc{F}^\star$ is a consistent and symmetric two--point flux, i.e.\ $\fnc{F}^\star(u,u)=\fnc{F}(u)$ and $\fnc{F}^\star(u_i,u_j) = \fnc{F}^\star(u_j,u_i)$. 
When $\fnc{F}^\star$ is smooth, \eqref{eq_flux_difference} is accurate to the same order of accuracy as $\mat{D}$~\cite{fisher_2013, ranocha_2018}. It can also be recast in a telescopic form~\cite{fisher_2013}, establishing discrete conservation. Consider \eqref{eq_flux_difference} applied to the variable-coefficient linear advection equation \eqref{eq_lin_conv}, where $\fnc{F} = a(x) \fnc{U}$. In the following, we will introduce several choices for $\fnc{F}^\star$.

In \S\ref{sec_continuous_VCLCE} we introduced the $a$-energy as a conserved quantity of \eqref{eq_lin_conv}. Let the discrete analogue of the $a$-energy be defined as $\norm{\bu}_{a \Hnrm}^2 \defn \sum_{i=1}^N \Hnrm_{ii} a_i u_i^2$, where $a_i \defn a(x_i)$. The entropy-stability theory of~\cite{tadmor_2003} then states that the central flux is the unique two-point flux that conserves the $a$-energy,
\begin{equation} \label{eq_central_flux}
\begin{gathered}
\fnc{F}^\star(u_i,u_j) = \mean{ a u }_\text{arith}^{ij} \defn \frac{1}{2} \left( a_i u_i + a_j u_j \right) \\
\Rightarrow \quad 2 \sum_{j=1}^N \mat{D}_{ij} \fnc{F}^\star(u_i, u_j) = \sum_{j=1}^N \mat{D}_{ij} a_j u_j.
\end{gathered}
\end{equation}
Furthermore, entropy-stability theory guarantees that the $a$-energy satisfies a telescoping local conservation law~\cite{tadmor_2003, fisher_2013}. For example, using a second order central difference stencil, one obtains
\[
\begin{gathered}
\der[u_i]{t} + \frac{a_{i+1}u_{i+1}-a_{i-1}u_{i-1}}{2h} = 0 \\
\Rightarrow \quad
\der[(a_i u_i^2)]{t} + \frac{\fnc{G}^\star(u_i,u_{i+1})-\fnc{G}^\star(u_{i-1},u_i)}{h} , \quad  \fnc{G}^\star(u_i,u_j) \defn a_i u_i a_j u_j .
\end{gathered}
\]
Since the two-point flux $\fnc{F}^\star$ is bilinear, conservation of the $a$-energy implies that the semidiscretization Jacobian $\mat{J}$ is linear and skew-adjoint with respect to the $a$-inner product. It follows that the spectrum of $\mat{J}$ is purely imaginary, \ie the semidiscretization exhibits no modal growth, consistent with the continuous differential operator. The addition of linear surface or volume dissipation can only further decrease the real parts of the eigenvalues~\cite{dissipation}, establishing the desired spectral bound $\Re{(\lambda(\mat{J}))}_{\max} \leq 0$. 

Take by contrast some other consistent and symmetric two-point flux, for example the logarithmic flux,
\begin{equation} \label{eq_log_flux}
    \fnc{F}^\star(u_i,u_j) = \mean{ a u }^{ij}_\text{log} \defn \frac{a_j u_j - a_i u_i}{\log(a_j u_j) - \log(a_i u_i)} .
\end{equation}
The logarithmic mean is of particular interest because it forms the basis of entropy-stable flux-differencing schemes for the Euler and Navier--Stokes equations \cite{ranocha_2018, ranocha_2022_preventing_pressure}. For the linear advection problem \eqref{eq_lin_conv}, assuming $u>0$, a flux-differencing semidiscretization with the flux \eqref{eq_log_flux} conserves the nonlinear entropy $\fnc{S}(\bu) \defn \sum_{i=1}^N \Hnrm_{ii} \left( u_i \log\left(a_i u_i\right) - u_i \right)$, a discrete analogue of $\fnc{S}(\fnc{U}) \defn \int_0^L \fnc{U} \log\left(a \fnc{U}\right) - \fnc{U} \ \mathrm{d}x$. Although this may encode important physics, since the entropy is not quadratic and is not defined on a linear space (it requires positivity of $a$ and $u$), it does not result in a spectral bound. An argument goes, therefore, that since the semidiscretization Jacobian can have eigenvalues with positive real parts, perturbations can grow unphysically~\cite{ranocha_2022_preventing_pressure}. In later sections, however, we will show that this modal (or local linear) instability does not contribute to significant perturbation growth. We therefore argue that the logarithmic flux remains equally suitable for this problem as the central flux, provided that $u$ and $a$ are sufficiently bounded away from 0.

Unfortunately, the logarithmic flux is difficult to analyze from a theoretical perspective due to the lack of a compact matrix representation. For this reason, we introduce the geometric flux,
\begin{equation} \label{eq_geom_flux}
    \begin{gathered}
\fnc{F}^\star(u_i,u_j) = \mean{ a u }^{ij}_\text{geom}  \defn \sqrt{ a_i u_i a_j u_j } \\
\Rightarrow \quad 2 \sum_{j=1}^N \mat{D}_{ij} \fnc{F}^\star(u_i, u_j) = 2 \sqrt{a_i u_i} \sum_{j=1}^N \mat{D}_{ij} \sqrt{a_j u_j}.
    \end{gathered}
\end{equation}
The geometric mean has been used as an algebraic approximation for the logarithmic mean in discretizations of the compressible Euler equations~\cite{de_michele_coppola_2023}. For the linear advection equation \eqref{eq_lin_conv}, assuming $u>0$, the flux \eqref{eq_geom_flux} conserves the nonlinear entropy $\fnc{S}(\bu) \defn -2 \sum_{i=1}^N \Hnrm_{ii} \sqrt{ u_i / a_i } $,
a discrete analogue of $\fnc{S}(\fnc{U}) \defn -2 \int_0^L \sqrt{ \fnc{U} / a } \ \mathrm{d}x$. Since this entropy is not quadratic, the same arguments used for the logarithmic flux regarding unphysical perturbation growth can be applied. A well-known useful result (\eg \cite{ALGinequality}) is the relation
\[
\mean{ a u }^{ij}_\text{geom} < \mean{ a u }^{ij}_\text{log} < \mean{ a u }^{ij}_\text{arith} ,
\]
which implies that the logarithmic flux can be expressed as a convex (albeit nonlinear) combination of the arithmetic and geometric fluxes. This will allow us to generalize certain results obtained for the geometric scheme to the logarithmic scheme.

Finally, consider the product flux,
\begin{equation} \label{eq_product_flux}
    \begin{gathered}
\fnc{F}^\star(u_i,u_j) = \mean{ a \cdot u }^{ij}_\text{prod} \defn \frac{1}{2} \left( a_i u_j + a_j u_i \right) \\
\Rightarrow \quad 2 \sum_{j=1}^N \mat{D}_{ij} \fnc{F}^\star(u_i, u_j) = a_i \sum_{j=1}^N \mat{D}_{ij} u_j + u_i \sum_{j=1}^N \mat{D}_{ij} a_j.
    \end{gathered}
\end{equation}
This semidiscretization appears in linearizations of entropy-stable discretizations of the Burgers equation (\eg see Appendix~\ref{app_lin_burgers}), and therefore governs the evolution of perturbations in such schemes~\cite{gassner_stability_2022}. Although the flux is symmetric, consistent, and conservative, \eqref{eq_product_flux} conserves neither the $a$-energy, nor any known convex entropy. Therefore, the semidiscretization Jacobian $\mat{J}$ can have eigenvalues with unphysical positive real parts.

\subsection{Establishing Stability and Convergence for Non-Central Fluxes} \label{sec_stability_proof}

Unlike the central flux, the entropy-stable fluxes \eqref{eq_log_flux} and \eqref{eq_geom_flux} do not immediately yield an $L^2$ energy estimate for the conservative variables. Fortunately, we can appeal to the well-established entropy-stability literature to construct such a bound \cite{Dafermos2010, svard_2015, sbpbook}, which we detail in Appendix \ref{app_entropy_l2_bound}. A key assumption that makes this possible is the numerical enforcement of positivity. In fact, one can easily verify that the fluxes \eqref{eq_log_flux} and \eqref{eq_geom_flux} are ill-defined when $\fnc{U} \leq 0$. The consequences of approaching such a regime are discussed later in \S\ref{sec:nearvacuum}.

To prove stability of the product flux \eqref{eq_product_flux}, we generalize the stability proof of \cite{kopriva_energy_2014}, which directly bounds the $L^2$ energy norm but allows for exponential growth in time according to physical transient growth. Therefore, it is weaker than the $a$-energy bound of the central flux \eqref{eq_central_flux}, yet sharper and more physically interpretable than the bound of \cite{Hairer2016}, which assumes skew-symmetry of $\mat{D}$ and Lipschitz continuity of $a(x)$. For simplicity, we work with periodic SBP operators $\mat{Q} + \mat{Q}^\T = 0$, though the results can be extended to finite domains by choosing appropriate SATs (\eg see \cite{ranocha_variable_advection}). Consider a split-form discretization of \eqref{eq_lin_conv} with an arbitrary combination of the central and product fluxes,
\begin{equation} \label{eq_general_linear_split_form}
    \der[\bu]{t} + \alpha \mat{D} \mat{A} \bu + (1-\alpha)  \left( \mat{A} \mat{D} \bu + \mat{U} \mat{D} \vec{a} \right) = \vec{0} , \quad \mat{U} = \diag{\vec{u}} , \ \mat{A} = \diag{\vec{a}} .
\end{equation}
We refer to this as the central-product split-form scheme. Left-multiplying by $\bu^\T \Hnrm$ and using continuity in time to extract the temporal derivative, we obtain
\begin{equation*}
    \frac{1}{2} \der{t} \norm{\bu}^2_\Hnrm  + (\alpha - \tfrac{1}{2})  \bu^\T \left( \mat{Q} \mat{A} - \mat{A} \mat{Q} \right) \bu + (1-\alpha) \bu^\T  \mat{U} \mat{Q} \vec{a} = 0 .
\end{equation*}
In \cite{kopriva_energy_2014}, $\alpha=\frac{1}{2}$ is chosen to eliminate the second term. The third term is then bounded through
\begin{align*}
    \abs{\bu^\T  \mat{U} \mat{Q} \vec{a} } = \abs{\sum_i \Hnrm_i u_i^2 ( \mat{D} \vec{a} )_i} &\leq \sum_i \Hnrm_i u_i^2 ( \abs{\mat{D} \vec{a}} )_i \\
    &\leq \max_i{(\abs{\mat{D} \vec{a}})_i} \sum_i \Hnrm_i u_i^2 = \norm{ \mat{D} \vec{a} }_\infty \norm{\bu}^2_\Hnrm .
\end{align*}
This choice is unnecessary, however, as the second term is also a consistent discretization of the zeroth-order transient growth, as seen through
\[
\bu^\T \left( \mat{Q} \mat{A} - \mat{A} \mat{Q} \right) \bu \approx \int \fnc{U} \left[ \pder[a(x) \fnc{U}]{x}  -  a(x) \pder[\fnc{U}]{x} \right] \mathrm{d}x = \int \fnc{U}^2 \pder[a(x)]{x}  \mathrm{d}x . 
\]
Therefore, by introducing a product-rule defect term 
\[
\mat{\Theta} \defn \mat{D} \mat{A} - \mat{A} \mat{D} - \diag{\mat{D}\vec{a}} ,
\]
we can combine both terms to achieve a bound through
\begin{align*}
    & \, \abs{(\alpha - \tfrac{1}{2})  \bu^\T \left( \mat{Q} \mat{A} - \mat{A} \mat{Q} \right) \bu + (1-\alpha) \bu^\T  \mat{U} \mat{Q} \vec{a}}\\ 
    = & \, \abs{\bu^\T \Hnrm \left( (\alpha - \tfrac{1}{2}) \left( \mat{D} \mat{A} - \mat{A} \mat{D} \right) + (1-\alpha) \diag{\mat{D} \vec{a}}  \right) \bu} \\
    = & \, \abs{\bu^\T \Hnrm^{1/2} \Hnrm^{1/2} \left( \tfrac{1}{2} \diag{\mat{D} \vec{a}} + (\alpha - \tfrac{1}{2}) \mat{\Theta}  \right) \Hnrm^{-1/2} \Hnrm^{1/2} \bu} \\
    \leq & \, \tfrac{1}{2} \norm{ \diag{\mat{D} \vec{a}} + \left( 2 \alpha - 1 \right) \mat{\Theta} }_\Hnrm \norm{\bu}^2_\Hnrm ,
\end{align*}
where the matrix norm is induced by the vector norm $\norm{\bu}_\Hnrm = \sqrt{\sum_i^N \Hnrm_{ii} u_i^2}$. This in turn allows us to establish boundedness of the solution $\bu$ through
\begin{equation} \label{eq_product_stability}
    \der{t} \norm{\bu}^2_\Hnrm  \leq \gamma \norm{\bu}^2_\Hnrm , \quad \gamma \defn \norm{\diag{\mat{D} \vec{a}} + (2\alpha - 1) \mat{\Theta}}_\Hnrm .
\end{equation}
Grönwall's inequality then results in bounded transient exponential growth $\norm{\bu}^2_\Hnrm \leq \norm{\bu_0}^2_\Hnrm e^{\gamma t}$. When ${\alpha = \tfrac{1}{2}}$, we recover $\gamma = \norm{\mat{D} \vec{a}}_\infty$, consistent with the continuous case where $\gamma = - \min_{x \in \Omega} \partial a / \partial x$ \cite{kopriva_energy_2014}. More generally, $\gamma$ is a mesh-dependent worst-case operator bound that will not necessarily approach the continuous growth coefficient under mesh refinement. Nevertheless, for any nodal restriction $\bu$ of a smooth function $\fnc{U}$, consistency of the defect term $\mat{\Theta}$ implies $\mat{\Theta} \bu = \order{p}$. Therefore, the effective growth rate of physical (\ie non-mesh-dependent) modes is still asymptotically governed by the discrete approximation of $\partial a / \partial x$.

The exponential stability bound \eqref{eq_product_stability} is sufficient to invoke the Lax Equivalence theorem~\cite{lax_equivalence} to establish convergence for the product flux~\eqref{eq_product_flux}. Unfortunately, the nonlinearity of the logarithmic and geometric fluxes means that we must instead appeal to linearization arguments, such as those from \cite{strang}, which bound perturbation growth via estimates on the semidiscrete Jacobian. In Appendix \ref{app_linearizations}, we explicitly linearize both the geometric and logarithmic flux-differencing schemes, and provide a stability bound for the linearized geometric scheme analogous to \eqref{eq_product_stability}. In Appendix \ref{app_sharp_geom_bound}, we further exploit a reformulation of the geometric scheme in terms of a square-root variable to obtain the following sharper perturbation bound,
\begin{equation} \label{eq_geom_pert_bound}
    \norm{\vec{v}(t)}_{\Hnrm}^2 \leq \frac{\max_{i}{(u_i(t))}}{\min_{i}{(u_i(0))}} \norm{\vec{v}(0)}_{\Hnrm}^2 ,
\end{equation}
where pointwise boundedness of the baseflow $\bu$ (\ie the states about which we linearize the geometric scheme to obtain the semidiscretization governing the evolution of perturbations $\bv$) is discussed in Appendix~\ref{app_entropy_l2_bound}. Deriving analogous perturbation bounds for the logarithmic flux is more challenging, as the lack of a compact matrix representation leads to cumbersome algebra. Instead, in the remainder of this work we often use the geometric scheme as a prototypical example of nonlinear flux-differencing schemes, and instead rely on numerical experiments to demonstrate perturbation stability of the logarithmic flux.

\subsection{Questioning the `Central Plus (Anti-)Dissipation' Interpretation} \label{sec_modified_pde}

One viewpoint in the literature \cite{gassner_stability_2022, tadmor_2003} is to interpret general numerical fluxes as a central flux plus some corrective term, often interpreted as dissipation that can be either positive or negative. To illustrate this perspective, we follow \cite{gassner_stability_2022} and rewrite the general central-product split-form scheme~\eqref{eq_general_linear_split_form} for second-order central differencing as a central scheme plus a corrective term,
\begin{gather*}
\der[u_i]{t} + \frac{a_{i+1} u_{i+1} - a_{i-1} u_{i-1}}{2h} \\
= \frac{1-\alpha}{2} h^2 \left[ \frac{a_{i+1} -2 a_i + a_{i-1}}{h^2} \frac{u_{i+1}-u_{i-1}}{2h} + \frac{a_{i+1}-a_{i-1}}{2h} \frac{u_{i+1} -2 u_i + u_{i-1}}{h^2} \right] ,
\end{gather*}
which is a consistent discretization of the modified PDE
\begin{equation} \label{eq_modified_pde_prod}
\pder[\fnc{U}]{t} + \pder[a(x)\fnc{U}]{x} = \pder{x} \left( \nu(x) \pder[\fnc{U}]{x} \right) , \quad \nu(x) \defn h^2 \frac{1-\alpha}{2} \pder[a(x)]{x} .
\end{equation}
Compared to the central scheme, splittings with $\alpha \neq 1$ therefore introduce additional design-order error in a form that resembles anti-dissipation if $\nu < 0$, \ie in regions of decreasing~$a(x)$. This may appear concerning, but it is worthwhile to consider two things. First, this anti-dissipation (or dissipation, depending on the sign of $\partial a / \partial x$) is still consistent with the physical energy production (or decay) induced by the variable coefficient~$a(x)$ in the continuous problem \eqref{eq_lin_conv}. Indeed, the central scheme also exhibits energy growth in regions of decreasing~$a(x)$. The distinction is that the central scheme does it in a manner such that local contributions telescope to the boundaries so that the global energy growth, as measured by the discrete $a$-norm, exactly matches that of the continuous PDE. For more general splittings $\alpha \neq 1$, we proved in \S\ref{sec_stability_proof} that the total growth of $\norm{\bu}_\Hnrm$ is bounded by something mimetic of the physical growth induced by $a(x)$. Therefore, we should not expect this corrective term to introduce catastrophic unphysical energy growth, but rather some addition to the already present physical growth (or decay) such that the global behaviour no longer exactly matches the expected norm. Later sections of this work will further demonstrate that this additional growth is typically easily suppressed, if not negligible. The second important point to consider is that the central scheme also introduces design-order error. Although the global norm $\norm{\bu}_{a\Hnrm}$ of the central scheme is exactly balanced, local errors remain design-order, just as any energy-conserving discretization of the constant-coefficient linear advection or Burgers' equations will also only be design-order, despite exactly conserving energy. Therefore, we should not immediately assume that the $\alpha=1$ discretization is more accurate than general splittings with $\alpha \neq 1$.

We now repeat this analysis with the geometric flux (and by extension the logarithmic flux), but will find that the results are not as clear. A second-order geometric flux-differencing scheme can be expressed as
\begin{equation} \label{eq_geom_scheme_2nd_order}
\begin{gathered}
\der[u_i]{t} + 2 \sqrt{a_i u_i} \frac{\sqrt{a_{i+1}u_{i+1}}-\sqrt{a_{i-1}u_{i-1}}}{2h}=0 \\
\Leftrightarrow \quad \der[u_i]{t} + \frac{a_{i+1} u_{i+1} - a_{i-1} u_{i-1}}{2h} = h^2 \left[ \frac{w_{i+1}-w_{i-1}}{2h} \frac{w_{i+1}-2w_i+w_{i-1}}{h^2} \right] ,
\end{gathered}
\end{equation}
where we have introduced a square-root variable $w_i \defn \sqrt{a_i u_i}$. Letting the continuous counterpart be $\fnc{W} \defn \sqrt{a(x) \fnc{U}}$ and recognizing that 
\[
\pder[\fnc{W}]{x} \pder[^2\fnc{W}]{x^2} = \frac{1}{2} \pder{x} \left( \pder[\fnc{W}]{x} \right)^2 = \frac{1}{8} \pder{x} \left( \frac{1}{a(x) \fnc{U}} \left( \pder[a(x) \fnc{U}]{x} \right)^2 \right) ,
\]
we see that \eqref{eq_geom_scheme_2nd_order} is a consistent discretization of the modified PDE
\begin{equation*}
\pder[\fnc{U}]{t} + \pder[a(x)\fnc{U}]{x} = \pder{x} \left( \frac{\varepsilon(x)}{\fnc{U}} \left( \pder[a(x) \fnc{U}]{x} \right)^2 \right) , \quad \varepsilon(x) \defn \frac{h^2}{8 a(x)} > 0.
\end{equation*}
In contrast to the central-product split-form, the additional term resulting from the correction to the central scheme is not easily interpreted as diffusive. Instead, it represents the conservative transport of some nonlinear quantity in~$\fnc{U}$. Performing a standard energy analysis results in an energy balance of
\[
\der{t} \norm{\fnc{U}}^2 = - \int_0^L \pder[a(x)]{x} \fnc{U}^2 \, \mathrm{d}x - \int_0^L \frac{\varepsilon(x)}{\fnc{U}} \pder[\fnc{U}]{x} \left( \pder[a(x)\fnc{U}]{x} \right)^2 \mathrm{d} x ,
\]
revealing that---similar to anti-dissipation---the correction can increase energy in regions where $\fnc{U}$ is decreasing (and vice-versa). However, this interpretation may be misleading, since not only would such a dissipation coefficient depend on the solution itself, but periodicity of $\fnc{U}$ ensures that there are always regions of both energy growth and decay. Another revealing perspective, which we detail in Appendix \ref{app_reformulation_geom}, is obtained by equivalently expressing the geometric scheme \eqref{eq_geom_scheme_2nd_order} as a central scheme in $\fnc{W}$,
\[
\der[u_i]{t} + 2 \sqrt{a_i u_i} \frac{\sqrt{a_{i+1}u_{i+1}}-\sqrt{a_{i-1}u_{i-1}}}{2h}=0 
\quad \Leftrightarrow \quad 
\der[w_i]{t} + a_i \frac{w_{i+1}-w_{i-1}}{2h}=0 ,
\]
which in \cite{manzanero} was shown to introduce no excess energy growth for the non-conservative form of the variable-coefficient linear advection equation. That is, just as the central scheme exhibits no solution growth as measured by conservation of the $a$-energy $\norm{\bu}^2_{a}$, the geometric scheme \eqref{eq_geom_scheme_2nd_order} exhibits no solution growth as measured by conservation of the $\tfrac{1}{a}$-energy $\norm{\vec{w}}_{1/a}^2$.

We have already proven via the bound in \eqref{eq_geom_pert_bound} that perturbations of the geometric flux differencing scheme cannnot grow in time, provided that the baseflow does not grow in time. The `central plus (anti-)dissipation' perspective therefore once again appears inconclusive when applied to linearizations of the geometric scheme. Computed explicitly in Appendix \ref{app_linearizations}, the linearized geometric scheme can be written in the central plus correction form as
\begin{equation} \label{eq_geom_scheme_lin_2nd_order}
\begin{gathered}
    \der[v_i]{t} + \frac{a_{i+1} v_{i+1} - a_{i-1} v_{i-1}}{2h} \\
    = \frac{h^2}{2} \left[ \frac{w_{i+1} -2 w_i + w_{i-1}}{h^2} \frac{z_{i+1} - z_{i-1} }{2h} + \frac{w_{i+1}-w_{i-1}}{2h} \frac{z_{i+1} -2 z_i +  z_{i-1}}{h^2} \right] ,
\end{gathered}
\end{equation}
where $z_i \defn w_i v_i / u_i = \sqrt{a_i / u_i} v_i$. By introducing the variable $\fnc{Z} \defn \sqrt{a(x) / \fnc{U}} \, \fnc{V}$ and performing tedious algebra (which we omit for brevity), the term in square brackets of the correction can be interpreted as
\begin{align*}
& \pder[^2\fnc{W}]{x^2} \pder[\fnc{Z}]{x} + \pder[\fnc{W}]{x} \pder[^2\fnc{Z}]{x^2} = \pder{x} \left( \pder[\fnc{W}]{x} \pder[\fnc{Z}]{x} \right) \\
= & \, \pder{x} \left( \left( \frac{1}{2} \frac{1}{ \sqrt{a(x) \fnc{U}}} \pder[a(x) \fnc{U}]{x} \right) \left( \sqrt{\frac{a(x)}{\fnc{U}}} \pder[\fnc{V}]{x} + \frac{\fnc{V}}{2} \sqrt{\frac{\fnc{U}}{a(x)}} \pder[a(x)/\fnc{U}]{x} \right) \right) \\
= & \, \frac{1}{2} \pder{x} \left( \frac{1}{\fnc{U}} \pder[a(x) \fnc{U}]{x} \pder[\fnc{V}]{x} + \frac{1}{2 a(x)} \left( \left( \pder[ a(x)]{x} \right)^2 - \frac{a(x)^2}{\fnc{U}^2} \left( \pder[\fnc{U}]{x} \right)^2 \right) \fnc{V} \right),
\end{align*}
demonstrating that \eqref{eq_geom_scheme_lin_2nd_order} is a consistent discretization of the modified PDE
\begin{equation} \label{eq_modified_pde_geom}
\begin{gathered}
\pder[\fnc{V}]{t} + \pder{x} \left[ \left( a(x) - \beta(x) \right)  \fnc{V} \right] = \pder{x} \left( \nu(x) \pder[\fnc{V}]{x} \right) ,\\
\beta(x) \defn \frac{h^2}{8} a(x) \left( \frac{1}{a(x)^2} \left(\pder[a(x)]{x}\right)^2 - \frac{1}{\fnc{U}^2} \left(\pder[\fnc{U}]{x}\right)^2 \right) , \quad \nu(x) \defn \frac{h^2}{4 \, \fnc{U}} \pder[a(x) \fnc{U}]{x}.
\end{gathered}
\end{equation}
While $\nu<0$ again identifies regions in which the modified equation contains an anti-diffusive contribution, similar to the analysis done for the baseflow, the coefficient $\nu$ depends nonlinearly on the baseflow and its derivative. Therefore, not only is the predictive power of this analysis limited in the context of a changing baseflow, but periodicity ensures that it will always induce regions of both perturbation growth and decay. Moreover, the modified equation contains an additional transport term $\beta(x)$, so the correction can no longer be interpreted as a purely diffusive perturbation of the central scheme. Consequently, unlike in the product-flux case, the `central plus (anti-)dissipation' perspective does not by itself provide a reliable prediction of perturbation growth or decay. The square-root variable reformulation of Appendix~\ref{app_sharp_geom_bound} and the resulting perturbation bound~\eqref{eq_geom_pert_bound} remain more informative.

\subsection{Frozen-Coefficient Analysis} \label{sec_frozen_coeff}

To gain insights into the nature of local linear instabilities---by which we mean the eigenmodes corresponding to eigenvalues with positive real parts, which have the potential to trigger the modal growth of perturbations in comparison to the central scheme---we analyze the modified PDE \eqref{eq_modified_pde_geom} which we obtained from the geometric scheme,
\[
\pder[\fnc{V}]{t} + \pder{x} \left[ \left( a(x) - \beta(x) \right)  \fnc{V} \right] = \pder{x} \left( \nu(x) \pder[\fnc{V}]{x} \right) ,
\]
but now leave the coefficients $\beta(x)$ and $\nu(x)$ unspecified. Note that the modified PDE \eqref{eq_modified_pde_prod} obtained from the product scheme is recovered by setting $\beta(x) = 0$. To analyze individual Fourier modes, we must make simplifications. We first consider a sufficiently small local region of space so that $a(x)$ and $\fnc{U}(x)$ are approximately linear, \ie
\[
a(x) \approx a_0 + a_1 \delta , \quad \fnc{U}(x) \approx \fnc{U}_0 + \fnc{U}_1 \delta , \quad \delta \defn x - x_0 .
\]
For the central-product split-form, we immediately find $\nu \approx h^2 \tfrac{1-\alpha}{2} a_1$. For convenience we define $\nu_0 \defn \tfrac{1-\alpha}{2} a_1$. The geometric scheme requires more effort, but Taylor expansions of the form $\beta(x) \approx h^2 (\beta_0 + \beta_1 \delta)$ and $\nu(x) \approx h^2 (\nu_0 + \nu_1 \delta)$ result in the following locally accurate modified PDE,
\[
\pder[\fnc{V}]{t} + (a_0 - h^2 \beta_0 - h^2 \nu_1) \pder[\fnc{V}]{x} + (a_1 - h^2 \beta_1) \delta \pder[\fnc{V}]{x} + (a_1 - h^2 \beta_1) \fnc{V}  = h^2 (\nu_0 + \nu_1 \delta) \pder[^2\fnc{V}]{x^2}  .
\]
The details of this derivation are provided in Appendix \ref{app_frozen_coeff}. 

We now wish to study Fourier mode solutions of the form $\fnc{V} = \hat{\fnc{V}} e^{i \kappa x + \lambda t}$ with spatial wavenumber $\kappa \in \mathbb{R}$ and temporal eigenvalue $\lambda \in \mathbb{C}$. To do this, we must make the additional assumption of a high wavenumber and localized $\fnc{V}$. That is, we consider modes such that the wavelength is on the same order as the local region we are considering, $\kappa \delta \sim 1$. We also assume that the local region is small enough such that the total variation in the baseflow is small in comparison to its total value, \ie $\abs{ a_1 \delta} \ll \abs{a_0}$ and $\abs{\fnc{U}_1 \delta} \ll \abs{\fnc{U}_0}$. With these assumptions, the modified PDE simplifies to
\begin{equation} \label{eq_modified_pde_frozen}
\pder[\fnc{V}]{t} + (a_0 - h^2 \beta_0 - h^2 \nu_1) \pder[\fnc{V}]{x} + (a_1 - h^2 \beta_1) \fnc{V}  = h^2 \nu_0 \pder[^2\fnc{V}]{x^2} .
\end{equation}
Once again, the justifications for dropping negligible terms are detailed in Appendix \ref{app_frozen_coeff}. Substituting $\fnc{V} = \hat{\fnc{V}} e^{i \kappa x + \lambda t}$ into this constant-coefficient linear equation yields
\[
\lambda + (a_0 - h^2 \beta_0 - h^2 \nu_1) i \kappa + a_1 - h^2 \beta_1 + h^2 \nu_0  \kappa^2 = 0 .
\]
Taking the imaginary and real parts, we find the growth rate and angular frequency
\[
\Re{(\lambda)} = - a_1 + h^2 \beta_1 - h^2 \kappa^2 \nu_0 , \quad 
\Im{(\lambda)} = - a_0 \kappa + h^2 \kappa \beta_0 + h^2 \kappa \nu_1 .
\]

Ignoring the $\order{2}$ corrections, a local Fourier mode exhibits growth in regions of negative baseflow ${a_1 < 0}$, consistent with the local growth induced by the variable coefficient of the continuous problem \eqref{eq_lin_conv}. It also oscillates at a frequency proportional to the wavenumber, consistent with convection. With regard to the corrections, which stem from the design-order correction term of the `central plus (anti-)dissipation' decomposition, we consider the following two important scenarios:
\begin{enumerate}
    \item If $\fnc{V}$ is a physical and well-resolved mode, then $\kappa$ is independent of the mesh spacing $h$. As a result, when $h \rightarrow 0$, $h^2 \kappa \rightarrow 0$, and $h^2 \kappa^2 \rightarrow 0$. Since $\beta_0$, $\beta_1$, $\nu_0$, and $\nu_1$ are all also independent of $h$, $\Re{(\lambda)} \rightarrow - a_1$ and $\Im{(\lambda)} \rightarrow - a_0 \kappa$, \ie the effect of the correction vanishes and local growth becomes consistent with both the central scheme and the continuous problem.
    \item If $\fnc{V}$ is an unphysical or under-resolved mode (\eg near the Nyquist frequency), then ${\kappa = \order{-1}}$, so as $h \rightarrow0$, $h^2 \kappa \rightarrow 0$ but $h^2 \kappa^2 \rightarrow C > 0$. Consequently, $\Re{(\lambda)} \rightarrow - a_1 - C \nu_0$ (where $\nu_0 \defn \tfrac{1-\alpha}{2} a_1$ for the central-product split-form, and $\nu_0 \defn \tfrac{1}{4} (a_1 + \tfrac{a_0}{\fnc{U}_0} \fnc{U}_1)$ for the geometric scheme) and $\Im{(\lambda)} \rightarrow - a_0 \kappa$, \ie there is additional growth that is inconsistent with the central scheme, and by extension, a globally balanced $a$-energy norm.
\end{enumerate}

From this we draw several important conclusions regarding the nature of local linear instabilities. First, we can expect the most significant unphysical perturbation growth to come from under-resolved or highly oscillatory spurious modes. Numerical dissipation should therefore be effective in suppressing the majority of the local linear instabilities. Admittedly, modes associated with the mid-range may still excite some unphysical growth, so an assessment of the worst case scenario practical impact of local linear instabilities remains important. Second, although high-frequency modes excite additional local perturbation growth in regions of negative baseflow gradient, they also excite additional local perturbation \emph{damping} in regions of \emph{positive} baseflow gradient. Therefore, the strongest local linear instabilities will be associated with highly oscillatory but \emph{localized} modes in regions of negative baseflow gradients, such as those associated with the modified boundary stencils of SBP operators. Third, the perturbation growth should be small in comparison to the oscillatory component of the mode, since by assumption, $\abs{a_1} \ll a_0 / \abs{\delta} \sim a_0 \abs{\kappa}$, so $\abs{\Re{(\lambda)}} \ll \abs{\Im{(\lambda)}}$. Therefore, in a nonnormal or time-varying setting, even small phase differences between rapidly oscillating nonorthogonal eigenmodes can produce significant constructive or destructive interference, leading to transient growth or decay that can dominate the modal growth predicted by the real parts of the eigenvalues. In later sections, we will see that this effect is important for both time-varying problems and for nonlinear discretizations of time-invariant problems.

\subsection{Eigenspectra of the Central-Product Split-Form with a Spatially-Varying Coefficient} \label{sec:spectra}

\begin{figure}[t]
    \centering
    \captionsetup[subfigure]{justification=raggedright,singlelinecheck=false}

    \begin{subfigure}[t]{0.253\textwidth}
        \includegraphics[height=3.12cm,
        trim={3mm 3mm 28mm 3mm},clip,left]{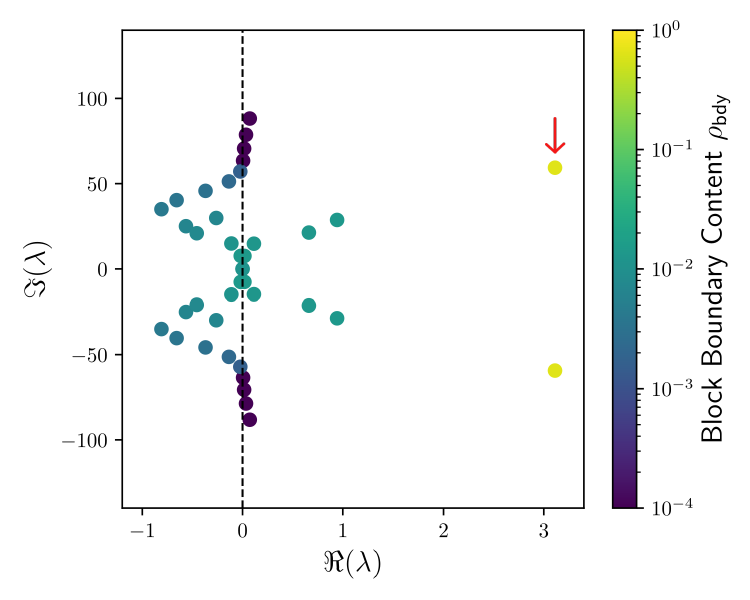}
        \captionsetup{margin={6mm,0mm}}
        \caption{CSBP $p=1$, \\ non-dissipative.\\ $\Re{(\lambda)}_{\max} = 3.1$}
        \label{fig:eigvals1}
    \end{subfigure}
    \hfill
    \begin{subfigure}[t]{0.208\textwidth}
        \includegraphics[height=3.12cm,
        trim={19mm 3mm 28mm 3mm},clip,center]{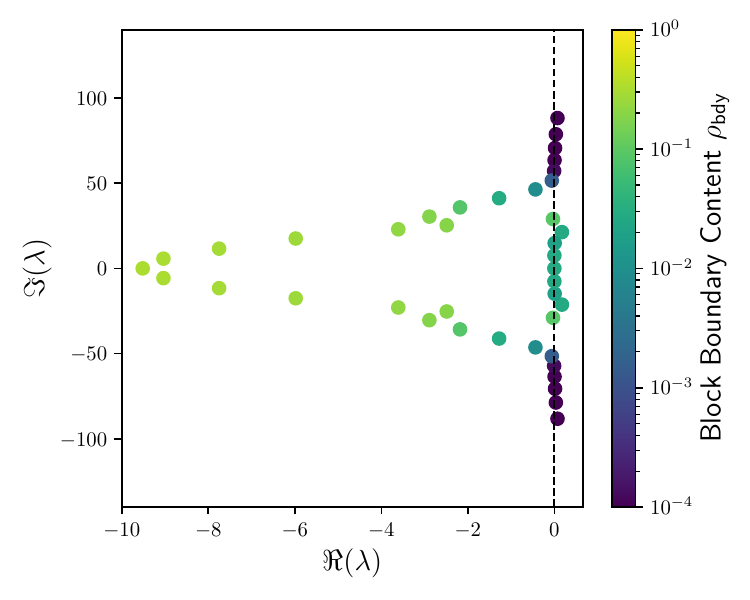}
        \captionsetup{width=1.2\linewidth}
        \caption{CSBP $p=1$, \\ dissipative SAT.\\ $\Re{(\lambda)}_{\max} = 0.17$}
        \label{fig:eigvals2}
    \end{subfigure}
    \hfill
        \begin{subfigure}[t]{0.208\textwidth}
        \includegraphics[height=3.12cm,
        trim={19mm 3mm 24mm 3mm},clip,center]{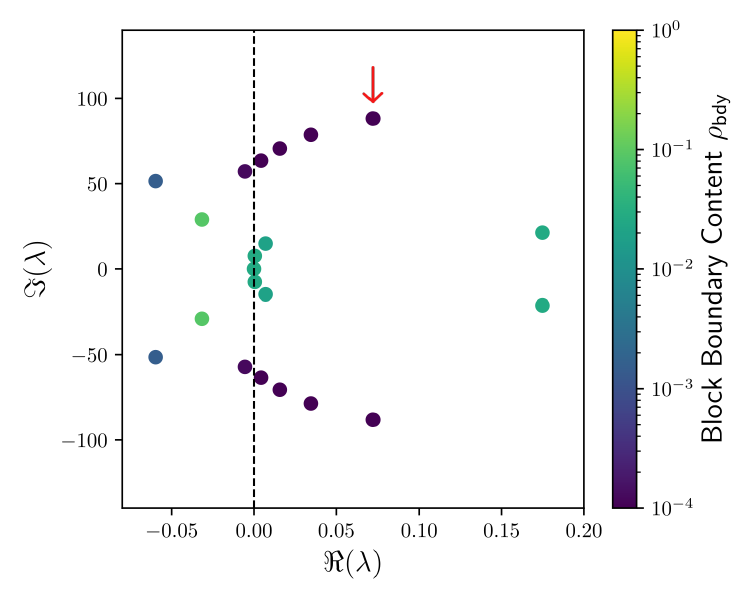}
        \caption{Zoom of (b)}
        \label{fig:eigvals3}
    \end{subfigure}
    \hfill
        \begin{subfigure}[t]{0.272\textwidth}
        \includegraphics[height=3.12cm,
        trim={19mm 3mm 3mm 3mm},clip,right]{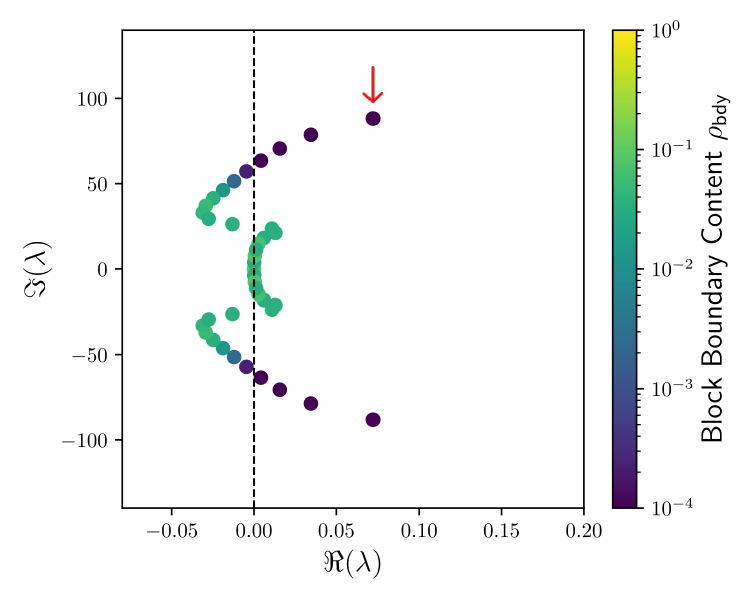}
        \captionsetup{margin={-1mm,0mm}}
        \caption{Circulant 2\textsuperscript{nd}-order,\\ non-dissipative.\\ $\Re{(\lambda)}_{\max} = 0.072$}
        \label{fig:eigvals4}
    \end{subfigure}

    \vspace{0.5em}
    
    \begin{subfigure}[t]{0.253\textwidth}
        \includegraphics[height=3.12cm,
        trim={3mm 3mm 28mm 3mm},clip,left]{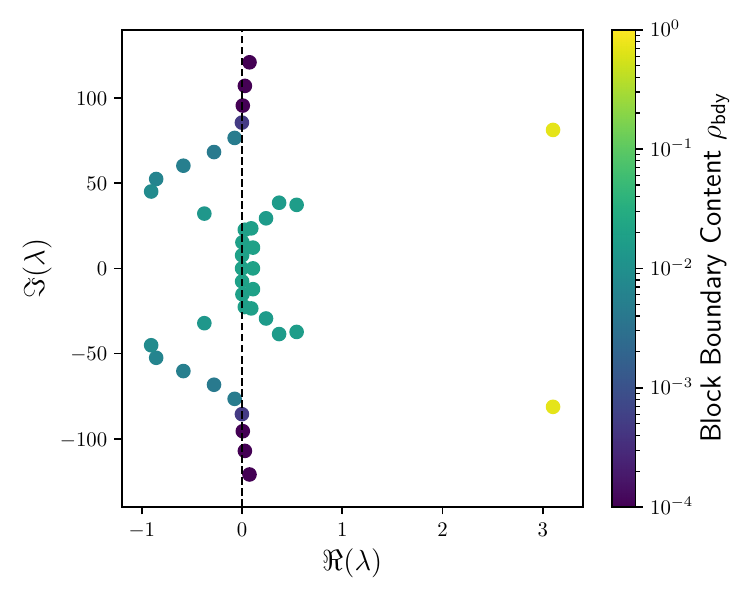}
        \captionsetup{margin={6mm,0mm}}
        \caption{CSBP $p=2$,\\ non-dissipative.\\ $\Re{(\lambda)}_{\max} = 3.1$}
        \label{fig:eigvals5}
    \end{subfigure}
    \hfill
    \begin{subfigure}[t]{0.208\textwidth}
        \includegraphics[height=3.12cm,
        trim={19mm 3mm 28mm 3mm},clip,center]{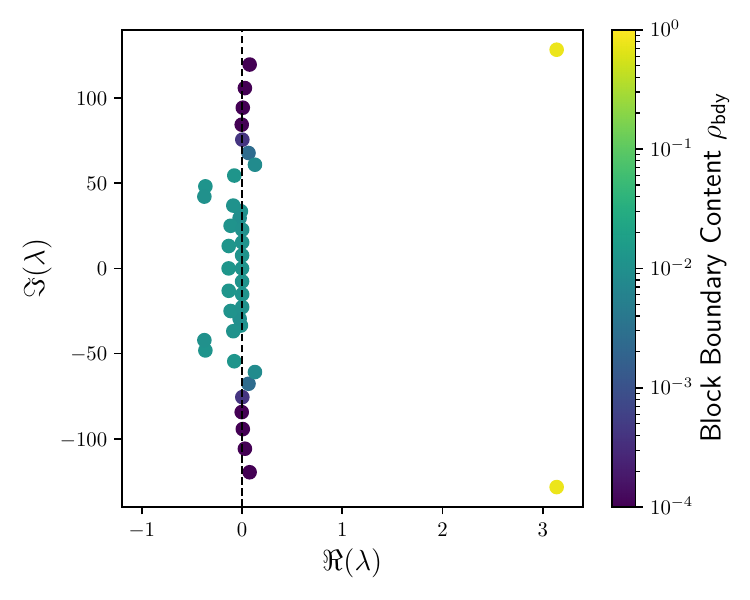}
        \captionsetup{width=1.1\linewidth, margin={-1mm,0mm}}
        \caption{Mattsson $p=2$,\\ non-dissipative.\\ $\Re{(\lambda)}_{\max} = 3.1$}
        \label{fig:eigvals6}
    \end{subfigure}
    \hfill
    \begin{subfigure}[t]{0.208\textwidth}
        \includegraphics[height=3.12cm,
        trim={19mm 3mm 24mm 3mm},clip,center]{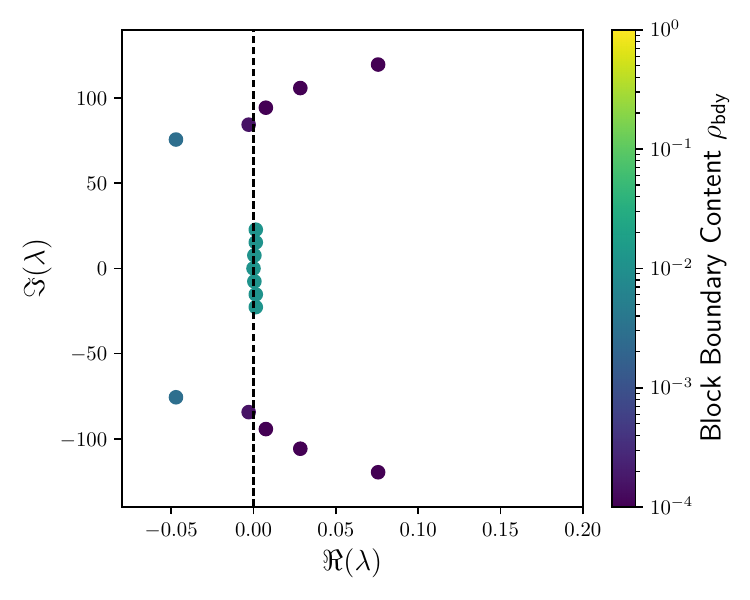}
        \captionsetup{margin={-1mm,0mm}}
        \caption{Mattsson $p=2$,\\ dissipative SAT (zoom).\\ $\Re{(\lambda)}_{\max} = 0.076$}
        \label{fig:eigvals7}
    \end{subfigure}
    \hfill
    \begin{subfigure}[t]{0.272\textwidth}
        \includegraphics[height=3.12cm,
        trim={19mm 3mm 3mm 3mm},clip,right]{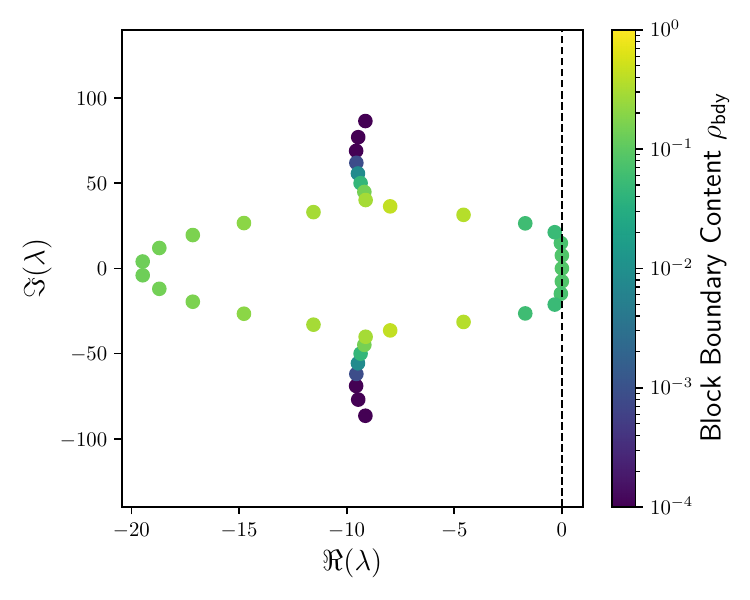}
        \captionsetup{margin={-1mm,0mm}}
        \caption{Circulant 2\textsuperscript{nd}-order,\\ dissipative $s=2$, $\varepsilon=\tfrac{1}{40}$.\\ $\Re{(\lambda)}_{\max} = 0$}
        \label{fig:eigvals8}
    \end{subfigure}

    \caption{Eigenspectra of the product scheme, \ie \eqref{eq_general_linear_split_form} with $\alpha=0$, and one block of $N=40$ nodes. The circulant operators use 39 nodes to match the interior stencil of $p=1$ CSBP operators. Interface dissipation is applied using upwind SATs, and volume dissipation is applied according to \cite{dissipation}. Red arrows indicate the modes analyzed in Figure \ref{fig:eigvecs}. Note the different scales in the $x$-axis.}
    \label{fig:eigvals}
\end{figure}

To verify the conclusions of the frozen-coefficient analysis, we numerically examine the eigenspectra of the semidiscretization Jacobians $\mat{J}$. We focus on the product scheme with a time-constant coefficient, \ie \eqref{eq_general_linear_split_form} with $\alpha = 0$, since the geometric scheme is nonlinear, introducing additional complexity. Therefore, the results discussed here for solution growth can equivalently be interpreted as perturbation growth for entropy-stable and split-form discretizations of the Burgers equation with a frozen baseflow~\cite{gassner_stability_2022} (see Appendix \ref{app_lin_burgers}). To ensure there are no fortuitous symmetry cancellations, we use a skewed sinusoidal variable coefficient that mimics a steepening wave,
\begin{equation} \label{eq_skewed_sin}
    a(x) = 1.5 +
\sum_{k=1}^{n}
\frac{\binom{2n}{n-k}}{\binom{2n}{n}\,k}
\sin\!\left( k \left(
\frac{2\pi (x-x_{\min})}{x_{\max} - x_{\min}} + 4
\right) \right), 
\end{equation}
with $n=5$, and  $x \in [ x_{\min} , x_{\max} ] \defn [0,1]$. In Figure~\ref{fig:eigvals} we plot the eigenvalues of $\mat{J}$, coloured by a metric quantifying the localization of an eigenvector to the block boundaries. Let $\vec{\phi}$ be an eigenvector of $\mat{J}$. We isolate the boundary content of $\vec{\phi}$ by setting all interior nodes to zero, then assess the contribution of this $\vec{\phi}_\text{bdy}$ to the total norm of $\vec{\phi}$ through
\[
\rho_\text{bdy}(\vec{\phi}) \defn \frac{\vec{\phi}_\text{bdy}^* \Hnrm \vec{\phi}_\text{bdy}}{\vec{\phi}^* \Hnrm \vec{\phi}} , \quad \left( \vec{\phi}_\text{bdy} \right)_i \defn \left\lbrace \begin{array}{ll}
   \phi_i  & \text{if node is at a block boundary,} \\
   0  & \text{otherwise.}
\end{array} \right.
\]
Figures \ref{fig:eigvals1}, \ref{fig:eigvals5}, and \ref{fig:eigvals6} include no dissipation. We can clearly identify a dominant unstable eigenmode associated with the block boundaries, for which the positive real part of the eigenvalue does not move as the mesh is refined. Such modes can be created by the modified boundary stencils of finite-difference SBP operators~\cite{beam_warming_1993}. When interface dissipation is applied in figure \ref{fig:eigvals2}, \ref{fig:eigvals3}, or \ref{fig:eigvals7}, however, this unstable mode is effectively suppressed. Circulant operators, shown in Figure \ref{fig:eigvals4}, have no block boundaries and consequently exhibit much weaker instabilities. Furthermore, the real parts converge to zero with mesh refinement, and on any given mesh, spectral stability can be enforced by introducing a small amount of volume dissipation as shown in Figure \ref{fig:eigvals8} (\eg using the framework of \cite{dissipation}). Mattsson `accurate' operators \cite{Mattsson2014} naturally suppress high-frequency spurious modes \cite{dissipation}, so with the exception of the dominant boundary-localized unstable mode, they exhibit smaller instabilities as compared to Classical SBP (CSBP) operators with uniform nodal distributions (compare Figures \ref{fig:eigvals5} and \ref{fig:eigvals6}). The addition of interface dissipation again suppresses the unstable boundary modes.

\begin{figure}[t]
    \centering
    \captionsetup[subfigure]{justification=centering,singlelinecheck=false}

    \begin{subfigure}[t]{0.32\textwidth}
        \centering
        \includegraphics[width=\linewidth, trim={10 12 10 10}, clip]{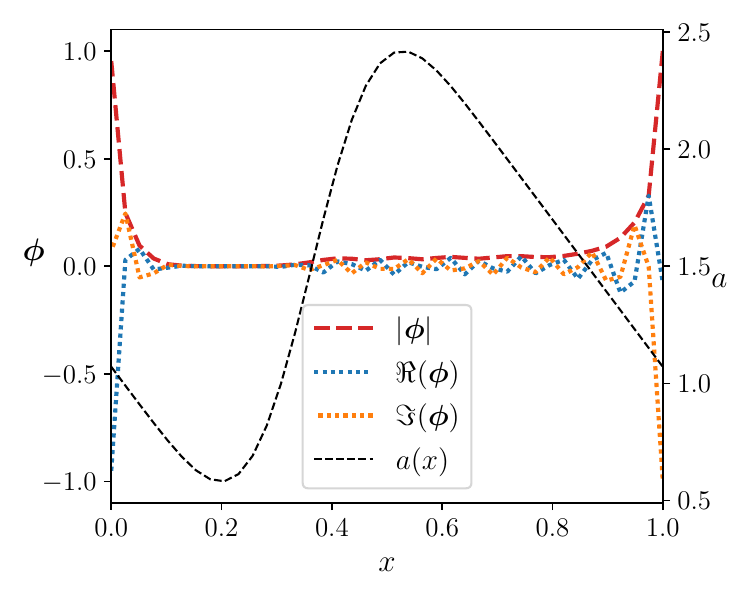}
    \end{subfigure}
    \hfill
    \begin{subfigure}[t]{0.32\textwidth}
        \centering
        \includegraphics[width=\linewidth, trim={10 12 10 10}, clip]{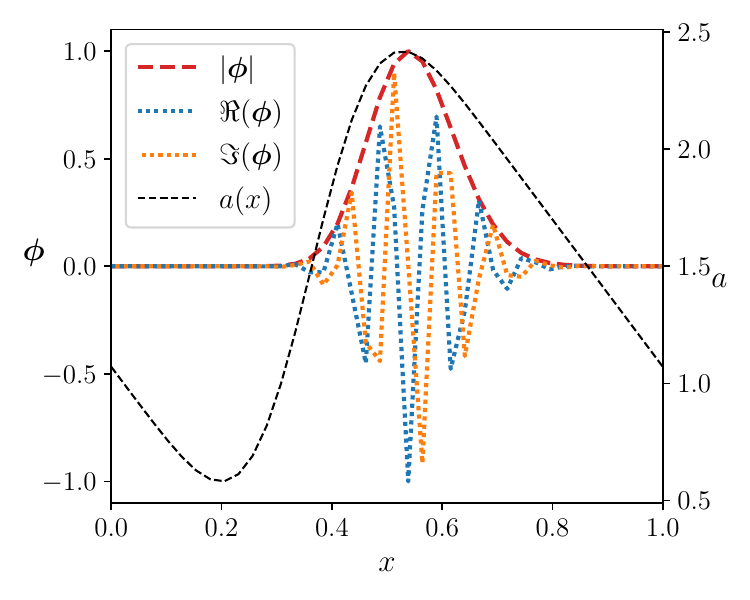}
    \end{subfigure}
    \hfill
    \begin{subfigure}[t]{0.32\textwidth}
        \centering
        \includegraphics[width=\linewidth, trim={10 12 10 10}, clip]{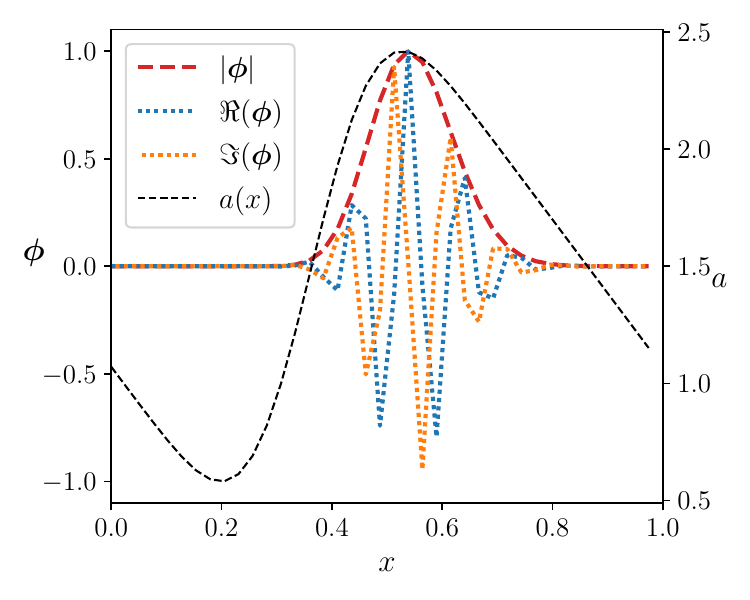}
    \end{subfigure}

    \vspace{0.5em}

    \begin{subfigure}[t]{0.32\textwidth}
        \centering
        \includegraphics[width=\linewidth, trim={10 12 10 10}, clip]{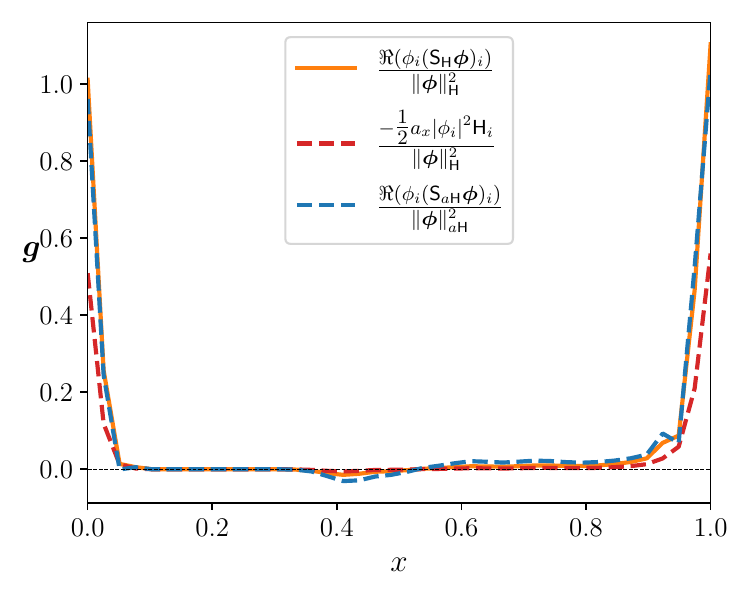}
        \caption{CSBP $p=1$, \\ non-dissipative.}
        \label{fig:eigvecs1}
    \end{subfigure}
    \hfill
    \begin{subfigure}[t]{0.32\textwidth}
        \centering
        \includegraphics[width=\linewidth, trim={10 12 10 10}, clip]{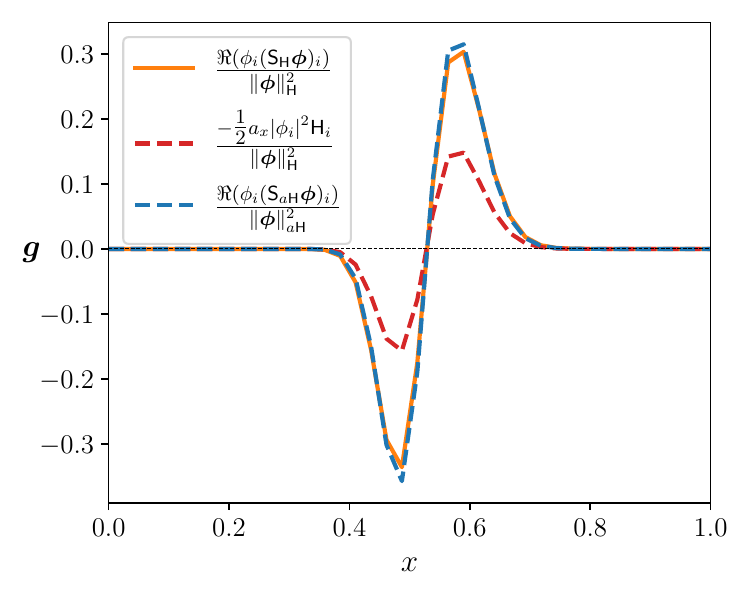}
        \caption{CSBP $p=1$, \\ dissipative SAT.}
        \label{fig:eigvecs2}
    \end{subfigure}
    \hfill
    \begin{subfigure}[t]{0.32\textwidth}
        \centering
        \includegraphics[width=\linewidth, trim={10 12 10 10}, clip]{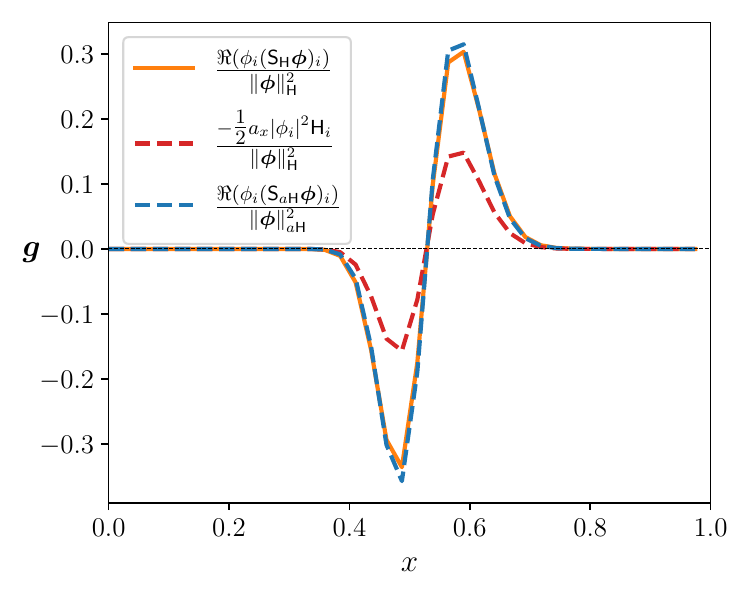}
        \caption{Circulant 2\textsuperscript{nd}-order, \\ non-dissipative.}
        \label{fig:eigvecs3}
    \end{subfigure}

    \caption{Eigenvectors (top row) and contributions to $\Re{(\lambda)}_{\max}$, or local Rayleigh quotients (bottom row), of the product scheme unstable modes indicated by red arrows in Figure~\ref{fig:eigvals}. Note that (b) plots the second-most unstable mode, instead of the first.}
    \label{fig:eigvecs}
\end{figure}

Alongside the variable coefficient \eqref{eq_skewed_sin}, in Figure \ref{fig:eigvecs} we plot the eigenvectors associated with three unstable eigenvalues identified by red arrows in Figure \ref{fig:eigvals} (\ie the dominant mode of CSBP, with and without interface dissipation, and non-dissipative circulant operators), as well as their local contributions to the real parts of the eigenvalues. Since $\norm{\vec{\phi}}_\Hnrm^2 \defn \vec{\phi}^* \Hnrm \vec{\phi}$, we can say that
\begin{align*}
 \quad \der{t} \norm{\vec{\phi}}_\Hnrm^2 &= \der[\vec{\phi}^*]{t} \Hnrm \vec{\phi} + \vec{\phi}^* \Hnrm \der[\vec{\phi}]{t} \\
 &= ( \mat{J} \vec{\phi})^* \Hnrm \vec{\phi} + \vec{\phi}^* \Hnrm \mat{J} \vec{\phi} = (\lambda^* + \lambda) \vec{\phi}^* \Hnrm \vec{\phi} = 2 \Re{(\lambda)} \norm{\vec{\phi}}_\Hnrm^2 ,
\end{align*}
which allows us to define a local contribution vector $\vec{g}$ as
\begin{equation} \label{eq_local_growth_vec}
g_i \defn \frac{\Re{(\phi_i^* (\mat{S}_\Hnrm \vec{\phi})_i)}}{\norm{\vec{\phi}}_\Hnrm^2} , \quad \mat{S}_\Hnrm \defn \tfrac{1}{2} \left( \mat{J}^\T \Hnrm + \Hnrm \mat{J} \right) \quad \text{such that} \quad \sum_i g_i = \Re{(\lambda)} .
\end{equation}
This can be interpreted as measuring the local energy growth of $\vec{\phi}$ determined by the PDE, seen through
\[
\int \fnc{G}(\fnc{V}) \mathrm{d}x = \frac{\tfrac{1}{2} \tfrac{\mathrm{d}}{\mathrm{d}t} \int \fnc{V}^2 \mathrm{d} x}{\int \fnc{V}^2 \mathrm{d}x} = \frac{\int \fnc{V} \fnc{V}_t \mathrm{d} x}{\int \fnc{V}^2 \mathrm{d}x}
= \frac{- \int \fnc{V} (a \fnc{V})_x \mathrm{d} x}{\int \fnc{V}^2 \mathrm{d}x} = \frac{- \tfrac{1}{2} \int a_x \fnc{V}^2 \mathrm{d} x}{\int \fnc{V}^2 \mathrm{d}x} .
\]
Indeed, $\vec{g}$ is precisely the Rayleigh quotient density of the symmetric operator $\mat{S}_\Hnrm$, which measures the local energy growth of $\norm{\vec{\phi}}_\Hnrm^2$. By replacing $\Hnrm$ with $\mat{A} \Hnrm$ in \eqref{eq_local_growth_vec}, we can also quantify the local $a$-norm energy growth, $\vec{g}_a$, which should be zero according to the continuous PDE \cite{manzanero}.

Figure \ref{fig:eigvecs1} confirms the boundary localization of the dominant unstable eigenmode. This mode is nearly identical between splittings $\alpha \in [0,1]$, however as $\alpha \rightarrow 1$, the mode is slightly rotated by $\mat{J}(\alpha)$ such that the real part of its eigenvalue (which is small in comparison to the imaginary component) becomes zero. For $\alpha \neq 0$, since this mode happens to exist in a region of negative coefficient gradient, it excites additional local energy growth, which we ultimately see as a linear instability. By contrast, in Figure \ref{fig:eigvecs3}, the dominant unstable mode of the circulant operator is a result of only interior dynamics, and so is not localized. Therefore, it exists in regions of both positive and negative coefficient gradient, resulting in more balanced total energy growth. To emphasize the importance of damping localized eigenmodes, we draw attention to the second-largest unstable mode of the CSBP discretization with upwind SATs in Figure \ref{fig:eigvecs2} (the largest unstable mode is still interface-localized---it remains unstable because the interface dissipation is only applied at the boundary nodes, which is insufficient to fully compensate for excitations at the next inner-most nodes). This mode is nearly identical to the circulant operator's unstable mode, consistent with \cite{beam_warming_1993}, which shows that the eigenspectrum of a linear operator with a repeating interior stencil but modified boundary stencils (\eg as found in finite-difference SBP discretizations), can be decoupled into the spectrum of the circulant operator with the same interior stencil, plus some boundary-dependent modes resulting from the modified boundary stencils.

Although we have only discussed finite-difference SBP operators so far, the same conclusions hold true for spectral-element operators such as LG and LGL. For non-dissipative schemes, the dominant unstable eigenmodes are localized to element boundaries, though since spectral-element schemes typically employ more elements than finite-difference schemes, we also see more of these boundary localized modes (here we have intentionally simplified our examples to use a single block). Interface dissipation likewise plays a stabilizing role. In fact, because the volume stencils are typically small in comparison to the total the number of interfaces, we have found that interface dissipation is almost always sufficient to achieve spectral stability. This is consistent with \cite{gassner_stability_2022}, where eigenvalues with positive real parts were observed in only the most extreme under-resolved scenarios (\eg degree $p=15$ LGL operators with an element-wise constant solution and $\fnc{O}(1)$ jumps between elements) after interface dissipation was included.

\begin{figure}[t]
    \centering

    \begin{subfigure}[t]{0.32\textwidth}
        \centering
        \includegraphics[width=\linewidth, trim={10 12 10 10}, clip]{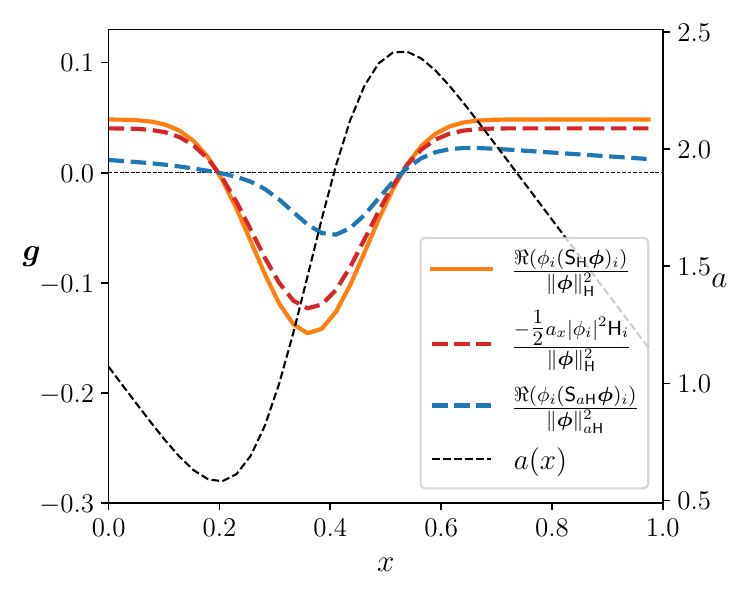}
        \caption{$\alpha=0.0$, $m=4$}
        \label{fig:rayleigh_mode1}
    \end{subfigure}
    \hfill
    \begin{subfigure}[t]{0.32\textwidth}
        \centering
        \includegraphics[width=\linewidth, trim={10 12 10 10}, clip]{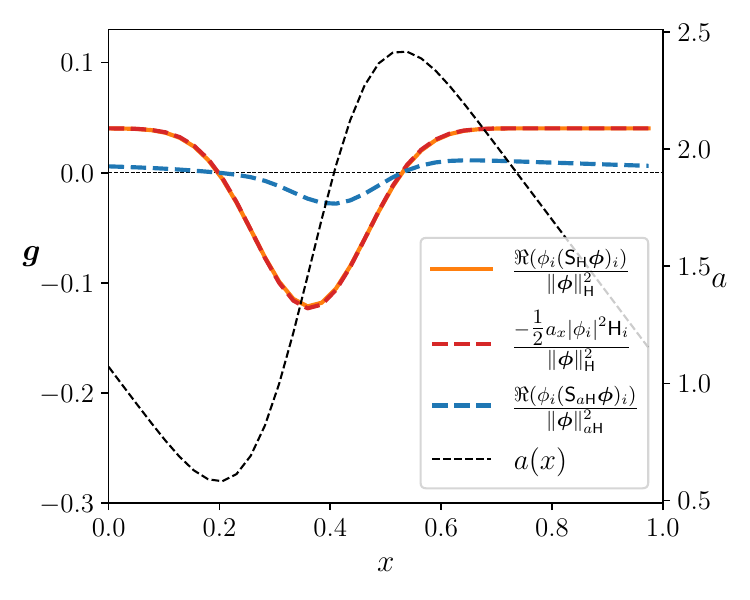}
        \caption{$\alpha=0.5$, $m=4$}
        \label{fig:rayleigh_mode2}
    \end{subfigure}
    \hfill
    \begin{subfigure}[t]{0.32\textwidth}
        \centering
        \includegraphics[width=\linewidth, trim={10 12 10 10}, clip]{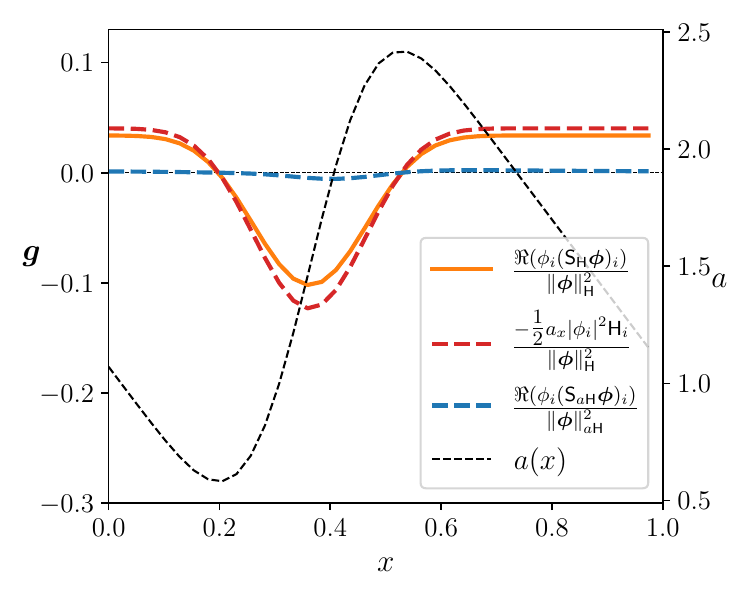}
        \caption{$\alpha=0.9$, $m=4$}
        \label{fig:rayleigh_mode3}
    \end{subfigure}

    \vspace{0.5em}

    \begin{subfigure}[t]{0.32\textwidth}
        \centering
        \includegraphics[width=\linewidth, trim={10 12 10 10}, clip]{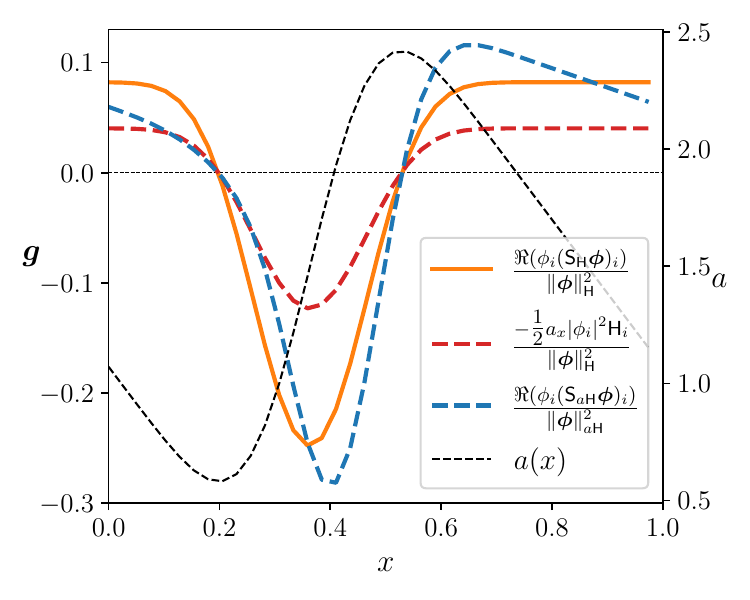}
        \caption{$\alpha=0.0$, $m=10$}
        \label{fig:rayleigh_mode4}
    \end{subfigure}
    \hfill
    \begin{subfigure}[t]{0.32\textwidth}
        \centering
        \includegraphics[width=\linewidth, trim={10 12 10 10}, clip]{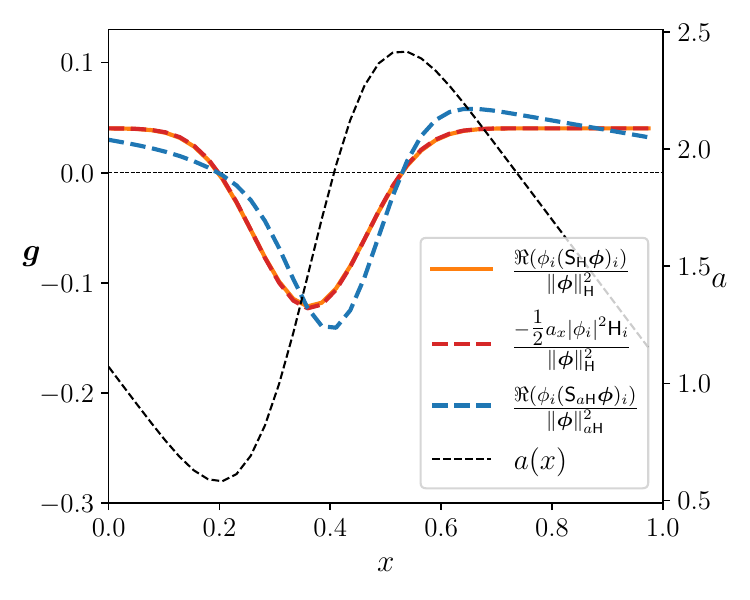}
        \caption{$\alpha=0.5$, $m=10$}
        \label{fig:rayleigh_mode5}
    \end{subfigure}
    \hfill
    \begin{subfigure}[t]{0.32\textwidth}
        \centering
        \includegraphics[width=\linewidth, trim={10 12 10 10}, clip]{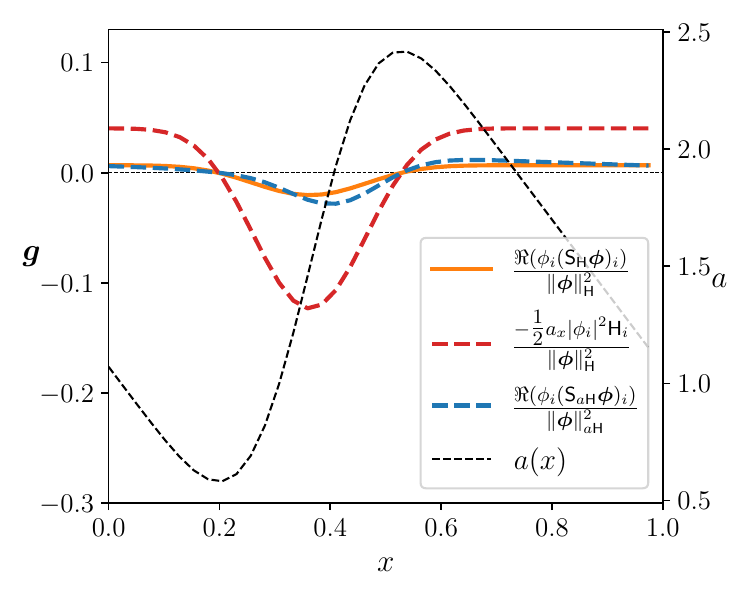}
        \caption{$\alpha=0.9$, $m=10$}
        \label{fig:rayleigh_mode6}
    \end{subfigure}

    \caption{Local Rayleigh quotients of the central-product split-form discretizations using a circulant 2\textsuperscript{nd}-order non-dissipative operator with $N=39$ nodes and Fourier modes $\vec{\phi} = e^{i\kappa x}$ with $\kappa = 2 \pi m$.}
    \label{fig:rayleigh_mode}
\end{figure}

When discussing ``additional local energy growth'' of the split-form schemes in Figure \ref{fig:eigvecs}, it is important to consider that, as explained in \S\ref{sec_modified_pde} and \S\ref{sec_frozen_coeff}, this characterization is in relation to the central scheme, rather than to the ``correct'' growth of the continuous PDE. In Figure \ref{fig:rayleigh_mode}, we again plot local Rayleigh quotients quantifying local energy and $a$-energy growth, but this time for individual Fourier modes $\vec{\phi} = e^{i \kappa x}$ and different values of splitting parameter $\alpha$. When the wavenumber is small (Figures \ref{fig:rayleigh_mode1}, \ref{fig:rayleigh_mode2}, and \ref{fig:rayleigh_mode3}), the local energy growth closely tracks the expected value of $- \tfrac{1}{2} a_x$ for all values of $\alpha$. Upon mesh refinement, the numerical growth further converges to the expected growth of the continuous PDE, and the numerical $a$-norm energy converges to zero, as expected. For $\alpha=0$, we observe \emph{additional} local energy growth when $a_x <0$ and additional energy dissipation when $a_x > 0$. As predicted in \S\ref{sec_frozen_coeff}, this additional growth is more pronounced for the higher wavenumber (less resolved) mode in Figure \ref{fig:rayleigh_mode4}. When $\alpha = 0.5$, however, the local energy growth is in very close agreement with the expected local growth, even for the high wavenumber mode. As we approach the central discretization with $\alpha = 0.9$, we see the opposite: the scheme underestimates both the local energy growth (when $a_x <0$) and dissipation (when $a_x > 0$). This may seem surprising, but is in fact consistent with the stability proof in \S\ref{sec_stability_proof}, which found the sharpest bound on energy growth when $\alpha=0.5$. 

To see this explicitly, consider the split-form semidiscrete operator from~\eqref{eq_general_linear_split_form},
\[
\mat{J}(\alpha) = - \alpha \mat{D} \mat{A} - (1-\alpha) \mat{A} \mat{D}  - (1-\alpha) \diag{\mat{D} \vec{a}} ,
\]
such that we can express the growth operators as
\begin{alignat*}{3}
    \mat{S}_\Hnrm &= - \tfrac{1}{2} \Hnrm \diag{\mat{D} \vec{a}} + (\alpha - \tfrac{1}{2}) \Hnrm \mat{\Theta} , \qquad 
    && \ \mat{\Theta} &&\defn \mat{D} \mat{A} - \mat{A} \mat{D} - \diag{\mat{D} \vec{a}} , \\
    \mat{S}_{a \Hnrm} &= (1 - \alpha) \Hnrm \mat{\Theta}_a , 
    && \mat{\Theta}_a &&\defn \mat{D} \mat{A}^2 - \mat{A}^2 \mat{D} - 2 \mat{A} \diag{\mat{D} \vec{a}} .
\end{alignat*}
Therefore, local energy growth is most consistent with the continuous PDE when $\alpha = \tfrac{1}{2}$; otherwise for $\alpha \neq \tfrac{1}{2}$ it is modified by the product-rule defect term~$\mat{\Theta}$. Alternatively, local $a$-norm energy growth is most consistent with the continuous PDE when $\alpha = 1$; otherwise for $\alpha \neq 1$ it is modified by a similar product-rule defect term~$\mat{\Theta}_a$. Prioritizing the $a$-norm, we used the central scheme ($\alpha =1$) as a reference for the modified PDE analysis in \S\ref{sec_modified_pde}, which unsurprisingly yielded a dissipation coefficient $\nu \propto 1 - \alpha$. However, we could just as well have prioritized consistency of the standard energy, after which a similar analysis performed with the $\alpha = \tfrac{1}{2}$ scheme as a reference would yield a dissipation coefficient $\nu \propto \tfrac{1}{2} - \alpha$. We can therefore ask: which is more important, consistency of local energy growth, or consistency of local $a$-energy growth? Perhaps some combination is ideal, such as $\alpha=\tfrac{2}{3}$, which is obtained by linearizing entropy-stable discretizations of the Burgers equation. As mentioned in \S\ref{sec_modified_pde}, if we are only concerned with numerical accuracy, there is no guarantee that $\alpha=1$ is optimal. The primary benefit of the central scheme is rather its $a$-norm stability estimate and its telescoping $a$-energy property, which translate to spectral stability. If these instabilities can be otherwise controlled, however, our results suggest that $\alpha \neq 1$ splittings may be just as suitable as the central scheme for solving the variable-coefficient advection equation.

Before proceeding with numerical experiments exploring the accuracy of the various discretizations and the practical impact of these local linear instabilities, we summarize here what we have shown so far.
\begin{itemize}
    \item The central scheme $\alpha=1$ exactly conserves the $a$-norm energy. Comparatively, the product scheme introduces additional local energy growth (or dissipation, depending on the sign of $a_x(x)$) that leads to a loss of exact $a$-norm energy conservation.
    \item For well-resolved physical modes, this additional growth vanishes under mesh refinement at design order. However, under-resolved unphysical modes, in particular boundary-localized eigenmodes, can excite unbalanced energy growth that persists under mesh refinement, leading to global instabilities.
    \item Unphysical oscillatory interior modes are less of a concern, because they typically benefit from global energy growth / dissipation cancellations. Furthermore, they can be effectively controlled by the inclusion of volume dissipation.
    \item For long-time simulations, splittings with $\alpha \neq 1$ may develop numerical instabilities if unstable modes are not properly controlled (\eg through numerical dissipation). However, until these (typically initially small) modes grow sufficiently, there is reason to expect negligible accuracy differences between different splittings $\alpha \in [0,1]$.
    \item At this point, it is unclear whether these conclusions, which were inspired by the `central plus (anti-)dissipation' perspective and tested numerically for the (linear) product scheme, also apply to nonlinear flux-differencing schemes such as the geometric and logarithmic fluxes. In \S\ref{sec:time-varying}, we will show that the modal perturbation-growth mechanism studied here does not appear to carry over directly to these nonlinear schemes.
\end{itemize}

\subsection{Numerical Accuracy Experiments of the Central-Product Split-Form with a Time-Invariant but Spatially-Varying Coefficient} \label{sec:time_invariant_accuracy_experiments}

With the theoretical understanding of the nature of local linear instabilities gained from the previous sections, we now aim to quantify the practical impact they have on simulations. That is, when does the additional unbalanced energy growth introduced by $\alpha \neq 1$ splittings affect solution accuracy? Once again, we defer experiments for the nonlinear geometric and logarithmic schemes to later sections as they will require more sophisticated tools. We solve~\eqref{eq_lin_conv} with the skewed-sinusoid variable coefficient~\eqref{eq_skewed_sin} and a Gaussian initial condition. We employ 8\textsuperscript{th}-order adaptive explicit time marching of~\cite{Hairer1993} so that temporal errors are negligible, and plot both the $\Hnrm$-norm and $L^\infty$ errors of the numerical solution as a function of $t$. We also plot the predicted local linear instability growth by projecting the initial condition onto the linearly unstable modes, then solving for modal growth in $t$.

\begin{figure}[t]
    \centering
    \captionsetup[subfigure]{justification=centering,singlelinecheck=false}

    \begin{subfigure}[t]{0.32\textwidth}
        \centering
        \includegraphics[width=\linewidth, trim={12 12 8 8}, clip]{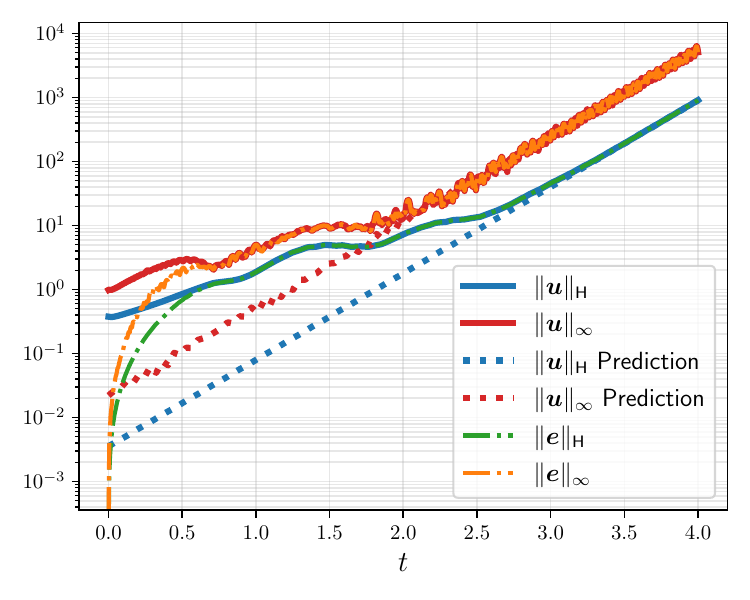}
        \caption{$\alpha=0$, $p=1$, $N=40$}
        \label{fig:error_csbp_nd1}
    \end{subfigure}
    \hfill
    \begin{subfigure}[t]{0.32\textwidth}
        \centering
        \includegraphics[width=\linewidth, trim={12 12 8 8}, clip]{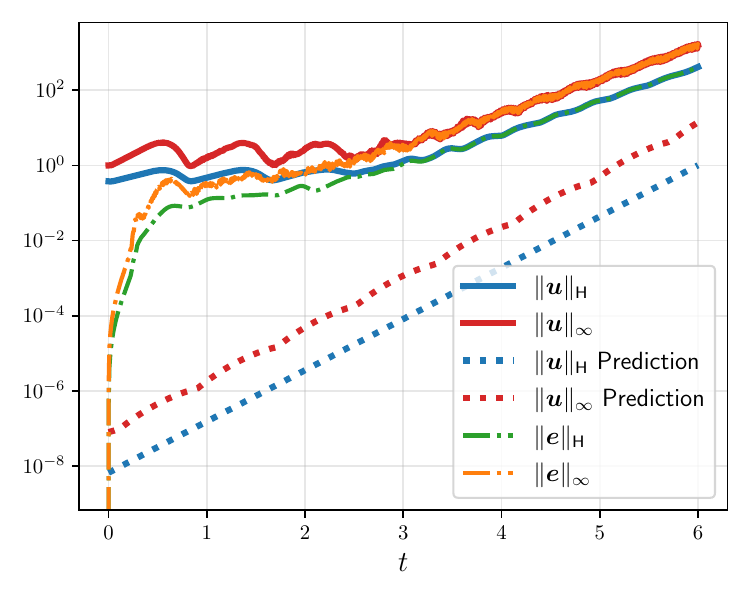}
        \caption{$\alpha=0$, $p=2$, $N=100$}
        \label{fig:error_csbp_nd2}
    \end{subfigure}
    \hfill
    \begin{subfigure}[t]{0.32\textwidth}
        \centering
        \includegraphics[width=\linewidth, trim={12 12 8 8}, clip]{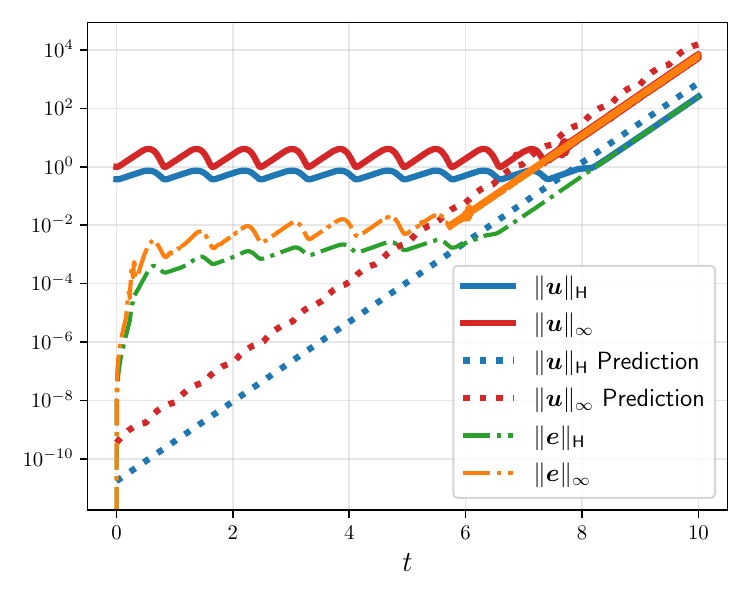}
        \caption{$\alpha=0$, $p=2$, $N=400$}
        \label{fig:error_csbp_nd3}
    \end{subfigure}

    \vspace{0.5em}

    \begin{subfigure}[t]{0.32\textwidth}
        \centering
        \includegraphics[width=\linewidth, trim={12 12 8 8}, clip]{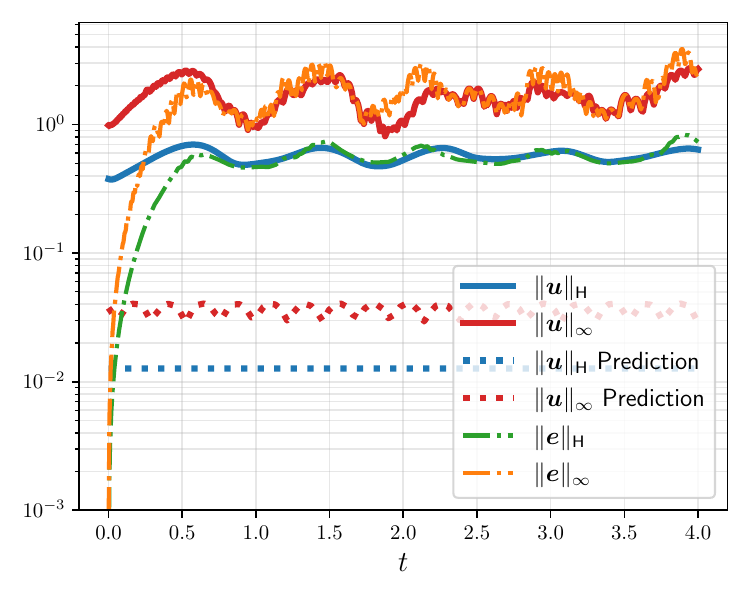}
        \caption{$\alpha=1$, $p=1$, $N=40$}
        \label{fig:error_csbp_nd4}
    \end{subfigure}
    \hfill
    \begin{subfigure}[t]{0.32\textwidth}
        \centering
        \includegraphics[width=\linewidth, trim={12 12 8 8}, clip]{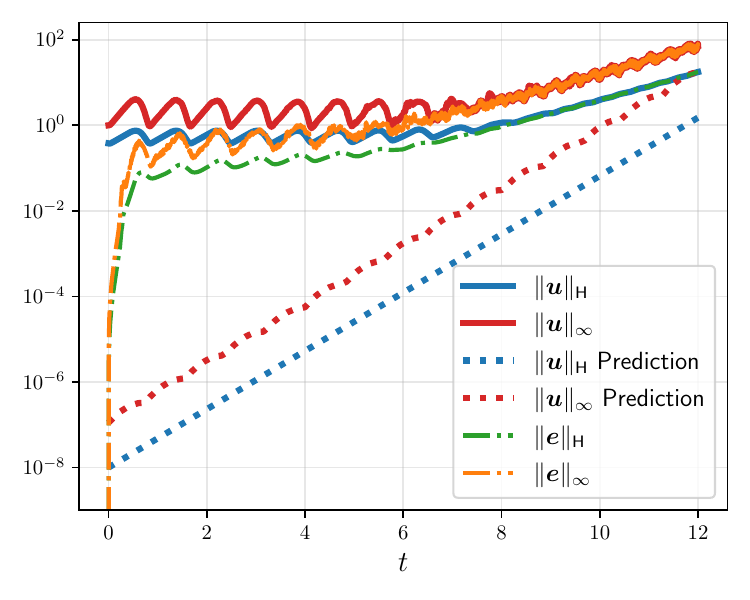}
        \caption{$\alpha=\tfrac{1}{2}$, $p=2$, $N=100$}
        \label{fig:error_csbp_nd5}
    \end{subfigure}
    \hfill
    \begin{subfigure}[t]{0.32\textwidth}
        \centering
        \includegraphics[width=\linewidth, trim={12 12 8 8}, clip]{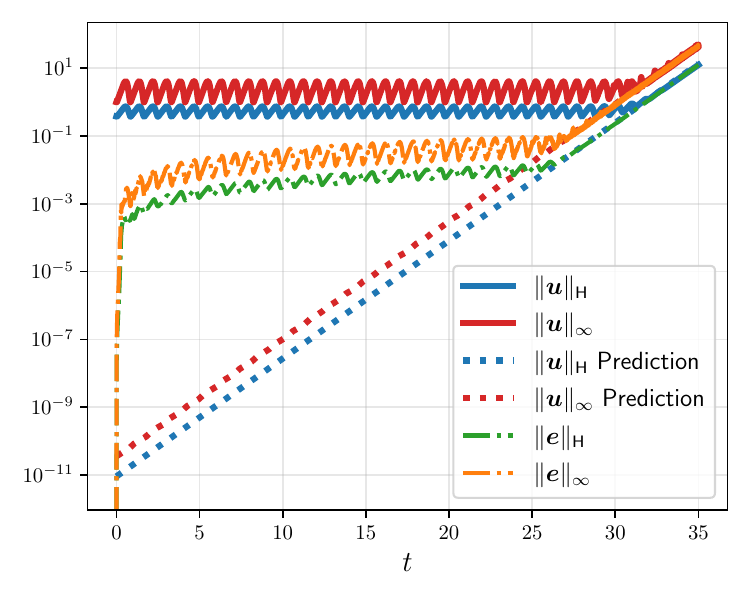}
        \caption{$\alpha=0$, continuous CSBP $p=2$, $N=400$}
        \label{fig:error_csbp_nd6}
    \end{subfigure}

    \caption{Worst case scenario central-product split-form scheme error growth with non-dissipative CSBP operators.}
    \label{fig:error_csbp_nd}
\end{figure}

In Figure~\ref{fig:error_csbp_nd} we consider the worst case scenario for local linear instability growth: CSBP operators with no interface or volume dissipation. As predicted in~\cite{gassner_stability_2022}, the unstable modes grow exponentially according to the real part of the eigenvalue, eventually becoming the dominant source of error. For coarse meshes, as in Figure~\ref{fig:error_csbp_nd1}, other sources of numerical error fully saturate the solution before the unstable modes have a chance to grow. This is true even for the $\alpha=1$ central scheme in Figure~\ref{fig:error_csbp_nd4}. For finer meshes, however, the exponential growth of linear instabilities eventually becomes a cause for concern. Shown in Figures~\ref{fig:error_csbp_nd2} and~\ref{fig:error_csbp_nd5}, the crossover point does not always come from the most unstable mode. Depending on the initial condition, a less unstable spurious mode can be more strongly excited such that in a short time it introduces more error. Figure~\ref{fig:error_csbp_nd5} shows that the strength of these instabilities (and hence error growth) decreases as $\alpha \rightarrow 1$, but the behaviour remains qualitatively similar. The fine mesh in Figure~\ref{fig:error_csbp_nd5} shows that the exponential growth of even initially small unstable modes eventually destroys the accuracy of a simulation. This is particularly concerning for high-order methods applied to long-time wave propagation, since their usual advantage of having small, slowly accumulating numerical errors is ultimately negated if this error component grows exponentially in time. Figure~\ref{fig:error_csbp_nd6} shows the same result, but using the continuous SBP operators of~\cite{hicken_multidimensional_2016, Hicken2020}, indicating that local linear instabilities are not necessarily generated by the presence of SATs, but are rather strongly enhanced by the modified stencils at block interfaces.

\begin{figure}[t]
    \centering
    \captionsetup[subfigure]{justification=centering,singlelinecheck=false}

    \begin{subfigure}[t]{0.32\textwidth}
        \centering
        \includegraphics[width=\linewidth, trim={12 12 8 8}, clip]{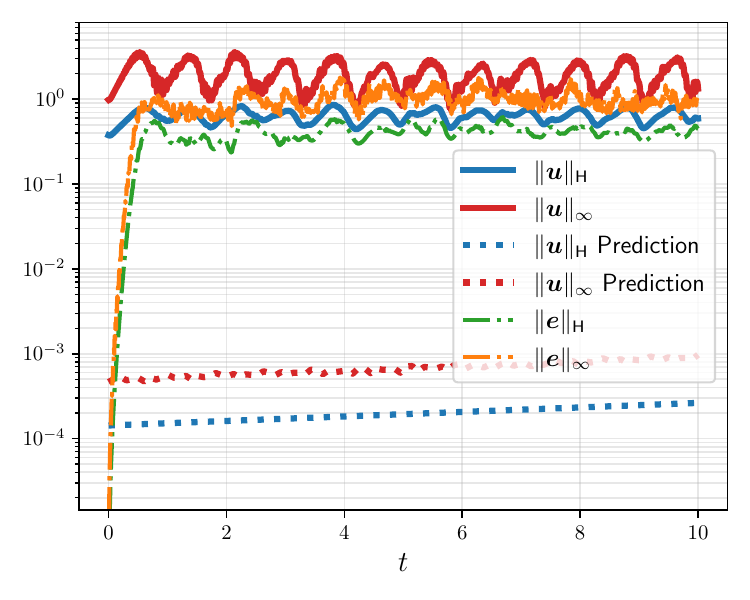}
        \caption{$N=39$, $\Re{(\lambda)} = 6.1 \times 10^{-2}$}
        \label{fig:error_circulant_nd1}
    \end{subfigure}
    \hfill
    \begin{subfigure}[t]{0.32\textwidth}
        \centering
        \includegraphics[width=\linewidth, trim={12 12 8 8}, clip]{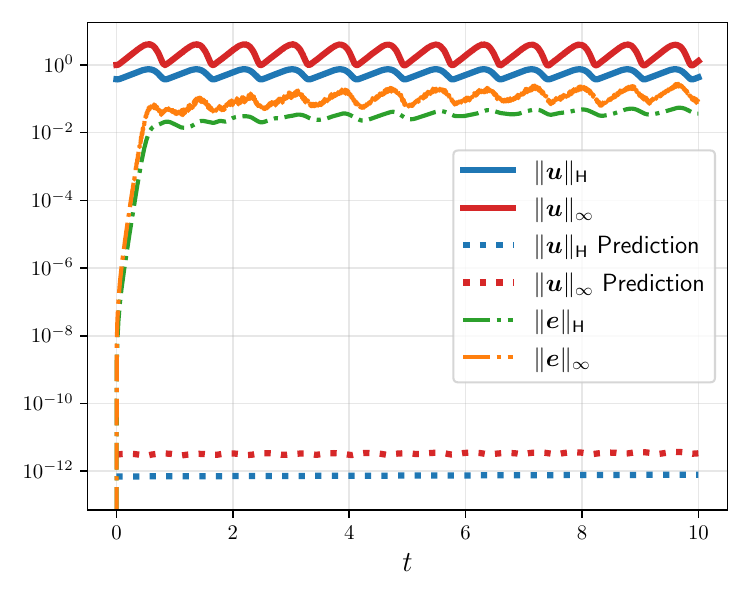}
        \caption{$N=99$, $\Re{(\lambda)} = 1.2 \times 10^{-2}$}
        \label{fig:error_circulant_nd2}
    \end{subfigure}
    \hfill
    \begin{subfigure}[t]{0.32\textwidth}
        \centering
        \includegraphics[width=\linewidth, trim={12 12 8 8}, clip]{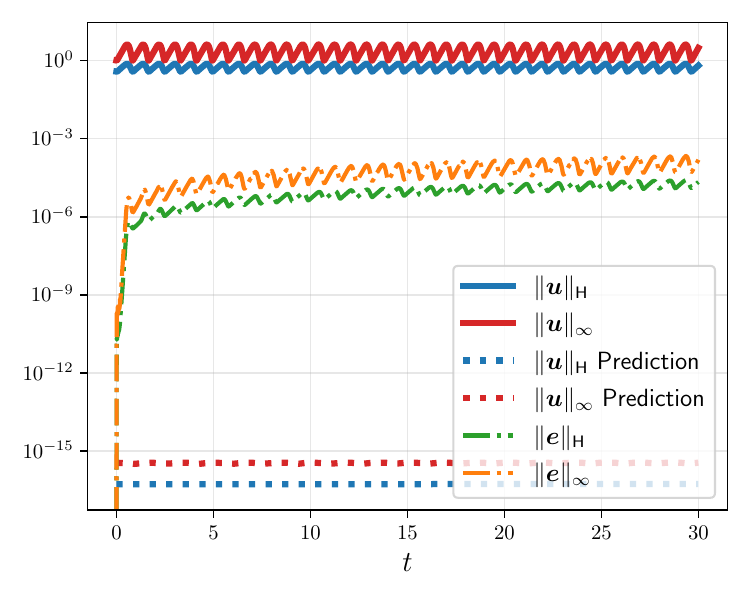}
        \caption{$N=399$, $\Re{(\lambda)} = 8.2 \times 10^{-4}$}
        \label{fig:error_circulant_nd3}
    \end{subfigure}

    \caption{Product scheme ($\alpha = 0$) error growth with non-dissipative 8\textsuperscript{th}-order circulant operators.}
    \label{fig:error_circulant_nd}
\end{figure}

\begin{figure}[t]
    \centering
    \captionsetup[subfigure]{justification=centering,singlelinecheck=false}

    \begin{subfigure}[t]{0.32\textwidth}
        \centering
        \includegraphics[width=\linewidth, trim={12 12 8 8}, clip]{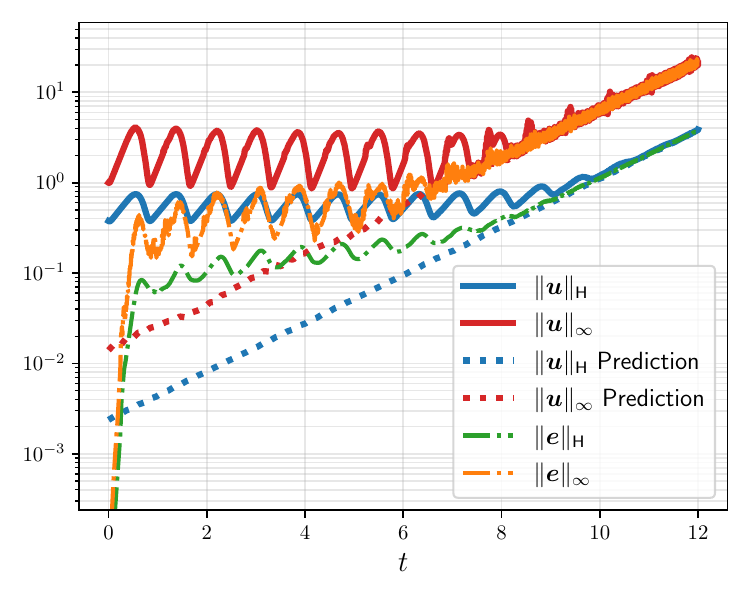}
        \caption{CSBP, $N=100$, dissipative SAT. $\Re{(\lambda)}_{\max} = 6.1 \times 10^{-1}$}
        \label{fig:error_diss1}
    \end{subfigure}
    \hfill
    \begin{subfigure}[t]{0.32\textwidth}
        \centering
        \includegraphics[width=\linewidth, trim={12 12 8 8}, clip]{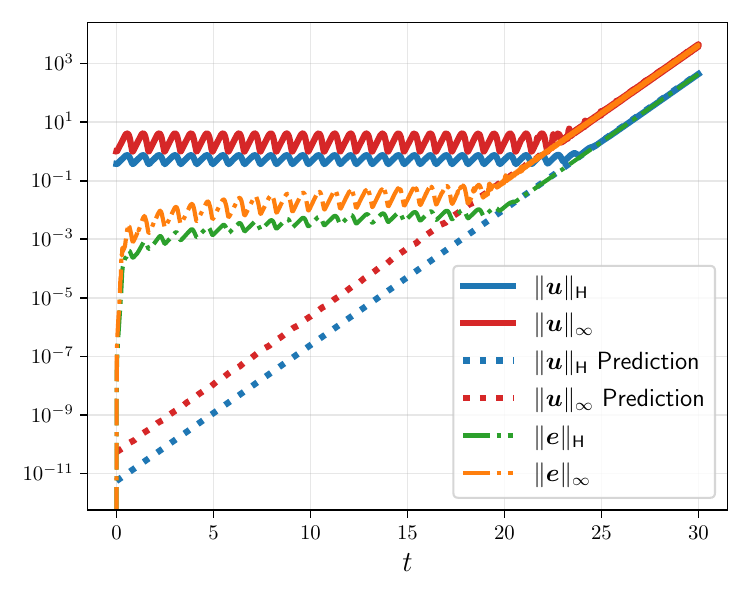}
        \caption{CSBP, $N=400$, dissipative SAT. $\Re{(\lambda)}_{\max} = 1.1$}
        \label{fig:error_diss2}
    \end{subfigure}
    \hfill
    \begin{subfigure}[t]{0.32\textwidth}
        \centering
        \includegraphics[width=\linewidth, trim={12 12 8 8}, clip]{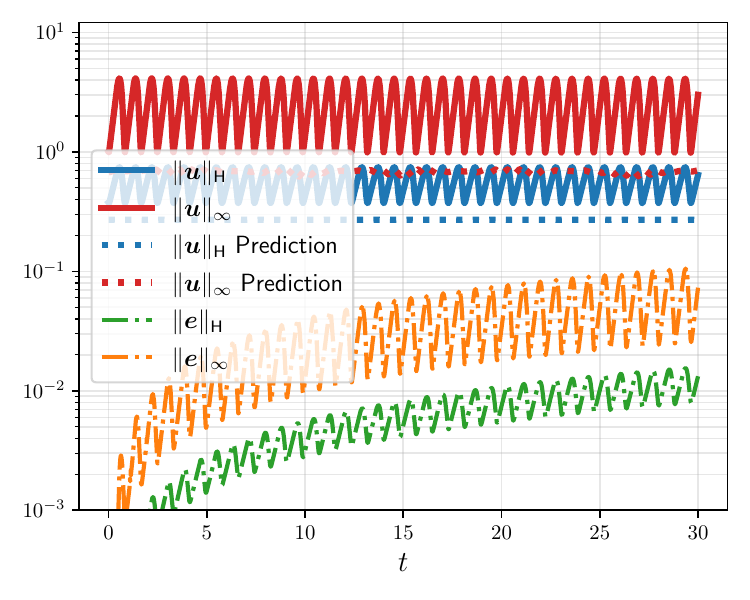}
        \caption{CSBP, $N=400$, volume and SAT dissipation. $\Re{(\lambda)}_{\max} = 1.9 \times 10^{-7}$}
        \label{fig:error_diss3}
    \end{subfigure}

    \vspace{0.5em}

    \begin{subfigure}[t]{0.32\textwidth}
        \centering
        \includegraphics[width=\linewidth, trim={12 12 8 8}, clip]{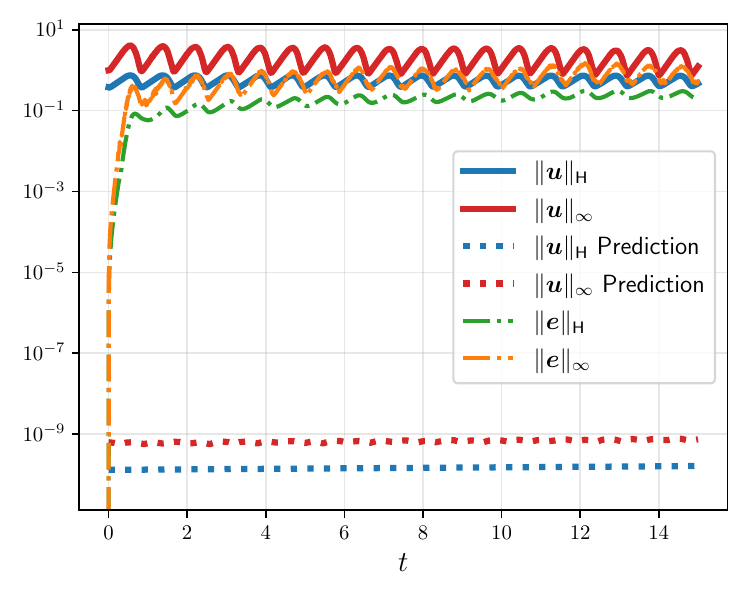}
        \caption{Mattsson, $N=100$, dissipative SAT. $\Re{(\lambda)}_{\max} = 1.4 \times 10^{-2}$}
        \label{fig:error_diss4}
    \end{subfigure}
    \hfill
    \begin{subfigure}[t]{0.32\textwidth}
        \centering
        \includegraphics[width=\linewidth, trim={12 12 8 8}, clip]{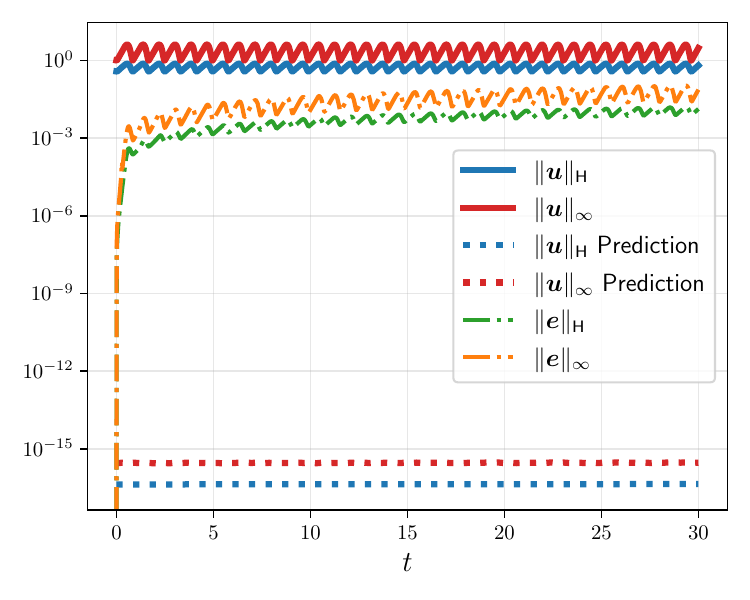}
        \caption{Mattsson, $N=400$, dissipative SAT. $\Re{(\lambda)}_{\max} = 9.9 \times 10^{-4}$}
        \label{fig:error_diss5}
    \end{subfigure}
    \hfill
    \begin{subfigure}[t]{0.32\textwidth}
        \centering
        \includegraphics[width=\linewidth, trim={12 12 8 8}, clip]{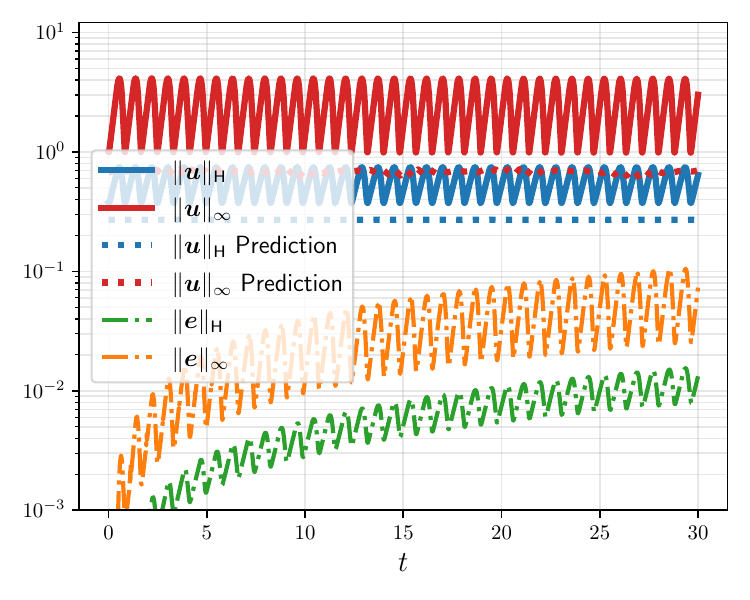}
        \caption{Continuous CSBP, $N=400$, volume dissipation. $\Re{(\lambda)}_{\max} = 1.9 \times 10^{-7}$}
        \label{fig:error_diss6}
    \end{subfigure}

    \caption{Product scheme ($\alpha = 0$, $p=2$) error growth with interface and volume dissipation.}
    \label{fig:error_diss}
\end{figure}

\begin{figure}[t]
    \centering
    \captionsetup[subfigure]{justification=centering,singlelinecheck=false}

    \begin{subfigure}[t]{0.32\textwidth}
        \centering
        \includegraphics[width=\linewidth,trim={12 12 10 10}, clip]{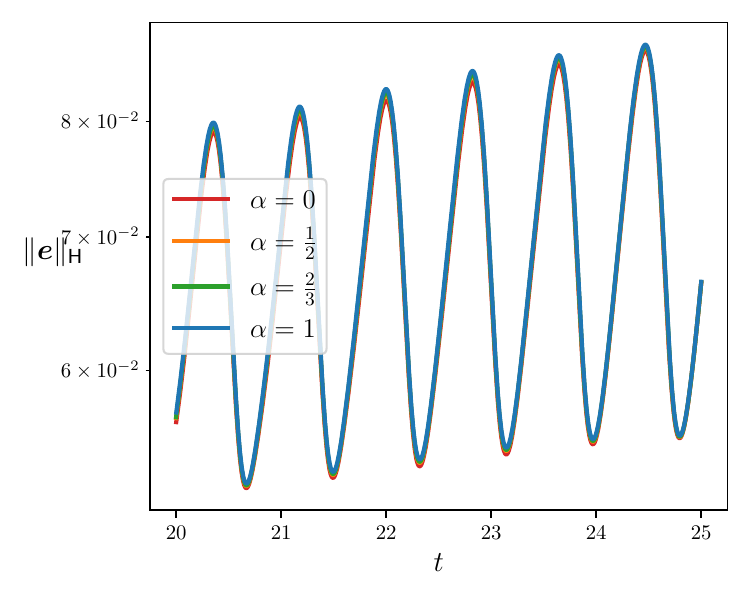}
        \caption{Mattsson $p=4$, $N=100$, dissipative SAT}
        \label{fig:error_compare1}
    \end{subfigure}
    \hfill
    \begin{subfigure}[t]{0.32\textwidth}
        \centering
        \includegraphics[width=\linewidth,trim={12 12 10 10}, clip]{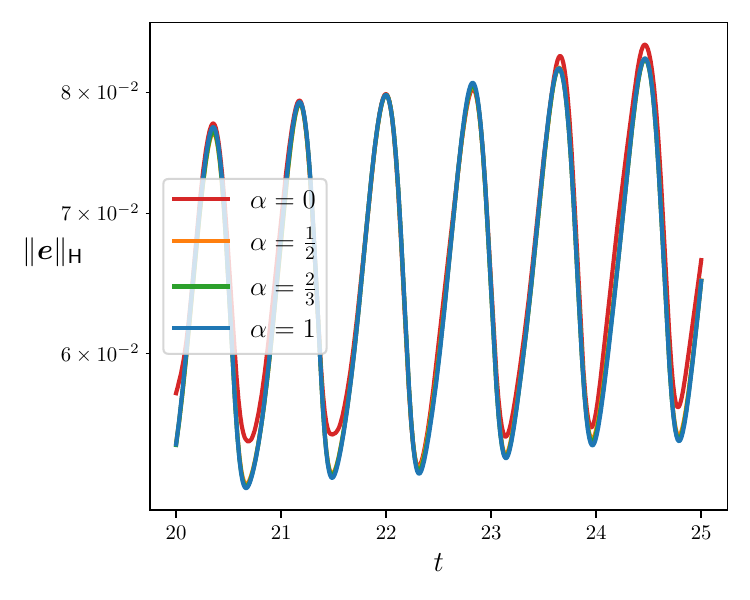}
        \caption{Circulant 8\textsuperscript{th}-order, $N=99$, non-dissipative}
        \label{fig:error_compare2}
    \end{subfigure}
    \hfill
    \begin{subfigure}[t]{0.32\textwidth}
        \centering
        \includegraphics[width=\linewidth,trim={12 12 10 10}, clip]{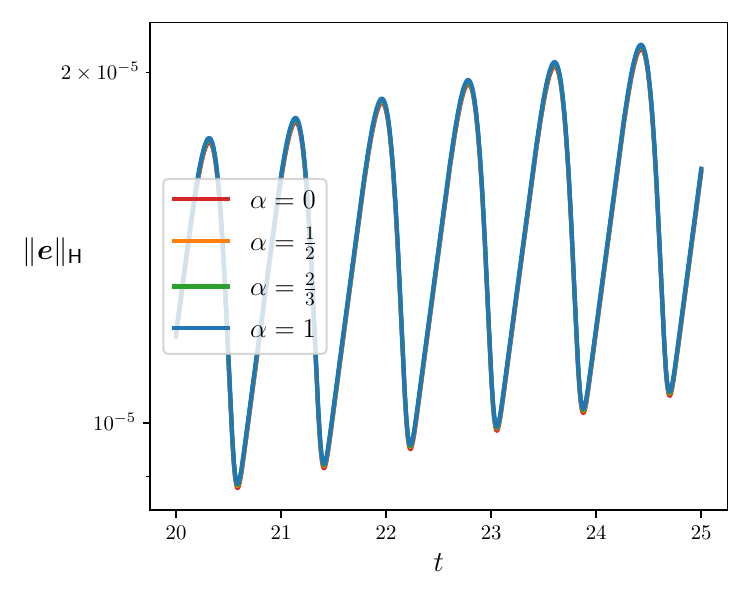}
        \caption{Circulant 8\textsuperscript{th}-order, $N=399$, non-dissipative}
        \label{fig:error_compare3}
    \end{subfigure}

    \caption{Comparison of $\Hnrm$-norm solution errors between different central-product splittings $\alpha \in [0,1]$.}
    \label{fig:error_compare}
\end{figure}

However, we have already shown that this worst-case scenario often overstates the practical impact of local linear instabilities. Therefore, we now use the lessons learned from previous sections to reassess their severity. We first showed that instabilities resulting from interior dynamics alone are weak and tend to converge to zero as the mesh is refined. This is confirmed in Figure \ref{fig:error_circulant_nd}, where we emphasize the importance of long-time solutions by taking an 8\textsuperscript{th}-order circulant operator with no volume dissipation. The linear instabilities are of the same magnitude observed for the 2\textsuperscript{nd}-order case in Figure~\ref{fig:eigvals4}, small enough that the error growth is effectively negligible compared to other sources of error. We then argued that interface and volume dissipation should effectively control linear instabilities. For CSBP operators, Figures~\ref{fig:error_diss1} and~\ref{fig:error_diss2} demonstrate that including interface dissipation alone only partially suppresses the boundary-localized unstable modes, consistent with behaviour observed in \S\ref{sec:spectra}. Consequently, local linear instabilities can still grow sufficiently to impact long-time simulations. However, Figures \ref{fig:error_diss4} and \ref{fig:error_diss5} show that if we instead use Mattsson operators, which naturally suppress spurious modes arising from the boundary, then interface dissipation does suppress all boundary-localized unstable modes. The remaining unstable modes are nearly identical to those of the circulant operators, and so once again exhibit negligible linear instability error growth. Alternatively, we can also suppress the remaining unstable modes by adding a small amount of volume dissipation using the framework of \cite{dissipation} with $s=3$ and $\varepsilon=10^{-3}$ for $p=2$. For both an SBP-SAT discretization in Figure~\ref{fig:error_diss3} and a continuous-SBP discretization in~\ref{fig:error_diss6}, the error growth from linear instabilities is negligible. Finally, in Figure~\ref{fig:error_compare} we show that accuracy differences between different splitting parameters $\alpha$ are negligible. That is, once the linear instabilities are suppressed, all splittings $\alpha \in [0,1]$ have nearly identical numerical errors. 

Once again, we omit results for spectral-element operators, noting that while simulations with no interface dissipation show the same concerning exponential error growth as the worst case scenario of Figure~\ref{fig:error_csbp_nd} (consistent with the results of \cite{gassner_stability_2022}), including interface dissipation recovers linearly stable schemes. Then similar to Figure~\ref{fig:error_compare}, accuracy differences between different splitting parameters $\alpha$ become negligible.

Before moving on to the more interesting case of a time-varying baseflow (and hence nonlinear discretizations, such as the geometric scheme), we again summarize the findings of this section.
\begin{itemize}
    \item Local linear instabilities are primarily a concern for discretizations with modified boundary stencils, introducing boundary-localized spurious modes that excite unbalanced energy growth. Once the boundary-localized modes are suppressed, all splittings $\alpha \in [0,1]$ exhibit comparable numerical error.
    \item When instabilities are not effectively suppressed, error grows exponentially as predicted in \cite{gassner_stability_2022}. While this may not be a stability concern in the context of perturbation growth of entropy-stable schemes, this unphysical growth can still saturate the numerical error depending on the initial strength of the perturbation, the simulation time, and the relative magnitude of other sources of error.
    \item Suppressing boundary-localized modes is sufficient to make the error introduced by local linear instabilities negligible. For spectral-element or Mattsson operators, interface dissipation alone can be sufficient. For CSBP and continuous-SBP discretizations, volume dissipation can be necessary.
\end{itemize}

\section{Time-Varying Coefficients: Floquet Analysis} \label{sec:time-varying}

For the geometric and logarithmic fluxes, the Jacobians $\mat{J}(\tilde{\bu})$ acting on perturbations are time-varying linear evolution operators, even if $a$ is constant in both space and time. In contrast to the numerical experiments of \S\ref{sec:time_invariant_accuracy_experiments}, therefore, the baseflow about which one linearizes is changing in time. In \cite{gassner_stability_2022} and \cite{ranocha_2022_preventing_pressure}, time-varying baseflows were analyzed by considering snapshots in time. That is, if at any point in time the eigenvalues of $\mat{J}(\tilde{\bu})$ exhibited positive real parts, perturbations were deemed to grow unphysically, similar to the worst case scenario experiments in \S\ref{sec:time_invariant_accuracy_experiments}. However, the instantaneous spectra of time-varying linear operators do not determine the growth of solutions. It is possible for $\mat{J}(\tilde{\bu})$ to have eigenvalues with positive real parts for all $t$, yet still have decaying solutions (and vice versa)~\cite{nonautonomous_systems}. Over short times, transient growth (due to the interference of nonorthogonal eigenmodes) can dominate the asymptotic behaviour, and if the eigenvectors also rotate in time, then long-time perturbation dynamics can differ substantially from predictions based solely on instantaneous eigenspectra. The frozen-coefficient analysis in~\S\ref{sec_frozen_coeff} showed that this scenario is plausible, since the linearly unstable eigenmodes oscillate on time scales much shorter than those associated with their growth rates. In general, time-varying linear systems are not completely understood~\cite{Chicone2006}. However, for the special case where $\mat{J}(\tilde{\bu})$ is periodic in time, we can use Floquet analysis to generalize the techniques of the time-invariant case while including the aforementioned transient behaviour.

\subsection{Floquet Theory} \label{sec:floquet}

In this section we briefly introduce relevant concepts of Floquet theory. For a more comprehensive treatment, we refer the reader to Appendix \ref{app_floquet}, or further to \cite{Chicone2006}. Consider the time-varying linear system
\begin{equation} \label{eq_time_periodic_system}
    \der[\bv]{t} = \mat{J}(t) \bv , \quad \bv \in \mathbb{R}^N , \quad \mat{J}(t) = \mat{J}(t + T) \in \mathbb{R}^{N \times N} .
\end{equation}
Floquet theory dictates that \eqref{eq_time_periodic_system} can be equivalently expressed as a constant-coefficient linear system in a periodically rotating frame. More precisely, there exists an
invertible \(T\)-periodic matrix \(\mat{P}(t)\in\mathbb{C}^{N\times N}\) with $\mat{P}(T) = \mat{P}(0) = \Idty$ and a constant matrix \(\mat{B}\in\mathbb{C}^{N\times N}\) such that
\[
    \bv(t)=\mat{P}(t)\tilde{\bv}(t)
    \quad \text{and} \quad 
    \der[\tilde{\bv}]{t}=\mat{B}\tilde{\bv}
    \quad \Rightarrow \quad 
    \bv(t)=\mat{P}(t) e^{\mat{B} t}\bv(0) .
\]
Constructing \(\mat{P}(t)\) and \(\mat{B}\) is nontrivial, as
they are not obtained pointwise from \(\mat{J}(t)\), but rather through the time-integrated evolution operator $\mat{\Psi}(t) = \mat{P}(t) e^{\mat{B} t} \in \mathbb{R}^{N \times N}$, which satisfies $\bv(t) = \mat{\Psi}(t)\bv(0)$. 

Of particular importance is the \emph{monodromy matrix} $\mat{\Psi}(T) = e^{\mat{B} T}$, because repeated application of $\mat{\Psi}(T)$ yields
\[
\bv(nT) = \mat{\Psi}(T)^n \bv(0) \quad \text{for all} \ n \in \mathbb{N}_0 .
\]
It follows that $\mat{\Psi}(T)$ characterizes the asymptotic growth of $\bv(t)$. Its eigenvalues $\rho \in \mathbb{C}$ are called \emph{characteristic multipliers} of \eqref{eq_time_periodic_system}, because they quantify the amount by which the eigenvectors of $\mat{\Psi}(T)$---which serve as the initial conditions of the \emph{Floquet modes}---grow after one period. Assuming that $\mat{\Psi}(T)$ is diagonalizable, the condition $\abs{\rho} \leq 1$ implies stability of the solution $\bv$. Equivalently, we can examine the eigenvalues $\lambda$ of the matrix $\mat{B}$, called \emph{Floquet exponents}, to asses stability via $\Re{(\lambda)} \leq 0$, noting that $\lambda = \tfrac{1}{T} \log(\rho)$, up to the additive shifts $\lambda = \lambda + 2 \pi i k / T, \ k \in \mathbb{Z}$.

In the special case that $\mat{J}$ is time-invariant, the system \eqref{eq_time_periodic_system} is trivially $T$-periodic and we recover $\mat{B}=\mat{J}$, meaning the Floquet exponents are simply the eigenvalues of $\mat{J}$ (up to the additive shifts). In this sense, Floquet theory generalizes eigenvalue analysis of constant-coefficient systems to time-varying periodic systems. Though the asymptotic behaviour of $\vec{v}$ is determined by the Floquet exponents, significant growth may occur even if $\Re{(\lambda)} \leq 0$ due to transient dynamics, which is described by the singular values $\sigma$ of $\mat{\Psi}(t)$.

Finding the monodromy matrix $\mat{\Psi}(T)$ can be difficult, as it requires solving a matrix initial value problem. We choose to numerically approximate $\mat{\Psi}(T)$ using the ``exponential midpoint'' method described in Appendix \ref{app_floquet}, which involves composing $K$ snapshots $\mat{J}(t_j)$ taken at discrete times $t_j \in [0,T]$, similar in spirit to time marching. We can approximate $\mat{\Psi}(T)$ arbitrarily well by taking a sufficiently large $K$. In all reported numerical experiments, $K$ was chosen large enough such that further increases resulted in negligible changes to the Floquet multipliers.

\subsection{Numerical Experiments of the Geometric and Logarithmic Schemes for the Constant-Coefficient Linear Advection Equation} \label{sec_logarithmic_LCE_experiments}

The simplest scenario to which we can apply Floquet theory to assess perturbation growth with a time-varying baseflow is the entropy-stable geometric and logarithmic schemes applied to the constant-coefficient linear advection equation. The `central plus (anti-)dissipation' argument of \S\ref{sec_modified_pde} suggests that for the geometric scheme, even when $a(x) = a$, perturbations evolve according to a modified PDE with anti-dissipation. With central second-order operators, we showed in \eqref{eq_modified_pde_geom} that this modified PDE is
\begin{equation} \label{eq_modified_pde_geom_const}
    \begin{gathered}
\pder[\fnc{V}]{t} + \pder{x} \left[ a \left( 1 + \beta(x,t) \right)  \fnc{V} \right] = \pder{x} \left( \nu(x,t) \pder[\fnc{V}]{x} \right) , \\
 \beta(x,t) \defn \frac{h^2}{8} \frac{1}{\fnc{U}^2} \left(\pder[\fnc{U}]{x}\right)^2  , \quad \nu(x,t) \defn a \frac{h^2}{4} \frac{1}{\fnc{U}} \pder[\fnc{U}]{x},
    \end{gathered}
\end{equation}
such that antidissipation is introduced in regions of decreasing baseflow $\fnc{U}$. The argument can then be followed to assert that the linearized Jacobians will exhibit eigenvalues with positive real parts, indicating unphysical perturbation instabilities. Indeed, in \cite{ranocha_2022_preventing_pressure} it was shown that the logarithmic entropy-stable scheme exhibits local linear instabilities when applied to the constant-coefficient advection equation. This problem is of particular interest because it is analogous to the Euler density-wave problem posed in \cite{gassner_stability_2022}, where it was argued that entropy-stable schemes crash due to unstable perturbation growth. In the following, we alternatively analyze the problem using Floquet analysis, showing that the schemes exhibit no significant perturbation growth in time, and instead behave according to the sharp perturbation bound \eqref{eq_geom_pert_bound}. 

\begin{figure}[t]
    \centering
    \captionsetup[subfigure]{justification=centering,singlelinecheck=false}

    \begin{subfigure}[t]{0.28\textwidth}
        \includegraphics[height=3.0cm,
        trim={3mm 3mm 3mm 3mm},clip,left]{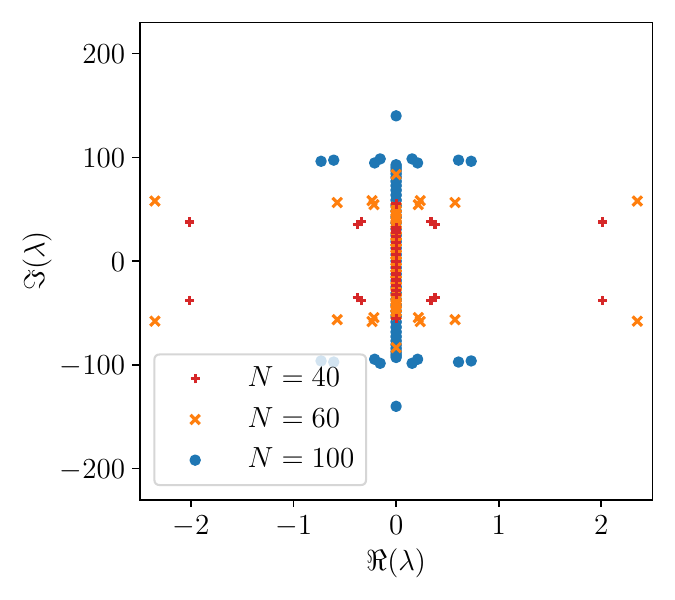}
    \end{subfigure}
    \hfill
    \begin{subfigure}[t]{0.225\textwidth}
        \includegraphics[height=3.0cm,
        trim={21mm 3mm 3mm 3mm},clip,center]{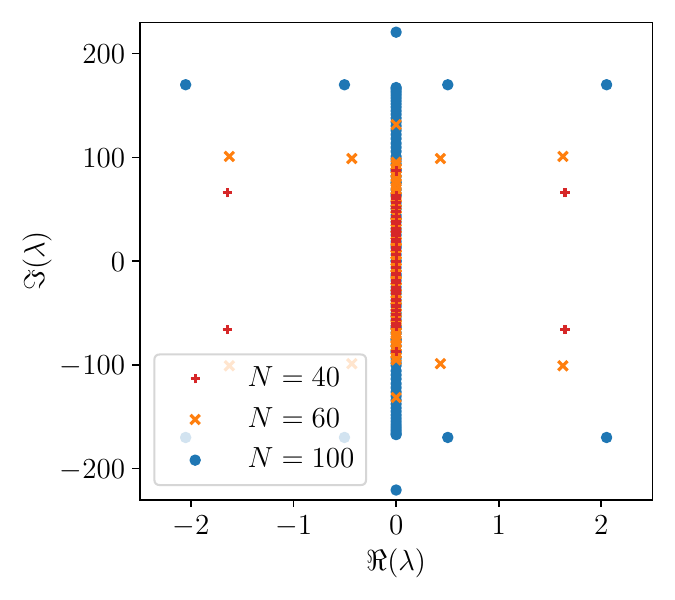}
    \end{subfigure}
    \hfill
        \begin{subfigure}[t]{0.225\textwidth}
        \includegraphics[height=3.0cm,
        trim={21mm 3mm 3mm 3mm},clip,center]{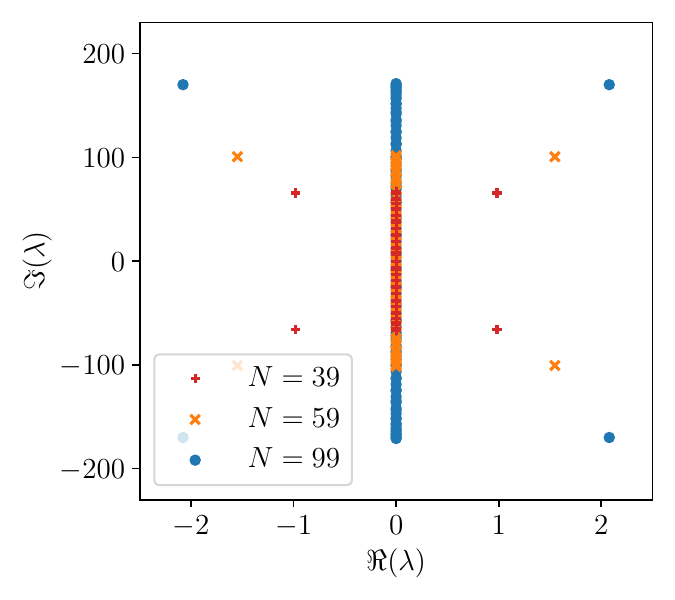}
    \end{subfigure}
    \hfill
        \begin{subfigure}[t]{0.23\textwidth}
        \includegraphics[height=3.0cm,
        trim={21mm 3mm 3mm 3mm},clip,right]{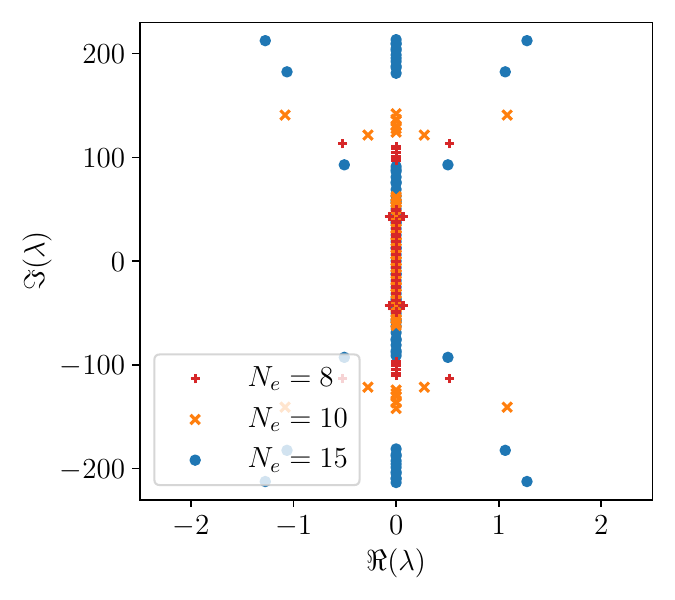}
    \end{subfigure}

    \vspace{0.5em}
    
    \begin{subfigure}[t]{0.272\textwidth}
        \includegraphics[height=2.91cm,
        trim={3mm 3mm 3mm 3mm},clip,left]{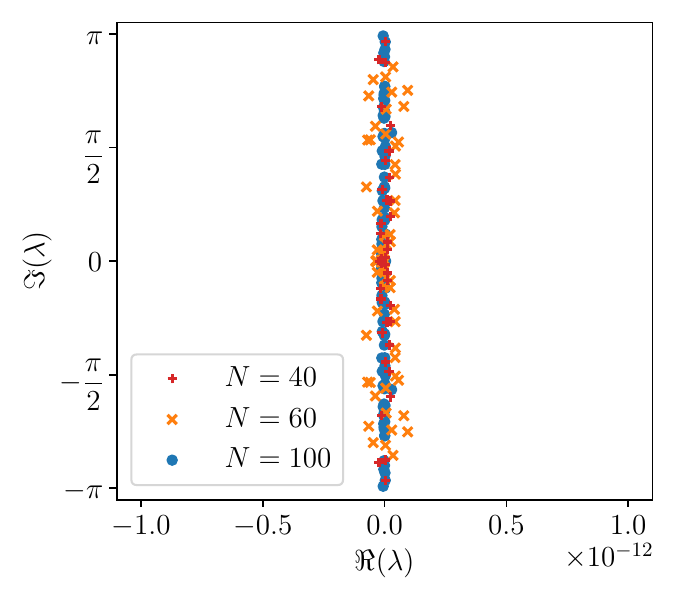}
        \caption{CSBP $p=1$}
        \label{fig:floquet1}
    \end{subfigure}
    \hfill
    \begin{subfigure}[t]{0.23\textwidth}
        \includegraphics[height=2.91cm,
        trim={17mm 3mm 3mm 3mm},clip,center]{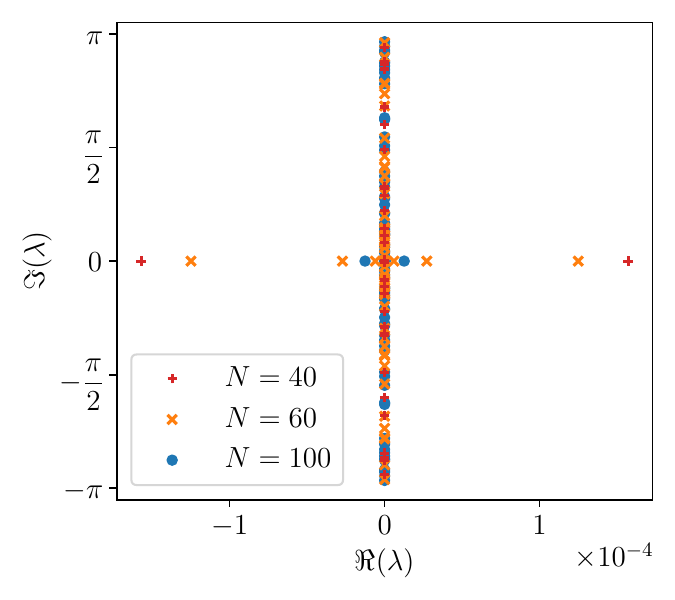}
        \caption{CSBP $p=4$}
        \label{fig:floquet2}
    \end{subfigure}
    \hfill
    \begin{subfigure}[t]{0.23\textwidth}
        \includegraphics[height=2.91cm,
        trim={17mm 3mm 3mm 3mm},clip,center]{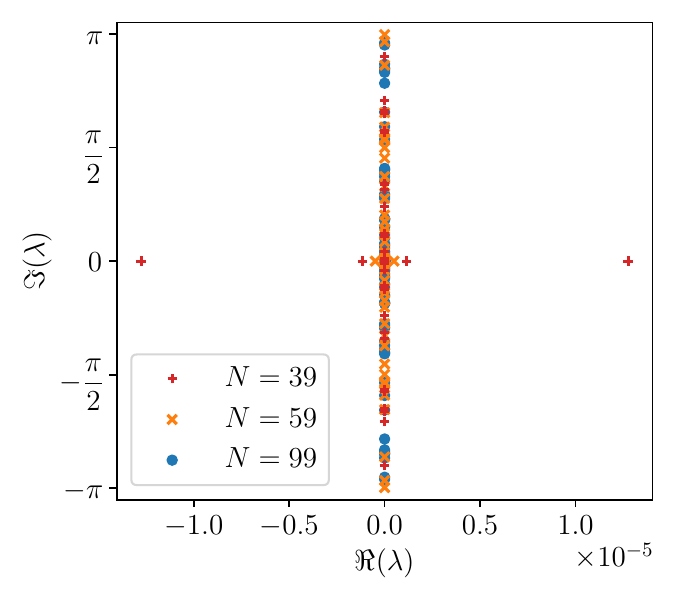}
        \caption{Circulant 8\textsuperscript{th}-order}
        \label{fig:floquet3}
    \end{subfigure}
    \hfill
    \begin{subfigure}[t]{0.235\textwidth}
        \includegraphics[height=2.91cm,
        trim={17mm 3mm 3mm 3mm},clip,right]{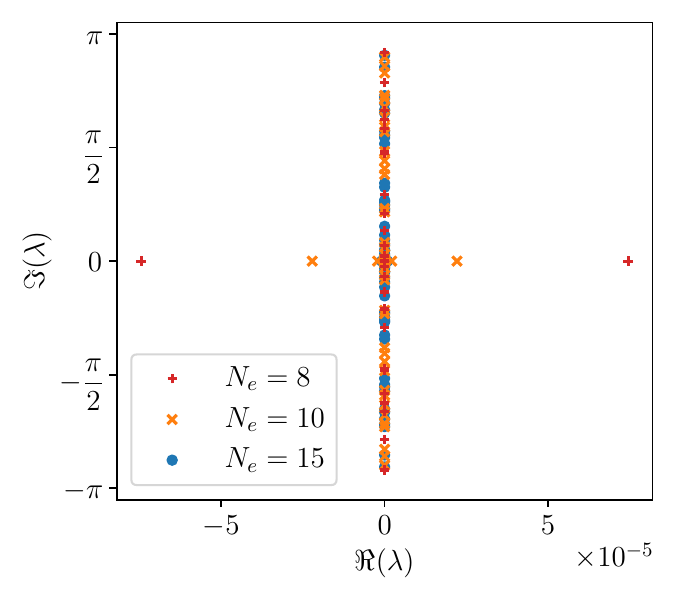}
        \caption{LGL $p=4$}
        \label{fig:floquet4}
    \end{subfigure}

    \caption{Instantaneous eigenspectra at $t=0$ (top row) and Floquet exponents for one period (bottom row) of the entropy-stable geometric flux-differencing scheme applied to the constant-coefficient linear advection equation. No interface or volume dissipation is used. Note the varying $x$-axes in the bottom row.}
    \label{fig:floquet}
\end{figure}

We initially focus on the geometric scheme, since the logarithmic flux can be viewed as a convex combination of the arithmetic and geometric fluxes. We take problem \eqref{eq_lin_conv} with $a(x)=1$, $L=1$, and a Gaussian initial condition of $\fnc{U}_0(x) = e^{0.5 \left( \frac{x-0.5}{0.08} \right)^2} + \sfrac{1}{2}$. To construct the monodromy matrix, we use the exact solution to \eqref{eq_lin_conv} as the baseflow since we wish to isolate perturbation growth, and furthermore, the numerical solution $\bu$ is not exactly periodic in time. In Figure \ref{fig:floquet} we plot the instantaneous eigenspectra and Floquet exponents of the geometric scheme with no volume or interface dissipation at $t=0$, which should both be purely imaginary to indicate convective perturbation transport. Consistent with \cite{gassner_stability_2022, ranocha_2022_preventing_pressure}, the entropy-stable schemes exhibit eigenvalues with large positive real parts. However, this does \emph{not} indicate that perturbations grow unphysically as we saw in \S\ref{sec:time_invariant_accuracy_experiments}, where the Jacobian $\mat{J}$ was time-invariant and its eigenvalues gave an accurate description of modal growth. The Floquet exponents reveal the true story: the perturbation growth rate over a full period is either zero to machine precision (as in Figure~\ref{fig:floquet1}), or a negligibly small value that converges to zero at design-order under mesh refinement (as in Figures~\ref{fig:floquet2}-\ref{fig:floquet4}). As explained in Appendix~\ref{app_sharp_geom_bound}, these small non-zero growth rates are essentially a numerical artifact of using the exact PDE solution instead of the numerical solution as the baseflow in the computation of the monodromy matrix, which introduces a small forcing term due to truncation error that is not present in actual simulations. Floquet analysis therefore reveals that perturbations experiencing instantaneous amplification at a given point in time will soon after be damped due to the rotating dynamics of the problem, precluding any meaningful long-time perturbation growth. This does not entirely preclude perturbation growth in finite time, however. Looking at the singular values of the monodromy matrices, we find $\sigma_{\max}(\mat{\Psi}(T)) \approx 1.7$, telling us that there exist combinations of periodically-stable Floquet modes that experience nearly a doubling in amplitude during a single period, even though their growth is bounded for long times. This is consistent with the sharp perturbation bound of \eqref{eq_geom_pert_bound}, which allows for bounded transient amplification, but precludes exponential modal growth. Similar results are obtained for the logarithmic flux.

\begin{figure}[t!]
    \centering
    \captionsetup[subfigure]{justification=centering,singlelinecheck=false}

    \begin{subfigure}[t]{\textwidth}
        \includegraphics[height=0.5cm,
        trim={5mm 7mm 5mm 7mm},clip,center]{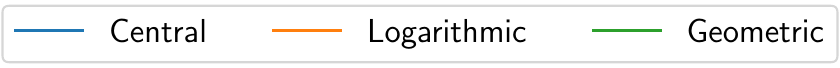}
    \end{subfigure}
    \vspace{-1em}

    \begin{subfigure}[t]{0.255\textwidth}
        \includegraphics[height=2.45cm,
        trim={3mm 3mm 3mm 3mm},clip,left]{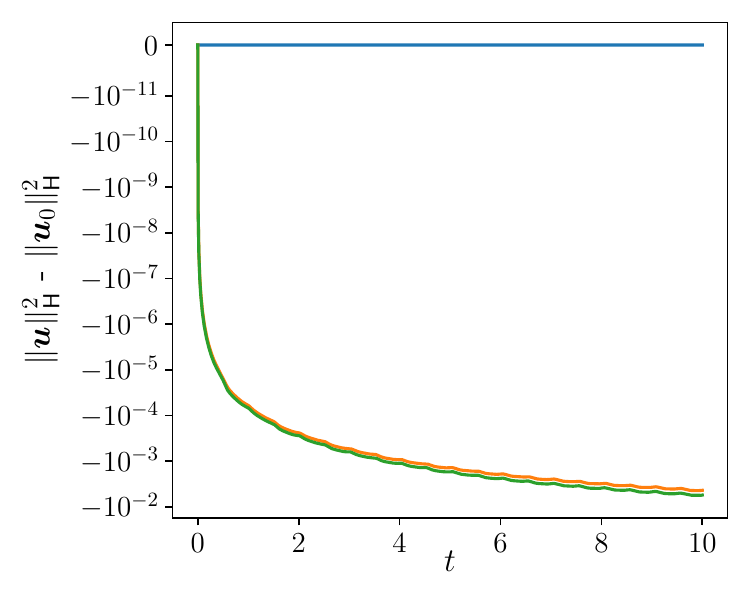}
    \end{subfigure}
    \hfill
    \begin{subfigure}[t]{0.236\textwidth}
        \includegraphics[height=2.45cm,
        trim={12mm 3mm 3mm 3mm},clip,center]{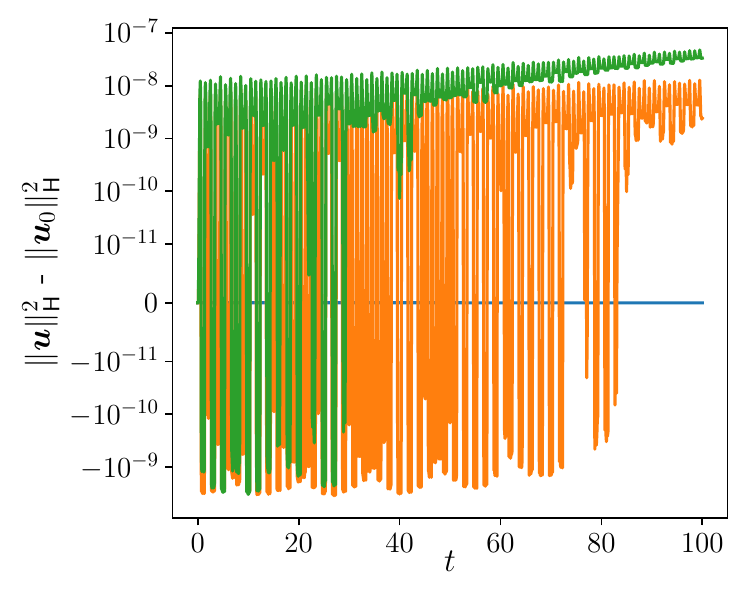}
    \end{subfigure}
    \hfill
        \begin{subfigure}[t]{0.236\textwidth}
        \includegraphics[height=2.45cm,
        trim={11mm 3mm 3mm 3mm},clip,center]{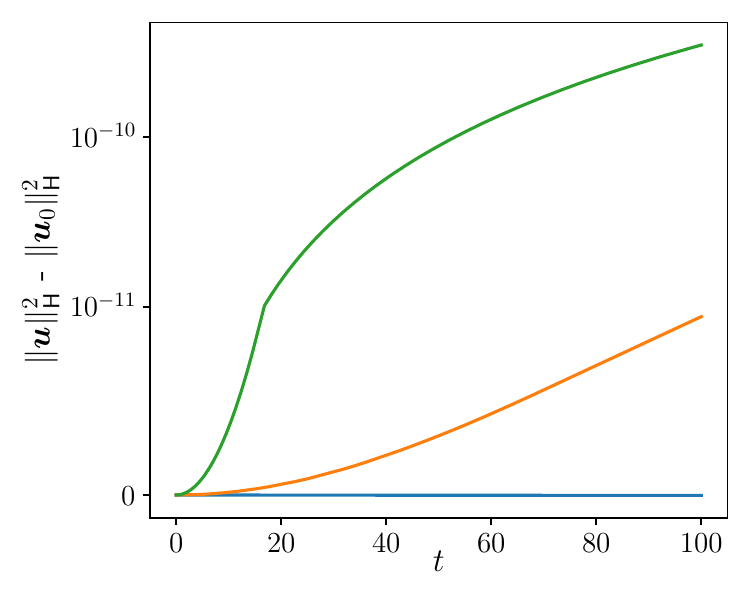}
    \end{subfigure}
    \hfill
        \begin{subfigure}[t]{0.24\textwidth}
        \includegraphics[height=2.45cm,
        trim={11mm 3mm 3mm 3mm},clip,right]{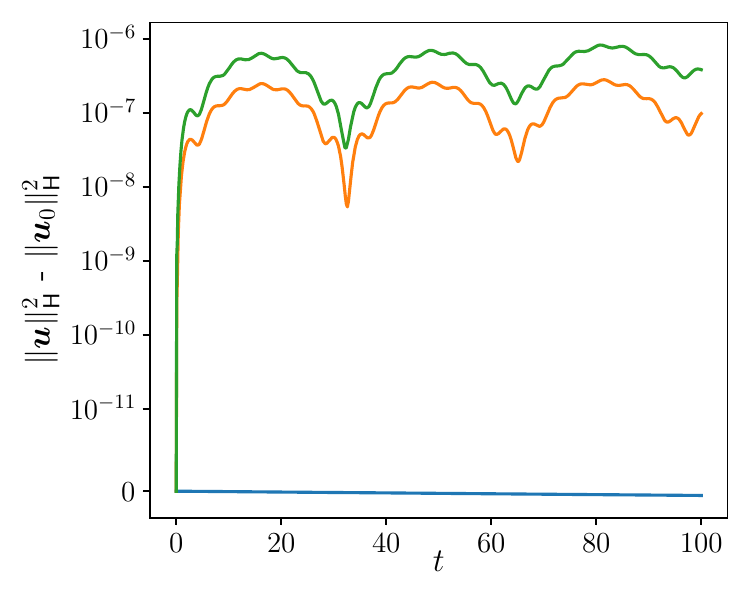}
    \end{subfigure}

    \vspace{0.5em}

    \begin{subfigure}[t]{0.255\textwidth}
        \includegraphics[height=2.45cm,
        trim={3mm 3mm 3mm 3mm},clip,left]{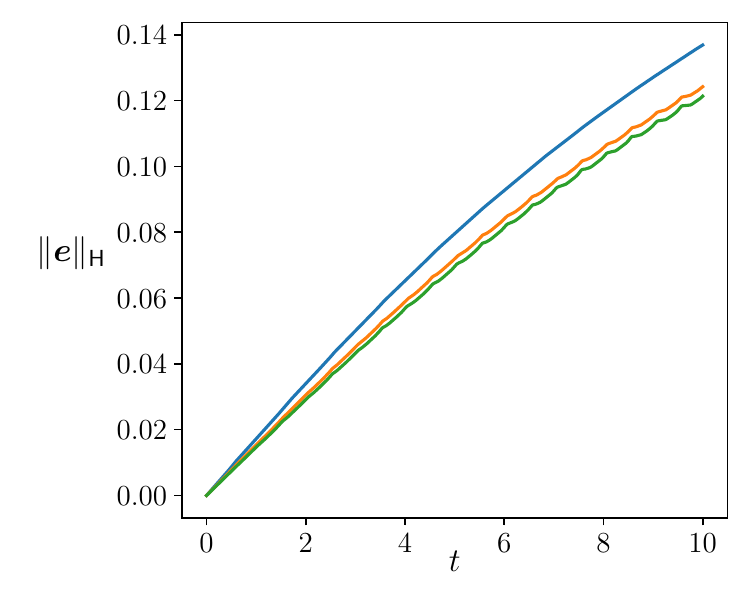}
    \end{subfigure}
    \hfill
    \begin{subfigure}[t]{0.236\textwidth}
        \includegraphics[height=2.45cm,
        trim={16mm 3mm 3mm 3mm},clip,center]{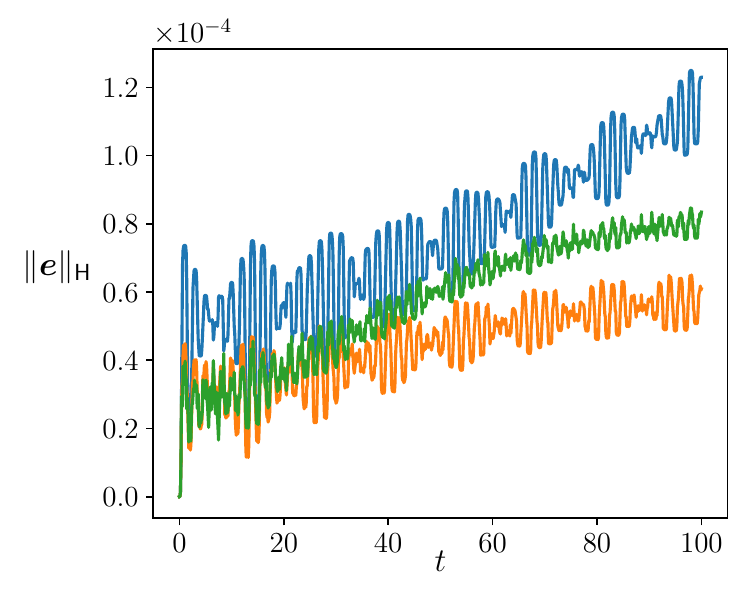}
    \end{subfigure}
    \hfill
        \begin{subfigure}[t]{0.236\textwidth}
        \includegraphics[height=2.45cm,
        trim={16mm 3mm 3mm 3mm},clip,center]{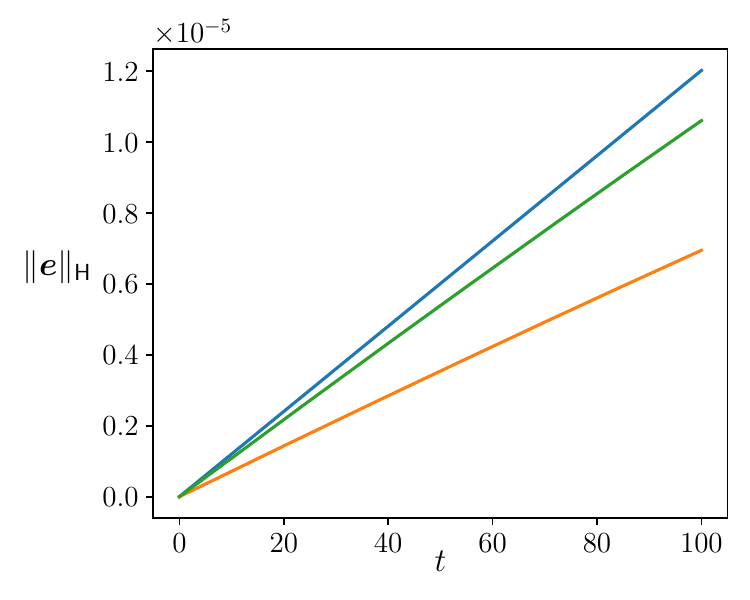}
    \end{subfigure}
    \hfill
        \begin{subfigure}[t]{0.24\textwidth}
        \includegraphics[height=2.45cm,
        trim={16mm 3mm 3mm 3mm},clip,right]{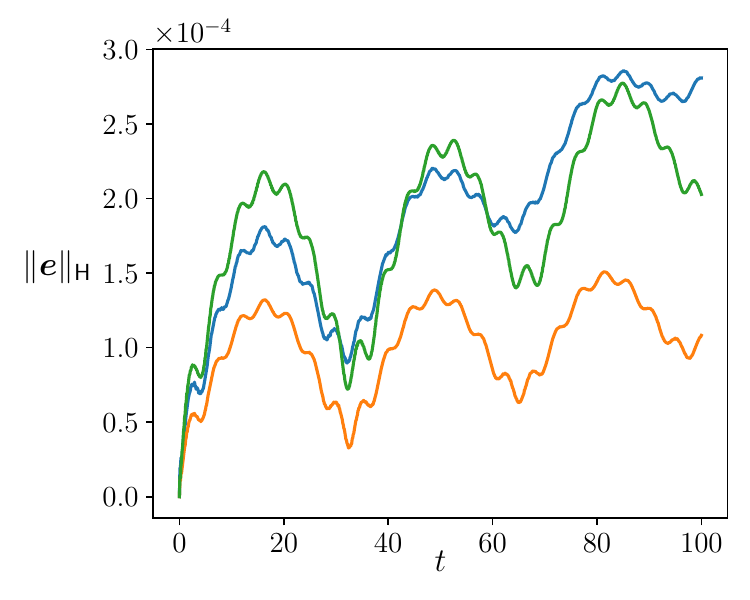}
    \end{subfigure}

    \vspace{0.5em}

    \begin{subfigure}[t]{0.255\textwidth}
        \includegraphics[height=2.45cm,
        trim={3mm 3mm 3mm 3mm},clip,left]{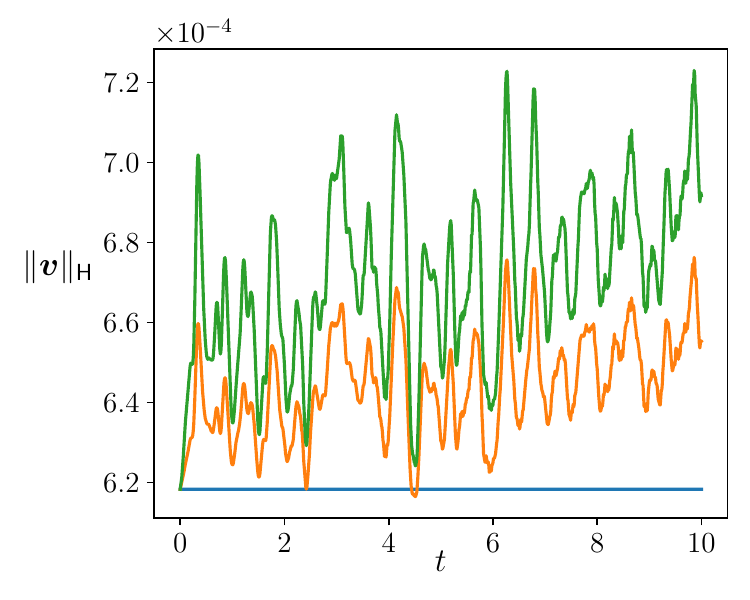}
    \end{subfigure}
    \hfill
    \begin{subfigure}[t]{0.236\textwidth}
        \includegraphics[height=2.45cm,
        trim={16mm 3mm 3mm 3mm},clip,center]{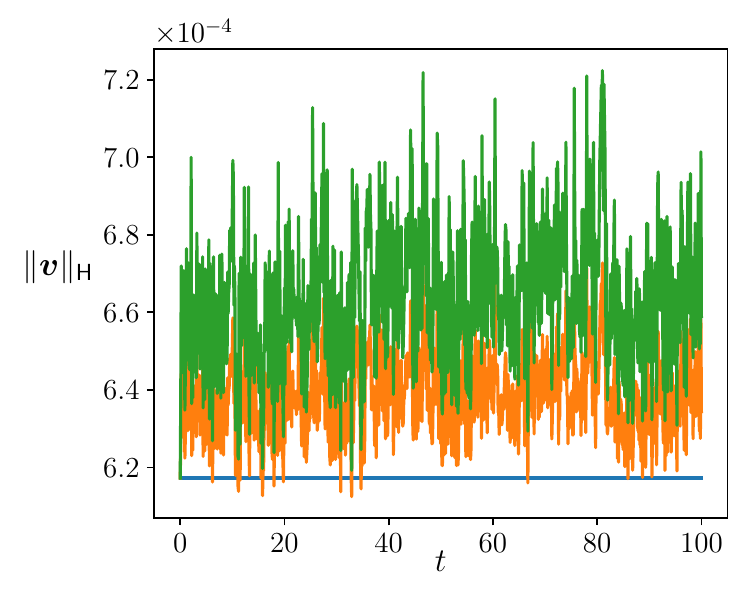}
    \end{subfigure}
    \hfill
        \begin{subfigure}[t]{0.236\textwidth}
        \includegraphics[height=2.45cm,
        trim={16mm 3mm 3mm 3mm},clip,center]{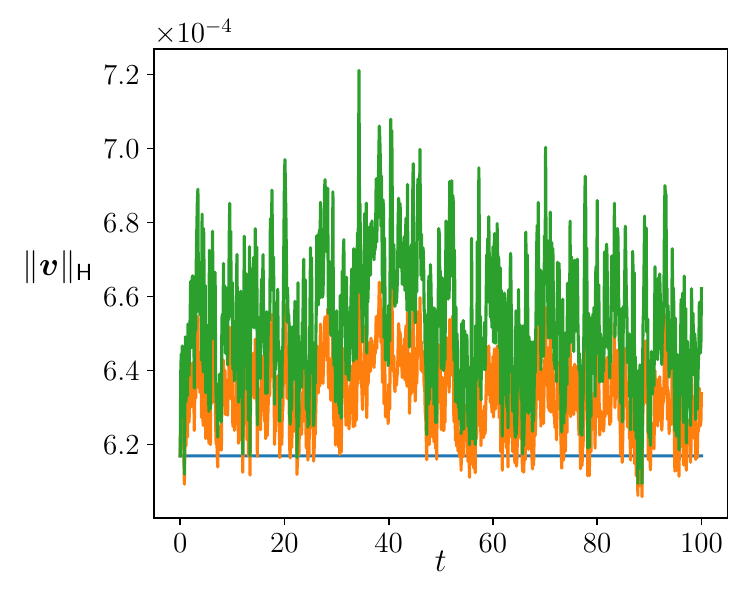}
    \end{subfigure}
    \hfill
        \begin{subfigure}[t]{0.24\textwidth}
        \includegraphics[height=2.45cm,
        trim={16mm 3mm 3mm 3mm},clip,right]{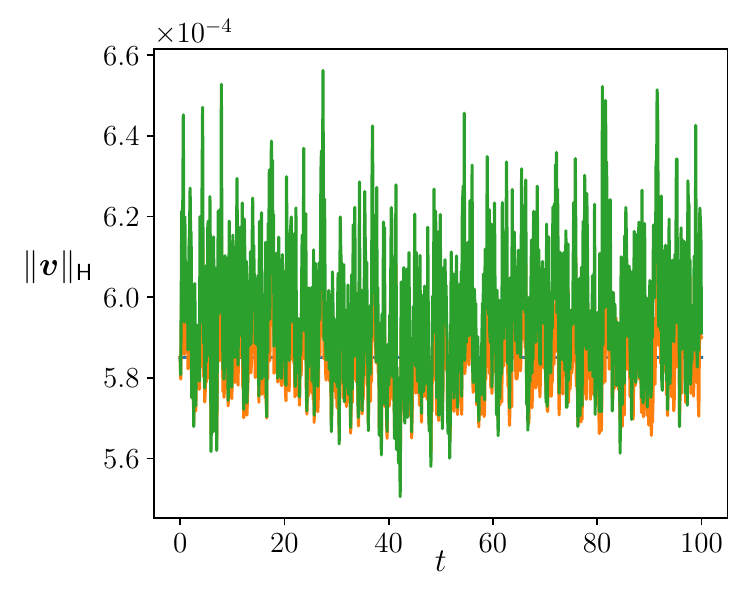}
    \end{subfigure}

    \vspace{0.5em}

    \begin{subfigure}[t]{0.255\textwidth}
        \includegraphics[height=2.45cm,
        trim={3mm 3mm 3mm 3mm},clip,left]{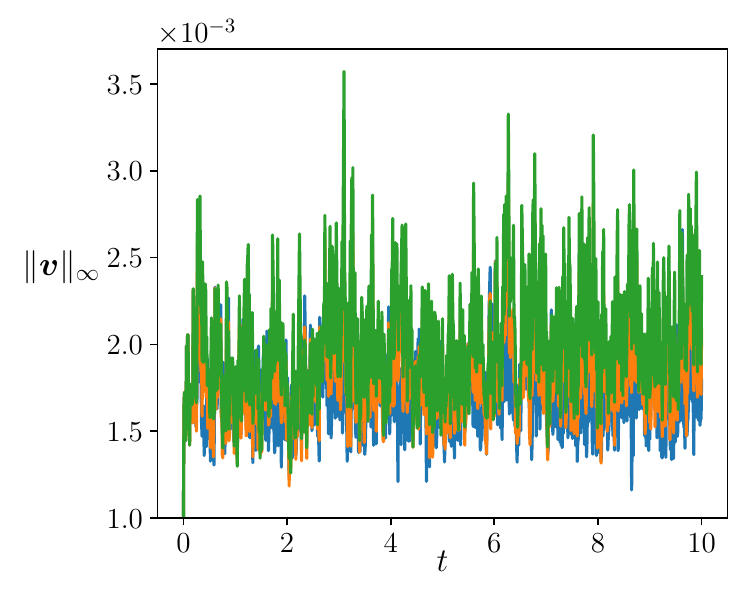}
        \captionsetup{margin={2mm,0mm}}
        \caption{CSBP $p=1$, $N=100$}
        \label{fig:floquetsim1}
    \end{subfigure}
    \hfill
    \begin{subfigure}[t]{0.236\textwidth}
        \includegraphics[height=2.45cm,
        trim={14mm 3mm 3mm 3mm},clip,center]{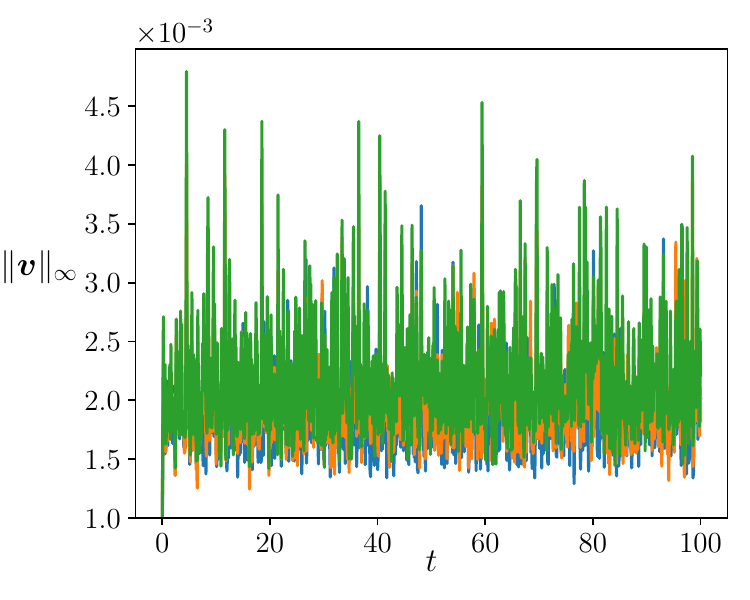}
        \caption{CSBP $p=4$, $N=100$}
        \label{fig:floquetsim2}
    \end{subfigure}
    \hfill
    \begin{subfigure}[t]{0.236\textwidth}
        \includegraphics[height=2.45cm,
        trim={18mm 3mm 3mm 3mm},clip,center]{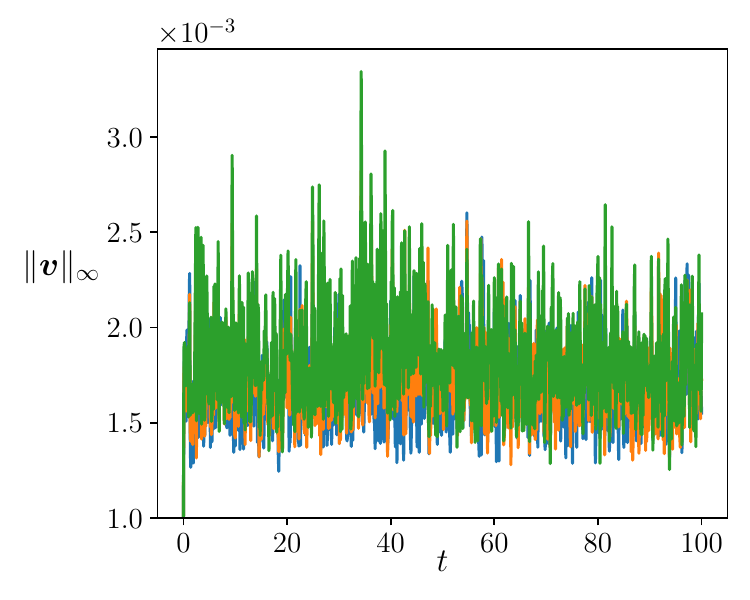}
        \captionsetup{margin={0mm,0mm}}
        \caption{Circulant 8\textsuperscript{th}-order, $N=99$}
        \label{fig:floquetsim3}
    \end{subfigure}
    \hfill
    \begin{subfigure}[t]{0.24\textwidth}
        \includegraphics[height=2.45cm,
        trim={18mm 3mm 3mm 3mm},clip,right]{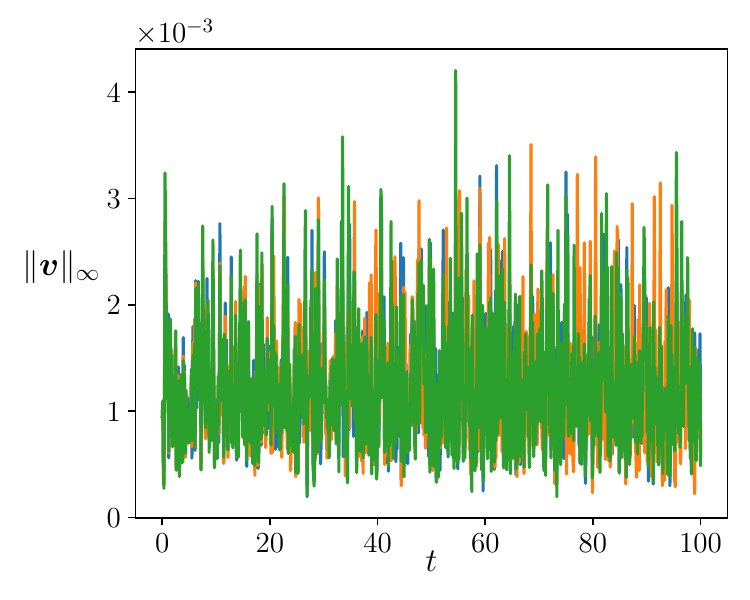}
        \captionsetup{margin={2mm,0mm}}
        \caption{LGL $p=4$, $N_e=20$}
        \label{fig:floquetsim4}
    \end{subfigure}

    \caption{Change in solution energy $\norm{\bu}^2_\Hnrm$ (top row), solution error $\norm{\vec{e}}^2_\Hnrm$ (second row), perturbation $\Hnrm$-norm $\norm{\bv}_\Hnrm$ (third row), and perturbation $L^\infty$-norm $\norm{\bv}_\infty$  (bottom row) of the flux-differencing schemes with central, logarithmic, and geometric fluxes applied to the constant-coefficient linear advection equation with a Gaussian initial condition and a random perturbation added. No interface or volume dissipation is used.}
    \label{fig:floquet_sim}
\end{figure}

We next verify the predictions of Floquet analysis by running the nonlinear flux-differencing schemes for the full problem~\eqref{eq_lin_conv} with $\bu_0 = \fnc{U}_0(\vec{x})$ as described above, then again for a perturbed initial condition $\bu^\epsilon_0 = \fnc{U}_0(\vec{x}) + \bv_0$, and finally isolating the perturbation through $\bv = \bu^\epsilon - \bu$. In this way, the baseflow $\tilde{\bu}$ can be identified as the unperturbed numerical solution $\bu$. We use a random noise initial perturbation $\bv_0 \in [-10^{-3},10^{-3}]^N$, and classical 4\textsuperscript{th}-order Runge--Kutta time marching with a sufficiently small timestep so that temporal errors are negligible. In the bottom two rows of Figure \ref{fig:floquet_sim} we plot perturbation growth as measured in the $\Hnrm$-norm and $L^\infty$ norms. Despite running the simulations for 100 periods, we observe no meaningful perturbation growth for either the geometric or logarithmic schemes, consistent with the Floquet analysis of Figure \ref{fig:floquet} and the bound \eqref{eq_geom_pert_bound}. As expected, the central scheme exactly preserves~$\norm{\bv}_\Hnrm^2$, but this does not indicate that the central scheme is more accurate than the geometric or logarithmic schemes. In fact, the opposite appears to be true. In the first two rows of Figure \ref{fig:floquet_sim} we plot the change in energy $\norm{\bu}_\Hnrm^2$ and solution error $\norm{\vec{e}}_\Hnrm \defn \norm{\bu - \fnc{U}(\vec{x},t)}_\Hnrm$ of the unperturbed simulations (\ie baseflow). The central scheme exhibits larger numerical errors across all tested discretizations. Though we do not claim this result will always hold, it can perhaps be explained by recalling the observation in \S\ref{sec_modified_pde} that the geometric scheme is in fact a central scheme in the square-root variables, and that the logarithmic scheme is a convex combination of the geometric and central schemes.

The Floquet analysis is sensitive to the initial condition $\bu_0$, as this determines the baseflow. In particular, we have found that when an under-resolved initial condition is used (\eg a square wave), the Floquet analysis can predict large Floquet multipliers, or equivalently, Floquet exponents with large positive real parts. This should not be concerning, however, since it is exactly in this regime that we expect the Floquet analysis to be inaccurate. After all, we construct the monodromy matrix using the exact solution $\fnc{U}(\vec{x},t)$ instead of the numerical solution $\bu$. For under-resolved initial conditions, significant error is introduced by the baseflow convection, meaning that the perturbation evolves according to a very different operator than that used in the Floquet analysis. The inaccurate Floquet multipliers then overpredict perturbation growth, as explained further in Appendix~\ref{app_sharp_geom_bound}. Instead, directly running a simulation as in Figure~\ref{fig:floquet_sim} will show negligible perturbation growth in time. 

One may question whether these results contradict those of \cite{ranocha_2022_preventing_pressure}, in which exponential perturbation growth was observed for the logarithmic scheme. Crucially, however, those results were obtained by ``freezing'' the baseflow, \ie by adding a source term to negate the baseflow evolution, allowing only the perturbation to evolve in time. This process is inherently unphysical, as it forces the linear operator $\mat{J}(\bu)$ to be frozen in time. For this nonlinear scheme, one cannot decouple the perturbation evolution from the time-varying baseflow, which we can see even from the modified PDE \eqref{eq_modified_pde_geom_const}. 

A separate (and reasonable) objection to the analysis performed here is that we have used a baseflow with $\min{(\fnc{U})} \approx \tfrac{1}{2}$, while in \eqref{eq_modified_pde_geom_const}, $\nu \rightarrow \infty$ as $\fnc{U} \rightarrow 0$ and $\abs{\partial \fnc{U} / \partial x} \sim 1$, possibly indicating unbounded anti-dissipation. Consequently, local linear instabilities could play a significant role in near-vacuum regions, which is precisely the argument made in \cite{gassner_stability_2022}. We explore this important scenario in the following sections.

\subsection{Issues With Nonlinear Flux-Differencing Schemes At Near-Vacuum States}
\label{sec:nearvacuum}

The previous section demonstrated through numerical experiments that, consistent with \eqref{eq_geom_pert_bound}, local linear instabilities do not lead to perturbation growth in logarithmic and geometric flux-differencing schemes. One may naturally question then why the 1D Euler density-wave problem of \cite{gassner_stability_2022} was shown to be problematic for logarithmic flux-differencing schemes. 

In Appendix \ref{app_nearvacuum}, we explicitly compute the truncation error of the flux-differencing volume term~\eqref{eq_flux_difference} for the logarithmic flux. The key result is that for operators of degree $p$, in near-vacuum regions where $\min_x \fnc{U} \rightarrow 0$, the truncation error contains terms that scale proportionally to $\fnc{U}^{-p}$. This implies that the flux-differencing volume term \eqref{eq_flux_difference} becomes extremely inaccurate as the solution enters near-vacuum regions, particularly for high-order discretizations. This result is stated in terms of a fixed wave profile of $\fnc{U}$ being shifted downwards towards $0$ such that the problem can not be avoided through amplitude rescalings: we consider a small solution in comparison with its local spatial variation. A similar issue can also be identified with respect to the ill-conditioning and nonnormality of the semidiscrete linearized Jacobians~$\mat{J}(t)$. Appendix \ref{app_nearvacuum} also demonstrates that as the numerical solution enters this regime, perturbation growth can become unbounded as measured through 
\[
\lambda_{\max}(\mat{J}_\mathrm{sym}) \rightarrow \infty \ , \qquad \mat{J}_\mathrm{sym} \defn \tfrac{1}{2} \left( \mat{J} + \mat{J}^\T \right) .
\]
This can be concerning because the worst-case instantaneous growth, often driven by transient growth, is controlled by the largest positive eigenvalues of the symmetric part of~$\mat{J}(\bu(t))$. To see this, left multiply \eqref{eq_time_periodic_system} by~$\bv^\T$ to obtain $\mathrm{d}/\mathrm{d}t \norm{\bv}^2 = \bv^\T \mat{J}_\mathrm{sym} \bv$. Furthermore, the spectral radius of~$\mat{J}$ can also increase unboundedly in this regime, \ie~$\rho(\mat{J}) \rightarrow \infty$. This introduces severe timestep stability restrictions for explicit time marching. Though these results do not contradict the perturbation bound \eqref{eq_geom_pert_bound}, as the bound itself becomes singular in this regime, they nonetheless have potentially catastrophic implications for logarithmic and geometric flux-differencing simulations attempting to simulate the near-vacuum regime.

To observe the effects of entering this problematic regime, we again employ a logarithmic flux-differencing discretization for the constant-coefficient linear advection equation with $a(x)=1$, now with an initial condition mimicking the density-wave problem in \cite{gassner_stability_2022},
\[
\fnc{U}_0(x) = 0.98 \sin(2 \pi x) + 1 , \quad x \in [-1,1].
\]
We perform an initial Floquet analysis as described in \S\ref{sec_logarithmic_LCE_experiments} and initiate the perturbation using the most unstable Floquet mode to maximize potential perturbation growth. To isolate the near-vacuum error mechanism in the simplest possible setting, we consider only circulant operators. This removes any impact of modified boundary stencils, noting that even in the worst case, \S\ref{sec:time_invariant_accuracy_experiments} demonstrated that the addition of interface and volume dissipation are sufficient to reduce the problem to this simpler scenario. Since the sensitivity of the truncation error to near-vacuum increases with degree (unlike the local linear instabilities explored in \S\ref{sec:spectra}), we use 8\textsuperscript{th}-order operators and $N=39$. We optionally include a small amount of volume dissipation (see \cite{dissipation}) applied to either the conservative variable $\bu$ or the entropy variable $\log(\bu)$, noting that volume dissipation applied to entropy variables has beneficial positivity-preservation properties~\cite{dissipation}. Alternatively, we also include a crude positivity-preserving filter, applied once per timestep after the final 4\textsuperscript{th}-order Runge--Kutta update, given by
\begin{equation} \label{eq_positivity_filter}
u_i \rightarrow u_i + \varepsilon \log\!\left(1+e^{(u_{\mathrm{floor}}-u_i)/\varepsilon)} \right), \quad \varepsilon = \tfrac{1}{12} \left( u_\mathrm{cut} - u_\mathrm{floor} \right) ,
\end{equation}
where we set $u_\textrm{floor} = 5 \times 10^{-4}$ and $u_\textrm{cut} = 5 \times 10^{-2}$. This filter indiscriminately increases any values of $u_i$ below the approximate threshold set by $u_\text{cut}$ to a value above $u_\text{floor}$. It is neither conservative nor entropy-dissipative. Furthermore, it does not preserve design-order. Therefore, in the absence of the near-vacuum error mechanism described above, one would expect it to increase numerical error.

\begin{figure}[t]
    \centering
    \captionsetup[subfigure]{justification=centering,singlelinecheck=false}

    \begin{subfigure}[t]{\textwidth}
        \includegraphics[height=0.5cm,
        trim={5mm 7mm 5mm 7mm},clip,center]{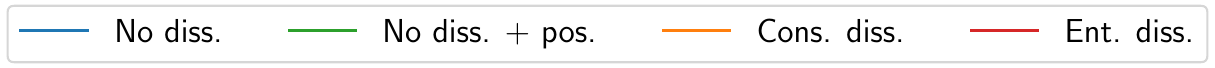}
    \end{subfigure}
    \vspace{-1em}

    \begin{subfigure}[t]{0.32\textwidth}
        \includegraphics[width=\textwidth,
        trim={3mm 3mm 3mm 3mm},clip]{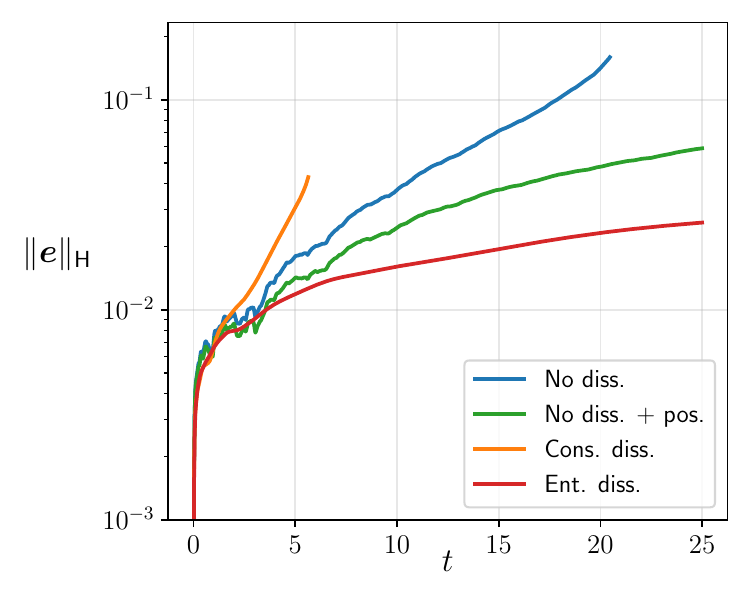}
        \caption{Baseflow error $\vec{e} \defn \bu - \fnc{U}(\vec{x})$}
        \label{fig:densitywave_near0a}
    \end{subfigure}
    \hfill
    \begin{subfigure}[t]{0.32\textwidth}
        \includegraphics[width=\textwidth,
        trim={3mm 3mm 3mm 3mm},clip]{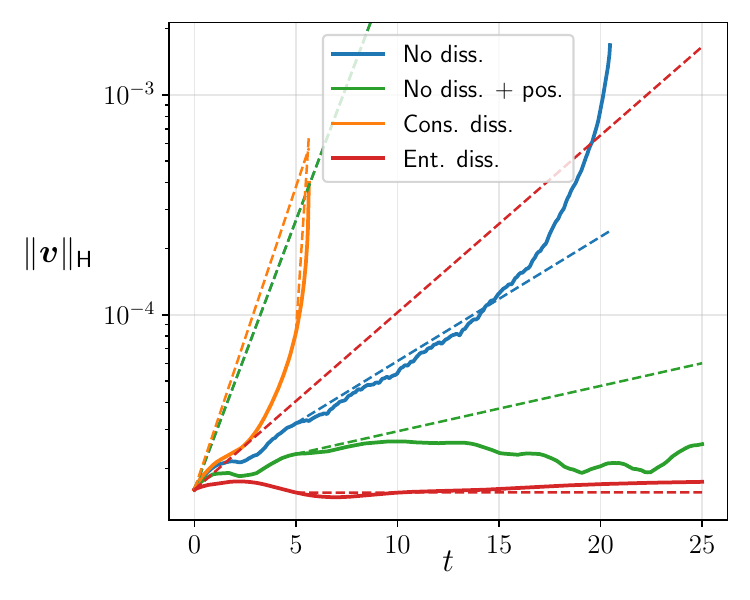}
        \caption{Perturbation growth and predictions at $t=0, 5$}
        \label{fig:densitywave_near0b}
    \end{subfigure}
    \hfill
        \begin{subfigure}[t]{0.32\textwidth}
        \includegraphics[width=\textwidth,
        trim={3mm 3mm 3mm 3mm},clip]{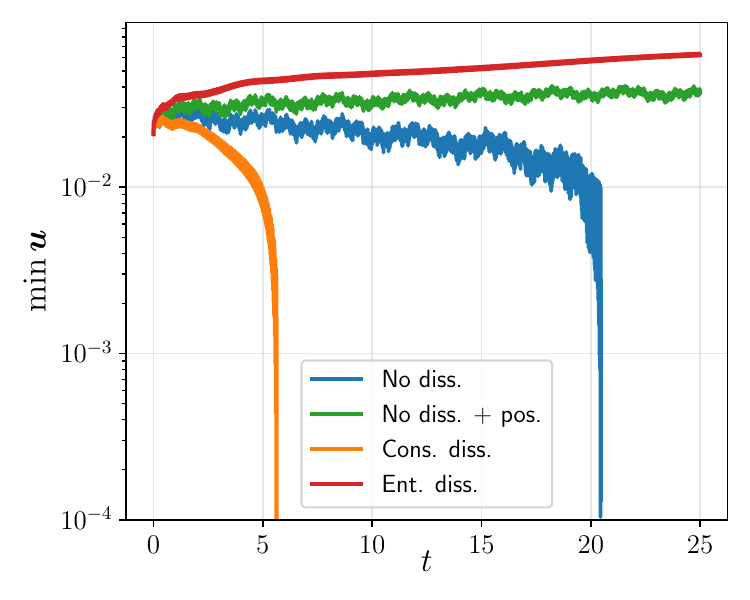}
        \caption{Baseflow minimum solution value $\bu$}
        \label{fig:densitywave_near0c}
    \end{subfigure}

    \vspace{0.5em}
    
    \begin{subfigure}[t]{0.32\textwidth}
        \includegraphics[width=\textwidth,
        trim={3mm 3mm 3mm 3mm},clip]{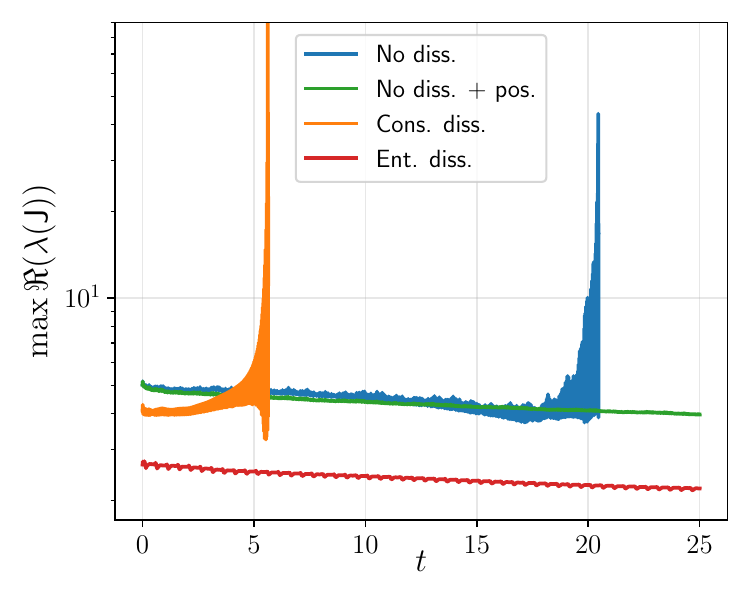}
        \caption{Max. instantaneous eigenvalues}
        \label{fig:densitywave_near0d}
    \end{subfigure}
    \hfill
    \begin{subfigure}[t]{0.32\textwidth}
        \includegraphics[width=\textwidth,
        trim={3mm 3mm 3mm 3mm},clip]{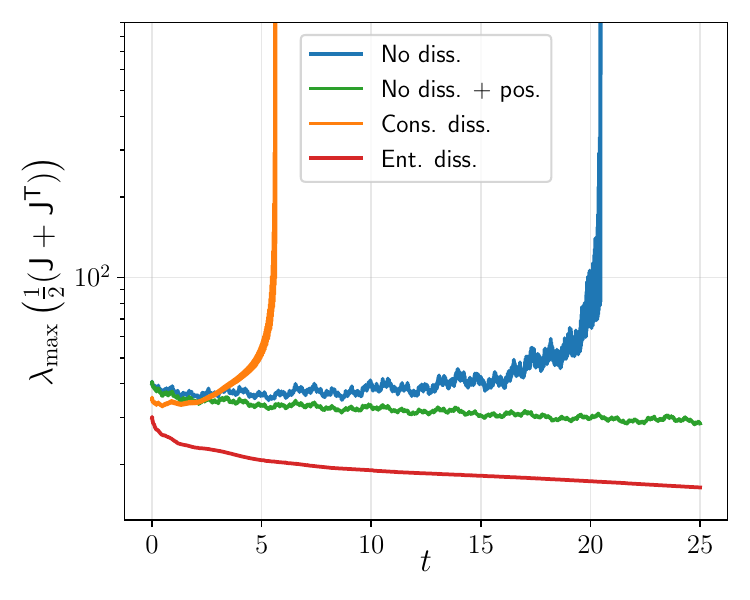}
        \caption{Max. instantaneous growth}
        \label{fig:densitywave_near0e}
    \end{subfigure}
    \hfill
        \begin{subfigure}[t]{0.32\textwidth}
        \includegraphics[width=\textwidth,
        trim={3mm 3mm 3mm 3mm},clip]{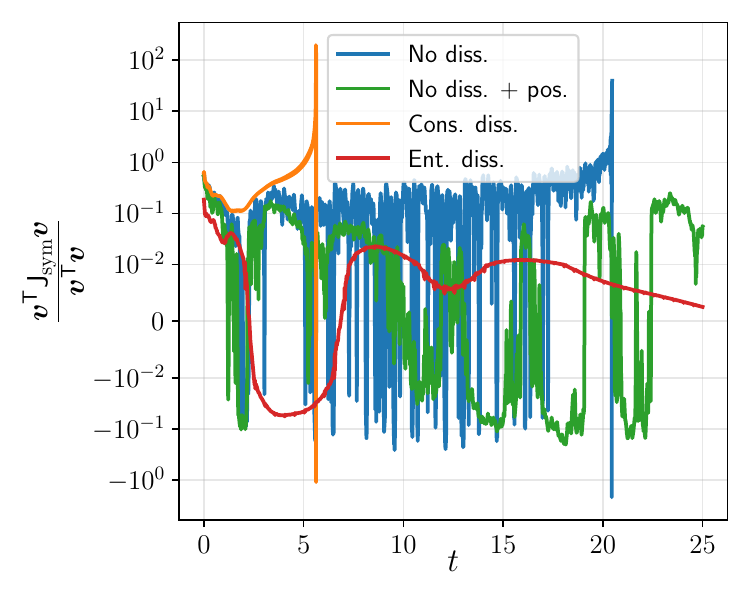}
        \caption{Instantaneous pert. growth}
        \label{fig:densitywave_near0f}
    \end{subfigure}

    \caption{A logarithmic flux-differencing discretization of a near-vacuum linear advection problem using 8\textsuperscript{th}-order circulant operators, $N=39$, and optionally including volume dissipation applied to either the conservative or entropy variable, or a crude positivity-preserving filter.}
    \label{fig:densitywave_near0}
\end{figure}

In Figure \ref{fig:densitywave_near0a}, the numerical error introduced in the baseflow convection is shown to be significant. For comparison (though we do not plot it here), in this near-vacuum regime, a central scheme maintains between three and four orders of magnitude smaller error than the logarithmic flux-differencing scheme. The numerical error increases rapidly near the beginning of the simulation as the initial condition projects onto the discrete solution space. Consequently, initial Floquet analyses are inaccurate. Shown in Figure~\ref{fig:densitywave_near0b} as dotted lines, the initial growth rates computed via Floquet analysis overpredict perturbation growth, whereas a second analysis performed at $t=5$ using the current numerical solution is shown to be much more accurate. 
Once again, we see no meaningful perturbation growth from the two discretizations that include some positivity preservation mechanism, either entropy-based dissipation or a crude positivity filter. In contrast, the two discretizations that are allowed to enter the near-vacuum regime experience some perturbation growth, though it is still much smaller than what one would expect if it were caused by instantaneous modal growth (\ie local linear instabilities), corresponding to the growth rates shown in Figure~\ref{fig:densitywave_near0d}. In Figure~\ref{fig:densitywave_near0c}, we observe that for the case with no dissipation and the case with conservative-variable dissipation, entering the near-vacuum regime instigates a series of catastrophic errors, including ill-conditioning of the system (seen in Figure~\ref{fig:densitywave_near0d}) and nonnormality (seen in Figure~\ref{fig:densitywave_near0e}). On the other hand, when the solution is prevented from approaching vacuum, $\mat{J}$ remains well-behaved, perturbations do not grow, and the overall error is decreased. This is particularly striking for the crude positivity filter, as any error introduced by this filter is seemingly compensated by a reduction in error from the inaccurate logarithmic mean in the near-vacuum region. These results demonstrate that any perturbation growth observed when running density-wave problems of~\cite{gassner_stability_2022, ranocha_2022_preventing_pressure} is not a result of local linear instabilities, but rather a side effect of the solution entering the problematic near-vacuum regime.

Without resorting to positivity preservation, it is still possible to push back the crash time through mesh refinement. This decreases the error introduced in the baseflow, and delays the time at which the solution enters the near-vacuum regime. Though we do not include the results here, we have also explored constructing flux-differencing schemes based on nonlinear means that avoid the problematic singularity at $\fnc{U}=0$, such as an inverse hyperbolic sine mean. Similar experiments demonstrate that these schemes exhibit much lower numerical error (\ie more consistent with a central discretization) and negligible perturbation growth, despite exhibiting eigenvalues with positive real parts.

\section{Nonlinear Equations} \label{sec:nonlinear}

So far, this paper has considered only linear schemes through analysis of the one-dimensional linear advection equation \eqref{eq_lin_conv}. The results of previous sections, however, extend directly to nonlinear equations of interest, including the inviscid Burgers and compressible Euler equations.

\subsection{The Inviscid Burgers Equation} \label{sec_burgers}

The variable-coefficient linear advection equation \eqref{eq_lin_conv} is directly relevant to the inviscid Burgers equation because it models the growth of perturbations in the nonlinear problem~\cite{kreiss_lorenz, gassner_stability_2022}. For completeness, this standard argument is repeated in Appendix~\ref{app_lin_burgers}.

Consider a split-form discretization of the Burgers equation with circulant operators (interface SATs are handled in Appendix~\ref{app_lin_burgers}) and periodic boundary conditions~\cite{DelReyFernandez2014_review},
\begin{equation} \label{burgers_disc_circulant}
    \der[\bu]{t} + \alpha \mat{D} \left( \tfrac{1}{2} \mat{U} \bu \right) + (1 - \alpha ) \mat{U} \mat{D} \bu   = \vec{0} , 
    \quad \mat{U} \defn \diag{\bu} .
\end{equation} 
Convergence of the entropy-conserving (or equivalently, energy-conserving) scheme with $\alpha = \sfrac{2}{3}$ can be established without considering the linearized problem~\cite{worku2026convergence}. At worst, therefore, local linear instabilities can only be a mechanism for introducing error. Since the $\alpha = \sfrac{2}{3}$ scheme has a nonlinear stability bound, once the error from linear instabilities saturates the numerical discretization, it can no longer grow as we saw in \S\ref{sec:time_invariant_accuracy_experiments}~\cite{gassner_stability_2022}. Regardless, to ensure accurate simulations, before they get a chance to grow, these instabilities should be sufficiently suppressed. We have already shown in \S\ref{sec:time_invariant_accuracy_experiments} that even if we ignore the rotating dynamics of a variable baseflow (which may naturally suppress any perturbation growth), these instabilities are effectively controlled by including interface and volume dissipation as needed. Furthermore, the inviscid Burgers equation develops shocks within a short time, far shorter than any time scale on which we observed significant perturbation growth in \S\ref{sec:time_invariant_accuracy_experiments}. Beyond this point, linear analysis cannot be used, and the arguments of local linear stability do not apply.

In the context of the nonlinear scheme, we can also question the `central plus dissipation' interpretation. Since the $\alpha=\sfrac{2}{3}$ split-form conserves energy while the central discretization $\alpha=1$ does not, in contrast to the linear case, there is a strong case to be made that the $\alpha=\sfrac{2}{3}$ split-form is the natural baseline scheme in this context. Just as we did in \S\ref{sec_modified_pde}, using 2\textsuperscript{nd}-order central differencing, one can write any splitting of \eqref{burgers_disc_circulant} as an entropy-conservative discretization plus some correction that looks like (anti-)dissipation,
\[
\der[u_i]{t} + \frac{f^\text{EC}_{i+1/2} - f^\text{EC}_{i-1/2}}{h} = - \frac{h^2}{2} \left(\alpha - \frac{2}{3} \right) \frac{(u_{i+1}-u_{i-1})}{2h} \frac{(u_{i+1}-2 u_i + u_{i-1})}{h^2} ,
\]
where $f^\text{EC}_{i+1/2} \defn \tfrac{1}{6}(u_i^2 + u_i u_{i+1} + u_{i+1}^2)$. This can be interpreted as a consistent discretization of the modified PDE
\[
\pder[\fnc{U}]{t} + \pder{x} \left( \frac{1}{2} \fnc{U}^2 \right) = \pder{x} \left( \nu(x,t) \pder[\fnc{U}]{x} \right) , \quad \nu(x,t) = - h^2 \frac{\left(\alpha - \tfrac{2}{3} \right)}{4} \pder[\fnc{U}]{x} .
\]
where for $\alpha = 1$, energy can grow unphysically in regions of increasing $\fnc{U}$. Therefore, although care should be taken to ensure that the error introduced by local linear instabilities is sufficiently controlled, our results suggest that their presence alone does not provide a compelling reason to favour the central scheme over the entropy-stable split form of the Burgers equation.

\subsection{The Compressible Euler Equations} \label{sec:Euler}

Consider the one-dimensional compressible Euler equations,
\begin{equation} \label{eq_euler}
    \pder[\vecfnc{U}]{t} + \pder[\vecfnc{F}]{x} = \vec{0}, \quad 
    \vecfnc{U} \defn (\rho, \rho \mathrm{v}, e), \quad
    \vecfnc{F} \defn (\rho \mathrm{v}, \rho \mathrm{v}^2 + p, (e+p) \mathrm{v}),
\end{equation}
where $\rho$ is the density, $\mathrm{v}$ is the velocity, $e$ is the total energy, and the pressure is determined via the ideal gas law, $p = (\gamma-1) (e - \tfrac{1}{2} \rho \mathrm{v}^2)$. We focus our study here on the one-dimensional density-wave problem from \cite{gassner_stability_2022}, 
\begin{equation} \label{eq_euler_densitywave}
    \begin{aligned}
        \rho(x,t=0) = 0.98 \sin(2\pi x) + 1, \quad
        \mathrm{v} = 0.1, \quad
        p = 20, \quad 
        x \in [-1,1],
    \end{aligned}
\end{equation}
which reduces the Euler equations to a decoupled system of linear advection equations similar to the problem studied in~\S\ref{sec:nearvacuum}. In recent publications (\eg \cite{gassner_stability_2022,ranocha_2022_preventing_pressure,dissipation}), it has been argued that the failure of entropy-conservative schemes to robustly solve this problem is due to local linear instabilities. That is, the presence of eigenvalues with positive real parts in the linearized semidiscrete Jacobian indicates that perturbations grow unphysically, eventually leading to spurious oscillations that drive the solution to violate positivity constraints. The results of \S\ref{sec_logarithmic_LCE_experiments} and \S\ref{sec:nearvacuum} however suggest that on the contrary, eigenvalues of the linearized Jacobian with positive real parts do not indicate unphysical perturbation growth. The failure of entropy-conservative schemes is instead primarily due to the inaccuracy of the logarithmic mean at near-vacuum states.

To assess the conclusions of~\S\ref{sec_logarithmic_LCE_experiments} and~\S\ref{sec:nearvacuum} in the context of the Euler equations, we repeat the numerical experiment of~\S\ref{sec:nearvacuum} for the density-wave problem. A simulation is performed using the initial condition \eqref{eq_euler_densitywave}, then again with a small added perturbation of the most long-term unstable Floquet mode, normalized to an initial amplitude of~$10^{-5}$. We employ the entropy-conserving two-point flux of Ranocha~\cite{ranocha_2018}, 8\textsuperscript{th}-order circulant operators with $N=39$, 4\textsuperscript{th}-order Runge--Kutta time marching with $\Delta t = 10^{-4}$, and optionally include matrix-based volume dissipation applied to either the conservative or entropy variables using the framework of \cite{dissipation}. We also optionally include the crude positivity-preserving filter \eqref{eq_positivity_filter} applied to both~$\rho$ and~$p$.

\begin{figure}[t]
    \centering

    \begin{subfigure}[t]{\textwidth}
        \includegraphics[height=0.5cm,
        trim={5mm 7mm 5mm 7mm},clip,center]{1dEuler_densitywave_near0_legend2.pdf}
    \end{subfigure}
    \vspace{-1em}

    \begin{subfigure}[t]{0.32\textwidth}
        \includegraphics[width=\textwidth,
        trim={3mm 3mm 3mm 3mm},clip]{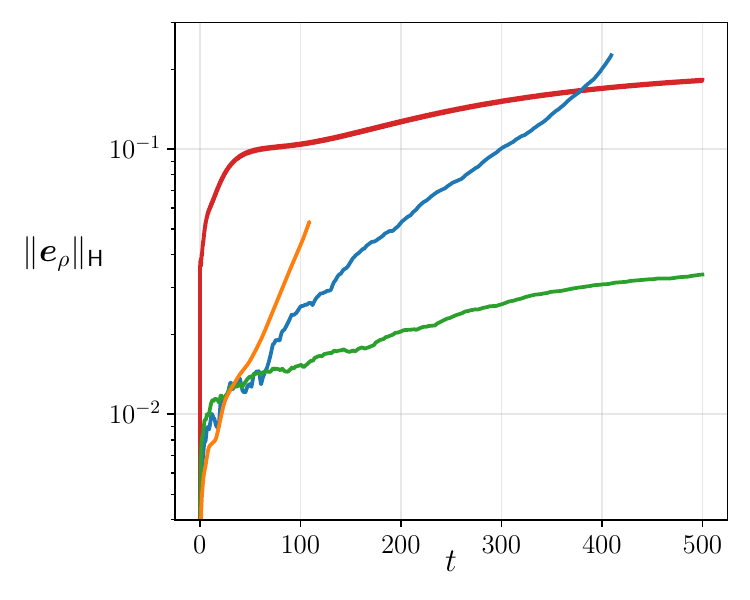}
        \caption{Baseflow $\rho$ error}
        \label{fig:Euler_densitywave_near0a}
    \end{subfigure}
    \hfill
    \begin{subfigure}[t]{0.32\textwidth}
        \includegraphics[width=\textwidth,
        trim={3mm 3mm 3mm 3mm},clip]{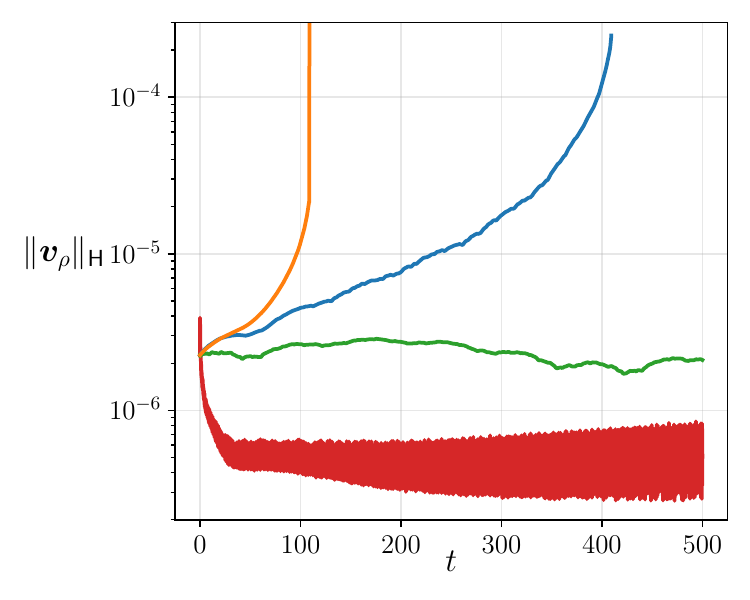}
        \caption{Perturbation $\rho$ growth}
        \label{fig:Euler_densitywave_near0b}
    \end{subfigure}
    \hfill
    \begin{subfigure}[t]{0.32\textwidth}
        \includegraphics[width=\textwidth,
        trim={3mm 3mm 3mm 3mm},clip]{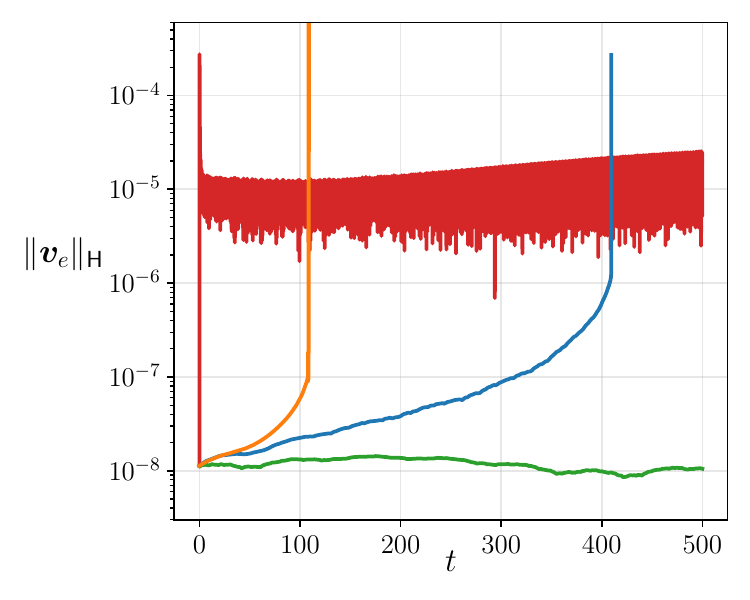}
        \caption{Perturbation $e$ growth}
        \label{fig:Euler_densitywave_near0c}
    \end{subfigure}

    \vspace{0.5em}
    
     \begin{subfigure}[t]{0.32\textwidth}
        \includegraphics[width=\textwidth,
        trim={3mm 3mm 3mm 3mm},clip]{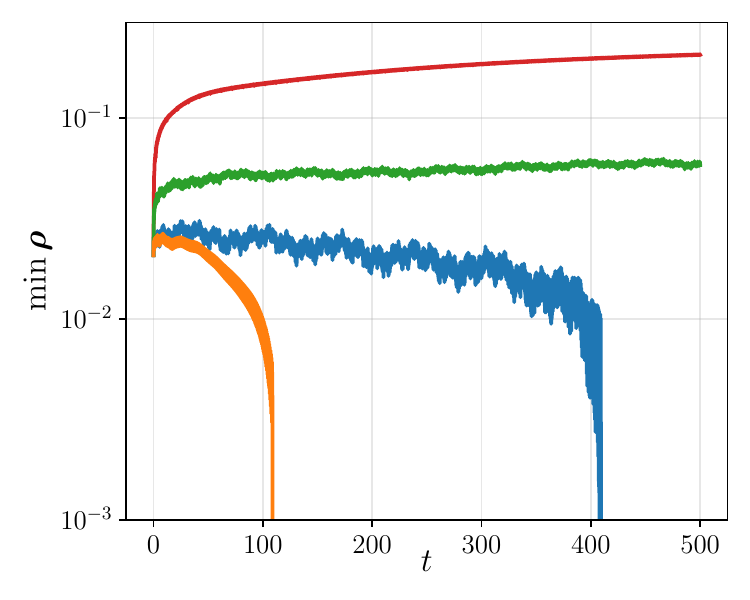}
        \caption{Baseflow minimum $\rho$}
        \label{fig:Euler_densitywave_near0d}
    \end{subfigure}
    \hfill
    \begin{subfigure}[t]{0.32\textwidth}
        \includegraphics[width=\textwidth,
        trim={3mm 3mm 3mm 3mm},clip]{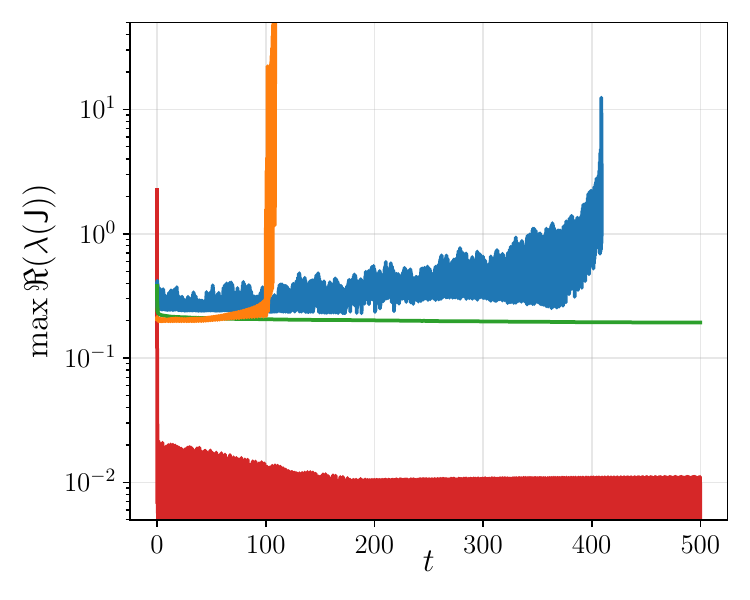}
        \caption{Max. inst. eigenvalues}
        \label{fig:Euler_densitywave_near0e}
    \end{subfigure}
    \hfill
    \begin{subfigure}[t]{0.32\textwidth}
        \includegraphics[width=\textwidth,
        trim={3mm 3mm 3mm 3mm},clip]{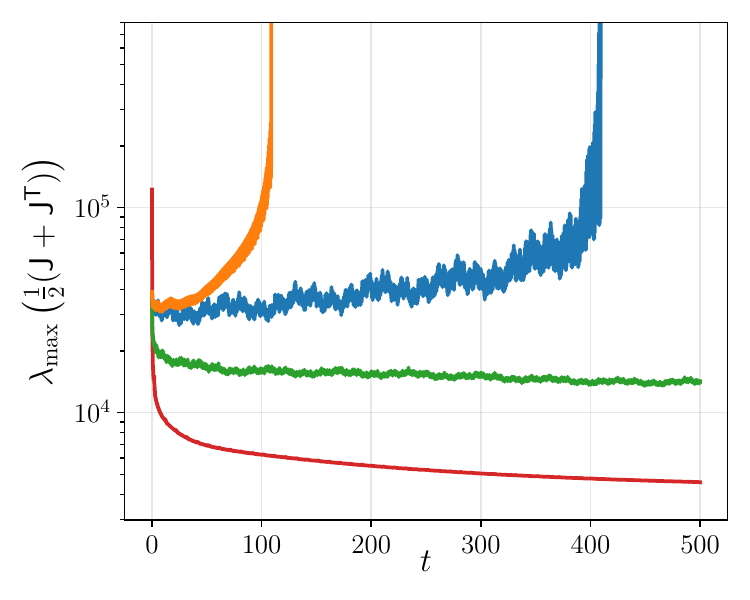}
        \caption{Max. inst. growth}
        \label{fig:Euler_densitywave_near0f}
    \end{subfigure}

    \caption{The 1D Euler density-wave problem using the entropy-conserving Ranocha two-point flux, 8\textsuperscript{th}-order circulant operators, $N=39$, and optionally including volume dissipation applied to either the conservative or entropy variable, or a crude positivity-preserving filter.}
    \label{fig:Euler_densitywave_near0}
\end{figure}

Figure \ref{fig:Euler_densitywave_near0a} shows the numerical error introduced in the baseflow convection to be significant. For comparison, the density error of flux-differencing schemes relying on arithmetic means (such as the central scheme) is much smaller for this problem, lying between $10^{-6}$ and $10^{-4}$. Because the entropy approaches non-convexity and the entropy variables become cusp-like at near vacuum-states~\cite{chan_entropyprojection,dissipation}, the entropy-based dissipation becomes excessively dissipative for this problem, introducing significant numerical error. Consistent with the results of \S\ref{sec:nearvacuum}, however, the crude positivity-preserving filter reduces the overall numerical error by preventing the solution from approaching vacuum where the logarithmic means become inaccurate. For the two schemes without positivity-preserving mechanisms, as the solution approaches the near-vacuum regime in Figure \ref{fig:Euler_densitywave_near0d}, the same catastrophic errors observed in \S\ref{sec:nearvacuum} appear. The numerical baseflow error grows rapidly (Figure \ref{fig:densitywave_near0a}), and both the eigenvalues of the linearized Jacobian (Figure \ref{fig:densitywave_near0e}) and the worst-case instantaneous perturbation growth (Figure \ref{fig:densitywave_near0f}) blow up. The latter, driven by nonnormality of the linearized operator, leads to significant transient growth that also causes the perturbations to blow up (Figures \ref{fig:densitywave_near0b} and \ref{fig:densitywave_near0c}). On the contrary, when the solution is prevented from entering the near-vacuum regime, no significant perturbation growth is observed, despite the presence of eigenvalues with large positive real parts in the linearized operator.

More work remains to be done, however, before we can confidently claim that the conclusions of \S\ref{sec_logarithmic_LCE_experiments} and \S\ref{sec:nearvacuum} are entirely applicable to density-wave problem, as the inaccuracy of the logarithmic mean is no longer the only issue that causes robustness issues. For example, the kinetic-energy-preserving numerical two-point fluxes of Kennedy and Gruber~\cite{kennedy_gruber_2008} and Kuya \etal~\cite{kuya_totani_kawai_2018} are also known to lack robustness for this problem, despite relying solely on arithmetic means~\cite{shima_flux,singh2021,coppola2024}. Likewise, repeating the numerical experiment shown in Figure \ref{fig:Euler_densitywave_near0} for the entropy-conservative two-point fluxes of Chandrashekar~\cite{chandrashekar2013} and Ismail and Roe~\cite{ismail_roe_2009} reveals that the crude positivity-preserving filter \eqref{eq_positivity_filter} is insufficient to prevent both positivity issues and perturbation growth driven by transient amplification, indicating that at least one other issue is likely at play. Future work should address the importance of satisfying additional properties at both the nonlinear level, such as the preservation of constant velocity and pressure~\cite{shima_flux,ranocha_2022_preventing_pressure}, and for the linearized problem, such as the conservation of a reduced energy for the acoustic subsystem~\cite{singh2021}. 

In the broader context of the Euler and Navier--Stokes equations, this challenge is further complicated by the lack of smooth benchmark problems for which the perturbation dynamics of the continuous equations are known in sufficient detail. Future work could investigate the one-dimensional time-periodic sound waves of Temple and Young~\cite{temple_young_sound}, which may provide a useful setting in which to apply the Floquet analysis introduced in this work. Another promising class of examples is compressible Couette flow, for which rigorous linearized transient-growth results are available in both the isentropic and non-isentropic settings~\cite{antonelli_dolce_marcati_couette,nonisentropic_couette}. Steady unstable flows may also provide useful test cases. For example, unstable eigenvalues have recently been identified for the incompressible Taylor--Green vortex~\cite{protas}, suggesting the possibility of comparing the spectra of flux-differencing discretizations (with or without numerical dissipation) to the physical growth rates of the continuous problem. Lastly, evaluating the effect of physical dissipation mechanisms on perturbation growth in the context of the Navier--Stokes equations is also an important consideration.

\section{Conclusions}

This work has examined the practical significance of local linear instabilities in entropy-stable SBP discretizations. Given the existence of nonlinear stability bounds, linear instabilities should be interpreted primarily as a mechanism of introducing numerical error, whose practical impact must be assessed relative to more dominant error sources. Analyzing the unphysical growth of such perturbations naturally leads to the study of a variable-coefficient linear advection problem. With this in mind we have identified two distinct scenarios. 

For linearizations of split-form discretizations of the Burgers equation, local linear instabilities are a genuine feature of the semidiscretization. For non-central splittings, the semidiscrete Jacobian can have eigenvalues with positive real parts, which, for frozen or slowly varying baseflows, can lead to unphysical exponential perturbation growth. However, we have shown that this growth is bounded by an estimate resembling the transient growth already present in the continuous variable-coefficient problem. Moreover, the modes that most strongly excite this growth are typically highly oscillatory, localized, and unphysical, often arising from modified boundary stencils. As a result, the practical impact of these instabilities on numerical accuracy is usually small compared with other sources of numerical error, and remains negligible on the time scales before the baseflow develops shocks. Even when this growth is not negligible, it is readily controlled by targeted numerical dissipation: interface dissipation alone for element-type or Mattsson operators, volume dissipation for circulant or continuous SBP operators, and both interface and volume dissipation for classical SBP operators.

Flux-differencing schemes with nonlinear two-point fluxes, such as those used in entropy-stable discretizations of the Euler equations, are more subtle, as their linearizations lead to time-varying operators even for linear time-invariant model problems. In this setting, frozen-baseflow eigenvalue analysis is a poor predictor of perturbation behaviour. Using the geometric mean for linear advection as a model problem, we obtained a sharp bound on perturbation growth. Numerical experiments with both the geometric and logarithmic means then confirmed that perturbations exhibit no meaningful growth, despite the presence of eigenvalues with large positive real parts. The time-varying linear operators evolve with the baseflow such that the growth and decay directions continuously rotate. Consequently, the baseflow dynamics cannot be separated from the perturbation dynamics: perturbations evolve under the accumulated effect of the rotating operators rather than any single frozen operator. As a result, transient dynamics dominate the modal growth predicted by instantaneous spectra. To accurately capture this behaviour, we introduced Floquet analysis as a diagnostic for the perturbation dynamics of time-periodic numerical discretizations. The Floquet exponents showed that the perturbation dynamics of the geometric and logarithmic schemes remain purely convective, even though frozen-time eigenvalue analyses suggest modal instability.

We have also explored the failure of entropy-stable schemes to robustly simulate density waves in near-vacuum regions. Whereas these failures were previously attributed to unphysical perturbation growth caused by local linear instabilities~\cite{gassner_stability_2022,ranocha_2022_preventing_pressure}, our results suggest that they are better understood as an accuracy and conditioning problem associated with logarithmic means near vacuum states. Entropy-conservative two-point fluxes become poorly behaved as the density approaches zero, leading to large numerical errors and highly nonnormal semidiscrete Jacobians with large spectral radii. This can indeed produce catastrophic short-time perturbation growth, but primarily through transient amplification rather than modal growth of frozen operators. Numerical experiments suggest that by adressing this issue directly, the accuracy and robustness of entropy-stable schemes in near-vacuum regions can be improved significantly.

Overall, these results indicate that local linear instabilities in entropy-stable schemes should not be viewed as a major obstruction to their practical use. While it is true that central schemes have valuable spectral stability properties, in comparison, entropy-stable schemes enforce additional structure at the nonlinear level that is often essential for robustly and accurately simulating nonlinear problems. An inability to respect the same frozen-time spectral bound as the continuous operator is therefore not, by itself, a sufficient reason to reject such schemes. In the examples considered here, local linear instabilities either have a negligible impact on numerical accuracy, or are artifacts of applying a time-invariant diagnostic to an intrinsically time-varying perturbation problem.

Several directions remain open for future work. First, the near-vacuum behaviour of entropy-conservative fluxes deserves further study. As suggested in \S\ref{sec:Euler}, the robustness of entropy-stable discretizations for the density-wave problem may depend on whether the scheme preserves additional structural properties of the continuous problem that were not considered here. Given the poor performance of the logarithmic mean in near-vacuum regions, it may also be interesting to investigate hybrid flux constructions. For example, alternating between entropy-conservative fluxes and robust arithmetic-mean-based fluxes (\eg Shima \etal~\cite{shima_flux}), or the asymptotically entropy-conservative fluxes of~\cite{de_michele_coppola_2023}, may improve both the robustness and accuracy of entropy-stable schemes in near-vacuum regions. Finally, analyzing numerical perturbation dynamics in more general Euler and Navier--Stokes settings remains an important extension of this work. We refer to the end of \S\ref{sec:Euler} for a discussion on possible test problems suitable for future experiments.

\begin{acknowledgements}
    The first author is grateful to Bartosz Protas for a helpful discussion in which he brought Floquet analysis to his attention.
\end{acknowledgements}

\section*{Statements and declarations}

\subsection*{Funding}
This research was supported by the Government of Ontario and the University of Toronto.


\subsection*{Conflict of interest}
The authors declare that they have no conflict of interest.

\subsection*{Data Availability}
The scripts to reproduce all of the numerical results are accessible via the open-source repository \url{https://github.com/alexbercik/ESSBP/tree/main/Paper-LinearStability}. 

\appendix


\section{Linearizations of the Flux-Differencing Schemes} \label{app_linearizations}

\subsection{General Flux-Differencing Schemes}

We consider the semidiscrete flux-differencing residual
\begin{equation}\label{eq:R_def}
  \frac{\mathrm{d}u_i}{\mathrm{d}t} = \mathcal{R}_i(\bu),
  \qquad
  \mathcal{R}_i(\bu) \defn - 2\sum_{j=1}^N \mat{D}_{ij}\,\fnc{F}^\star(u_i,u_j),
\end{equation}
where $\fnc{F}^\star$ is a smooth, symmetric, and consistent two-point flux.
We linearize $\bu$ about a positive baseflow $\tilde{\bu} > 0$ through $\bu = \tilde{\bu} + \bv$. The linearized perturbation equation then takes the form
\begin{equation}\label{eq:lin_general}
  \frac{\mathrm{d}\bv}{\mathrm{d}t} = \mat{J}(\tilde{\bu}(t))\,\bv ,
  \qquad
  \mat{J}_{ij}(\tilde{\bu}) \defn \left.\frac{\partial \mathcal{R}_i}{\partial u_j}\right|_{\bu=\tilde{\bu}}.
\end{equation}
The entries of $\mat{J}$ can be calculated explicitly by first defining $\partial_1\fnc{F}^\star$ and $\partial_2\fnc{F}^\star$ as the partial derivatives
with respect to the first and second arguments. Differentiating \eqref{eq:R_def} and dropping the tilde for simplicity yields
\begin{align}
  \mat{J}_{ik}(\bu)
  &= - 2\,\mat{D}_{ik}\,\partial_2\fnc{F}^\star(u_i,u_k) - 2\,\delta_{ik}\sum_{j=1}^N \mat{D}_{ij}\,\partial_1\fnc{F}^\star(u_i,u_j),
  \label{eq:L_ik_general}
\end{align}
or equivalently in matrix form,
\begin{equation}\label{eq:L_matrix_general}
  \mat{J}(\bu)
  = - 2\,\mat{D}\circ \partial_2 \mat{F} - 2\, \diag{  (\mat{D}\circ \partial_1 \mat{F})\,\vec{1} } , \quad \mat{F}_{ij} \defn \fnc{F}^\star(u_i,u_j)
\end{equation}
where $\circ$ denotes entry-wise multiplication, or the Hadamard product.

\subsection{The Logarithmic and Geometric Fluxes}

For the logarithmic flux \eqref{eq_log_flux}, we calculate the partial derivatives as
\begin{gather*}
  \left[ \partial_1 \mat{F} \right]_{ij} = \pder[\fnc{F}^\star(u_i,u_j)]{u_i} = \pder[a_i u_i]{u_i} \pder[\mean{a u}^{ij}_\text{log}]{a_i u_i} 
= a_i \frac{\mean{a u}^{ij}_\text{log} / u_i - 1}{\log(a_j u_j) - \log(a_i u_i)}, \\
\left[ \partial_2 \mat{F} \right]_{ij}
= a_j \frac{1 - \mean{a u}^{ij}_\text{log} / u_j}{\log(a_j u_j) - \log(a_i u_i)} .
\end{gather*}
For the diagonal entries $i=j$, we can directly appeal to consistency and smoothness of $\fnc{F}^\star$ to say
\[
  \left[ \partial_1 \mat{F} \right]_{ii} = \left[ \partial_2 \mat{F} \right]_{ii} = \pder[\fnc{F}^\star(u_i,u_i)]{u_i}
  = \frac{a_i}{2}.
\]
For the geometric flux \eqref{eq_geom_flux}, direct differentiation again gives
\[
 \left[ \partial_1 \mat{F} \right]_{ij} = \pder[\fnc{F}^\star(u_i,u_j)]{u_i} = \pder[a_i u_i]{u_i} \pder[\mean{a u}^{ij}_\text{geom}]{a_i u_i} 
= \frac{a_i}{2} \sqrt{ \frac{a_j u_j}{a_i u_i} },
\quad
\left[ \partial_2 \mat{F} \right]_{ij}
= \frac{a_j}{2} \sqrt{ \frac{a_i u_i}{a_j u_j} } .
\]
Substituting this into \eqref{eq:L_matrix_general}
yields a compact matrix representation, 
\[
\mat{J}_{\text{geom}}(\bu)
  = - \diag{\mat{D} \vec{w} } \mat{A} \mat{W}^{-1}  -  \mat{W} \mat{D} \mat{A} \mat{W}^{-1} , \quad w_i \defn \sqrt{a_i u_i}, \ \mat{W} \defn \diag{\vec{w}}.
\]

\subsection{Exponential Linearized Stability of the Geometric Flux} \label{app_linearizations_exp_bound}

The first observation we can make is that $\mat{J}_\text{geom}$ is a consistent discretization of the linearized equation. If we let $\vec{c} = \mat{A} \mat{W}^{-1} \bv$, $\fnc{C} = a(x) \fnc{V} / \fnc{W}$, and $\fnc{W} = \sqrt{a(x) \fnc{U}}$, then
\[
\mat{J}_\text{geom} \bv =  - \diag{\mat{D} \vec{w} } \vec{c} - \mat{W} \mat{D} \vec{c} \approx - \pder[\fnc{W}]{x} \fnc{C} - \fnc{W} \pder[\fnc{C}]{x} = - \pder[(\fnc{W} \fnc{C})]{x} = - \pder[(a(x) \fnc{V})]{x} .
\]
As we did in \S\ref{sec_stability_proof}, we can therefore introduce a geometric product-rule defect term
\[
\mat{\Theta}_\text{geom} \defn \left( \mat{D} \mat{W}  - \mat{W} \mat{D} - \diag{\mat{D}\vec{w}} \right) \mat{A} \mat{W}^{-1} \quad \Rightarrow \quad
\mat{J}_{\text{geom}} = - \mat{D}\mat{A} + \mat{\Theta}_{\text{geom}} = \mat{J}_{\text{cent}} + \mat{\Theta}_{\text{geom}}.
\]
Now that we have reduced the problem to a central scheme plus a defect term, we can re-use the stability proof from \S\ref{sec_stability_proof}. Following the same steps (setting $\alpha=1$), we find
\[
    \abs{\bv^\T \Hnrm \mat{J}_\text{geom} \bv} = \abs{\tfrac{1}{2}  \bv^\T \left( \mat{Q} \mat{A} - \mat{A} \mat{Q} \right) \bv + \bv^\T \Hnrm \mat{\Theta}_\text{geom} \bv }
    \leq \tfrac{1}{2} \norm{ \diag{\mat{D} \vec{a}} + \mat{\Theta} + 2 \mat{\Theta}_\text{geom} }_\Hnrm \norm{\bv}^2_\Hnrm ,
\]
where $\mat{\Theta}$ is the (arithmetic) product-rule defect introduced in \S\ref{sec_stability_proof}. We then establish boundedness through
\begin{align*}
    \der{t} \norm{\bv}^2_\Hnrm  &\leq \gamma \norm{\bv}^2_\Hnrm , \quad \gamma \defn \norm{\diag{\mat{D} \vec{a}} + \mat{\Theta} + 2 \mat{\Theta}_\text{geom}}_\Hnrm  .
\end{align*}
which again is structurally similar to the continuous case, where $\gamma = - \min_{x \in \Omega} \partial a / \partial x$ \cite{kopriva_energy_2014}. Convergence of the semidiscretization \eqref{eq:R_def} for the continuous problem \eqref{eq_lin_conv} then follows through the arguments of \cite{strang}.

\subsection{Sharp Linearized Stability of the Geometric Flux} \label{app_sharp_geom_bound}

Although in the previous section we established an exponential energy bound for perturbations of the geometric flux-differencing scheme, such a bound is not sharp. Alternatively, by considering weighted energy norms, similar to~\cite{manzanero}, we can construct a sharper bound. We begin with the linearized geometric scheme,
\[
\der[\bv]{t} + \diag{\mat{D} \vec{w} } \mat{A} \mat{W}^{-1} \bv  +  \mat{W} \mat{D} \mat{A} \mat{W}^{-1} \bv = \vec{0} , \quad w_i \defn \sqrt{a_i u_i}, \ \mat{W} \defn \diag{\vec{w}}.
\]
We introduce the intermediate variable $\vec{z} \defn \mat{W}^{-1} \bv$, such that the linearized scheme can be rewritten as
\[
\der[\bv]{t} + \diag{\mat{D} \vec{w} } \mat{A} \vec{z}  +  \mat{W} \mat{D} \mat{A} \vec{z} = \vec{0} .
\]
Using the relations 
\[
\der[z_i]{t} = \der[w_i^{-1}]{t}  v_i + w_i^{-1}  \der[v_i]{t} 
\quad \text{and} \quad
\der[w_i^{-1}]{t}  v_i =  - \der[w_i]{t} \frac{v_i}{w_i^2} = - \der[w_i]{t} \frac{z_i}{w_i},
\]
 we obtain $\displaystyle \der[v_i]{t} = w_i \der[z_i]{t} + z_i \der[w_i]{t}$, allowing us rewrite the geometric scheme in terms of $\vec{z}$ as
\begin{gather*}
    \mat{W} \der[\vec{z}]{t} + \diag{\der[\vec{w}]{t}} \vec{z} + \diag{\mat{D} \vec{w} } \mat{A} \vec{z}  +  \mat{W} \mat{D} \mat{A} \vec{z} = \vec{0} \\
    \Rightarrow \quad \der[\vec{z}]{t} + \mat{D} \mat{A} \vec{z} = - \mat{W}^{-1} \diag{\der[\vec{w}]{t} + \mat{A} \mat{D} \vec{w} } \vec{z} .
\end{gather*}
Therefore, the intermediate variable $\vec{z}$ evolves according to a central semidiscretization of the conservative linear advection equation, consistent with the continuous problem, but now with an additional source term. Notably, the term in brackets is the residual of a central discretization of the non-conservative linear advection equation in the variable $\fnc{W}$. However, this residual is directly related to the residual of the baseflow discretization. Since the baseflow geometric flux-differencing scheme satisfies
\begin{equation} \label{eq_geom_scheme_app_lin}
\der[\bu]{t} + 2 \diag{\sqrt{\mat{A} \bu}} \mat{D} \sqrt{\mat{A} \bu} = \vec{0} \quad \Leftrightarrow \quad \der[\bu]{t} + 2 \mat{W} \mat{D} \vec{w} = \vec{0} ,
\end{equation}
then as shown in Appendix \ref{app_reformulation_geom}, a reformulation of the baseflow geometric scheme in terms of $\vec{w}$ reveals that
\[
\vec{r} \defn \der[\vec{w}]{t} + \mat{A} \mat{D} \vec{w} = \vec{0} .
\]
Therefore, the source term vanishes, indicating that perturbations of the nonlinear geometric flux equation, when viewed in the intermediate variable $\fnc{Z} \defn \fnc{V} / \sqrt{a \fnc{U}}$, evolve exactly according to a central discretization. Since by assumption $a_i , u_i > 0$, we can immediately apply the weighted energy estimates of~\cite{manzanero} to say
\[
\der{t} \norm{\vec{z}}_{a \Hnrm}^2 = 0 , \quad \text{where}  \ \norm{\vec{z}}_{a \Hnrm}^2 \defn \sum_{i=1}^N \Hnrm_{ii} a_i z_i^2 .
\]
Converting from the intermediate variable back to the perturbation $\fnc{V}$, we obtain
\[
\der{t} \norm{\vec{v}}_{\Hnrm/\mat{U}}^2 = 0 , \quad \text{where}  \ \norm{\vec{v}}_{\Hnrm/\mat{U}}^2 \defn \sum_{i=1}^N \Hnrm_{ii} \frac{v_i^2}{u_i} .
\]
Therefore, given any baseflow state $u(x,t)$, which is bounded by the nonlinear entropy stability bound discussed in \S\ref{sec_stability_proof}, through equivalence of norms we find
\[
\frac{\norm{\vec{v}(t)}_\Hnrm^2}{\max_{i}{(u_i(t))}} \leq \norm{\vec{v}(t)}_{\Hnrm/\mat{U}}^2 \leq \frac{\norm{\vec{v}(t)}_\Hnrm^2}{\min_{i}{(u_i(t))}} \quad \Rightarrow \quad \norm{\vec{v}(t)}_{\Hnrm}^2 \leq \frac{\max_{i}{(u_i(t))}}{\min_{i}{(u_i(0))}} \norm{\vec{v}(0)}_{\Hnrm}^2 .
\]
If we consider a different baseflow $u_i(t) = \fnc{U}(x_i,t)$ than that of the numerical solution of the geometric flux-differencing scheme \eqref{eq_geom_scheme_app_lin}, then it will not in general be true that $\vec{r} = \vec{0}$. Consequently, we cannot set the right-side source term to zero. However, if $\fnc{U}(x_i,t)$ is a smooth and bounded function, such as the exact solution to the continuous PDE, then by consistency $\vec{r} = \order{p}$. Therefore, while we can expect to see some small exponentional perturbation growth on any given mesh, it will vanish rapidly as $\order{p}$ under mesh refinement. This occurs, for example, when performing Floquet analysis.

\subsection{Symmetric Part of the Jacobian} \label{app_linearizations_sym}

In \S\ref{sec:nearvacuum} we analyze the symmetric part of the semidiscrete linearized Jacobians, $\mat{J}_\mathrm{sym} \defn \tfrac{1}{2} (\mat{J} + \mat{J}^\T)$. For simplicity, hereafter we restrict ourselves to $a(x)=1$. For the logarithmic flux, we have already shown that
\begin{align*}
\text{for} \  i \neq j , \quad \mat{J}_{ij} &= - 2 \mat{D}_{ij} \frac{1 - \mean{u}^{ij}_\text{log} / u_j}{\log(u_j) - \log(u_i)} \\
&= - 2 \mat{D}_{ij} \left( \frac{1}{\log(u_j) - \log(u_i)} - \frac{1 - u_i / u_j}{(\log(u_j) - \log(u_i))^2} \right) .
\end{align*}
Letting $\mat{D} = - \mat{D}^\T$, or $\mat{D}_{ij} = - \mat{D}_{ji}$, we find
\begin{align*}
\text{for} \  i \neq j , \quad (\mat{J}_\mathrm{sym})_{ij} &= 
- \mat{D}_{ij} \left( \frac{2}{\log(u_j) - \log(u_i)} + \frac{u_i / u_j - u_j / u_i }{(\log(u_j) - \log^2(u_i)} \right) \\
&= - \mat{D}_{ij}  \frac{2 \log(u_j / u_i) + u_i / u_j - u_j / u_i }{(\log^2(u_j / u_i)} .
\end{align*}
For the diagonal components, $\mat{D}_{ii} = 0$, so
\[
\quad (\mat{J}_\mathrm{sym})_{ii} = \mat{J}_{ii} = - 2 \sum_{j \neq i}^N \mat{D}_{ij} \frac{\mean{u}_\text{log}^{ij} / u_i - 1}{\log(u_j) - \log(u_i)} 
= - 2 \sum_{j \neq i}^N \mat{D}_{ij} \frac{u_j / u_i - \log(u_j / u_i) - 1}{(\log^2(u_j / u_i)} .
\]
For the geometric flux, we can directly use the matrix representation,
\[
\mat{J}_\mathrm{sym} = - \diag{\mat{D} \vec{w} } \mat{W}^{-1}  -  \tfrac{1}{2} \left( \mat{W} \mat{D} \mat{W}^{-1} - \mat{W}^{-1} \mat{D} \mat{W} \right) , \quad \vec{w} \defn \sqrt{\bu}, \ \mat{W} \defn \diag{\vec{w}} ,
\]
or in index notation,
\[
\text{for} \  i \neq j , \quad (\mat{J}_\mathrm{sym})_{ij} = - \tfrac{1}{2} \mat{D}_{ij} \left( \sqrt{\frac{u_i}{u_j}} - \sqrt{\frac{u_j}{u_i}} \right) = \tfrac{1}{2} \mat{D}_{ij} \frac{u_j - u_i}{\sqrt{ u_i u_j}} , 
\]
and for the diagonal components,
\[
(\mat{J}_\mathrm{sym})_{ii} = - \sum_{j \neq i}^N \mat{D}_{ij} \sqrt{\frac{u_j}{u_i}} .
\]

\section{Reformulation of the Geometric Flux-Differencing Scheme} \label{app_reformulation_geom}

Here we provide details on the reformulation of the geometric flux-differencing scheme as a central discretization. Consider an SBP-SAT semidiscretization using the geometric flux \eqref{eq_geom_flux} with periodic boundary conditions and standard nondissipative SATs (\eg as one recovers from the entropy-stable framework of \cite{crean}),
\begin{equation} \label{eq_geom_scheme_app}
    \der[\bu]{t} + 2 \mat{W} \mat{D} \vec{w} = \Hnrm^{-1} \mat{W} \left[ \mat{E} \vec{w} - \tr \tlT \vec{w}_R + \tl \trT \vec{w}_L \right] , \quad \mat{W} \defn \diag{\vec{w}}, \ w_i \defn \sqrt{a_i u_i} ,
\end{equation}
where $\vec{w}_{R/L}$ correspond to the solutions $\bu_{R/L}$ of the right and left adjacent blocks, respectively. Continuous SBP discretizations, or those with circulant operators, can be included by setting the right-hand side to zero. Consider that
\[
\der[w_i]{t} = \der[\sqrt{a_i u_i}]{t} = \frac{1}{2 \sqrt{a_i u_i}}\der[(a_i u_i)]{t} = \frac{a_i}{2 w_i} \der[u_i]{t}, \quad \text{so} \quad \der[\bu]{t} = 2 \mat{W} \mat{A}^{-1} \der[\vec{w}]{t} ,
\]
where we used that $a(x)$ is independent of $t$. Plugging this into \eqref{eq_geom_scheme_app}, we obtain the equivalent formulation
\[
    \der[\vec{w}]{t} + \mat{A} \mat{D} \vec{w} = \tfrac{1}{2} \Hnrm^{-1} \mat{A} \left[ \mat{E} \vec{w} - \tr \tlT \vec{w}_R + \tl \trT \vec{w}_L \right] .
\]
Therefore, although the geometric flux-differencing scheme \eqref{eq_geom_scheme_app} is nonlinear in $\bu$, it is equivalent to a standard linear SBP-SAT semidiscretization of the non-conservative variable-coefficient linear advection equation in the square-root variable $\fnc{W} \defn \sqrt{a(x) \fnc{U}}$. Following the framework of \cite{manzanero}, we can premultiply by $\vec{w}^\T \mat{A}^{-1} \Hnrm$, use the SBP property, and sum over neighboring blocks to obtain the following weighted energy estimate,
\[
\der{t} \norm{\vec{w}}_{\Hnrm/a}^2 = 0 , \quad \text{where}   \ \norm{\vec{w}}_{\Hnrm/a}^2 \defn \sum_{i=1}^N \Hnrm_{ii} \frac{w_i^2}{a_i} .
\]
Converting from the square-root variable back to the conservative variable $\fnc{U}$, we obtain
\[
\der{t} \norm{\vec{u}}_{1,\Hnrm}^2 = 0 , \quad \text{where}  \ \norm{\vec{u}}_{1,\Hnrm}^2 \defn \sum_{i=1}^N \Hnrm_{ii} \abs{u_i} = \sum_{i=1}^N \Hnrm_{ii} u_i.
\]
This establishes discrete conservation of the conservative variable, but unfortunately does not on its own establish a stability bound in the conservative variable. For this we must rely on alternative arguments, as is discussed in Appendix~\ref{app_entropy_l2_bound}. Regardless, this reformulation in terms of the square-root variable is valuable as it allows us to establish sharp stability bounds for perturbations of the nonlinear scheme in the conservative variable, which is discussed in Appendix \ref{app_sharp_geom_bound}.

\section{Details from Frozen-Coefficient Analysis} \label{app_frozen_coeff}

We begin by analyzing the modified PDE introduced in \S\ref{sec_modified_pde}. For the geometric scheme, we found \eqref{eq_modified_pde_geom}, or
\begin{gather*}
\pder[\fnc{V}]{t} + \pder{x} \left[ \left( a(x) - \beta(x) \right)  \fnc{V} \right] = \pder{x} \left( \nu(x) \pder[\fnc{V}]{x} \right) ,\\
\beta(x) \defn \frac{h^2}{8} a(x) \left( \frac{1}{a(x)^2} \left(\pder[a(x)]{x}\right)^2 - \frac{1}{\fnc{U}^2} \left(\pder[\fnc{U}]{x}\right)^2 \right) , \quad \nu(x) \defn \frac{h^2}{4 \ \fnc{U}} \pder[a(x) \fnc{U}]{x}.
\end{gather*}
For the central-product split-form scheme, we found the simpler result 
\[
\beta(x) \defn 0 , \quad \nu(x) \defn \frac{1-\alpha}{2} h^2 \pder[a(x)]{x} .
\]
As explained in \S\ref{sec_frozen_coeff}, we now make a locality assumption such that
\[
\delta \defn x - x_0 , \quad a(x) \approx a_0 + a_1 \delta , \quad \fnc{U}(x) \approx \fnc{U}_0 + \fnc{U}_1 \delta .
\]
For the central-product split-form scheme, this immediately yields $\nu \approx h^2 \tfrac{1-\alpha}{2} a_1$. For the geometric scheme, however, this still yields nonlinear coefficients,
\[
\beta(x) \approx \frac{h^2}{8} ( a_0 + a_1 \delta ) \left( \frac{a_1^2}{( a_0 + a_1 \delta )^2} - \frac{\fnc{U}_1^2}{( \fnc{U}_0 + \fnc{U}_1 \delta )^2} \right) , \quad \nu(x) \approx \frac{h^2}{4} \left(  a_1 + \fnc{U}_1 \frac{ a_0 + a_1 \delta}{\fnc{U}_0 + \fnc{U}_1 \delta } \right).
\]
By assumption, $a_0$ and $\fnc{U}_0$ are bounded away from zero. Therefore, it is appropriate to take linear approximations through Taylor expansions of the above, obtaining
\begin{gather*}
\beta(x) \approx h^2 (\beta_0 + \beta_1 \delta) , \quad 
\beta_0 \defn \frac{a_0}{8} \left( \frac{a_1^2}{a_0^2} - \frac{\fnc{U}_1^2}{\fnc{U}_0^2} \right) , \quad 
\beta_1 \defn - \frac{a_1}{8} \left( \frac{a_1^2}{a_0^2} + \frac{\fnc{U}_1^2}{\fnc{U}_0^2} - 2 \frac{a_0}{a_1} \frac{\fnc{U}_1^3}{\fnc{U}_0^3} \right) , \\
\nu(x) \approx h^2 (\nu_0 + \nu_1 \delta) , \quad 
\nu_0 \defn \frac{1}{4} \left(  a_1 + a_0 \frac{\fnc{U}_1}{\fnc{U}_0 } \right) , \quad 
 \nu_1 \defn \frac{1}{4} \left( a_1 \frac{ \fnc{U}_1}{\fnc{U}_0 } - a_0 \frac{\fnc{U}_1^2}{\fnc{U}_0^2 } \right).
\end{gather*}
With these linear simplifications, the modified PDE becomes
\[
\pder[\fnc{V}]{t} + (a_0 - h^2 \beta_0 - h^2 \nu_1) \pder[\fnc{V}]{x} + (a_1 - h^2 \beta_1) \delta \pder[\fnc{V}]{x} + (a_1 - h^2 \beta_1) \fnc{V}  = h^2 (\nu_0 + \nu_1 \delta) \pder[^2\fnc{V}]{x^2}  .
\]
For the central-product split-form scheme, we have $\nu_0 \defn \frac{1-\alpha}{2} a_1$ and $\beta_0 = \beta_1 = \nu_1 = 0$. As explained in \S\ref{sec_frozen_coeff}, we then make the assumption that $\abs{ a_1 \delta} \ll \abs{a_0}$ and $\abs{\fnc{U}_1 \delta} \ll \abs{\fnc{U}_0}$. This allows us to drop certain terms from the modified PDE because their contributions are negligible. For example,
\[
\abs{a_1 \delta \pder[\fnc{V}]{x}} \leq  \abs{a_1  \delta} \abs{\pder[\fnc{V}]{x}} \ll \abs{a_0} \abs{\pder[\fnc{V}]{x}} = \abs{a_0 \pder[\fnc{V}]{x}} .
\]
For the central-product split-form scheme, this is the only simplification needed. For the geometric scheme, we recognize that
\[
\beta_0 + \nu_1 =  \frac{a_0}{8} \frac{a_1^2}{a_0^2} - \frac{a_0}{8} \frac{\fnc{U}_1^2}{\fnc{U}_0^2} - \frac{ a_0}{4} \frac{\fnc{U}_1^2}{\fnc{U}_0^2} + \frac{a_1}{4} \frac{ \fnc{U}_1}{\fnc{U}_0 } , \quad
\beta_1 \delta = - \frac{a_1}{8} \frac{a_1^2}{a_0^2} \delta - \frac{a_1}{8} \frac{\fnc{U}_1^2}{\fnc{U}_0^2} \delta + \frac{a_0}{4} \frac{\fnc{U}_1^3}{\fnc{U}_0^3} \delta .
\]
We wish to show that the terms of $\beta_1 \delta$ are negligible compared to those of $\beta_0 + \nu_1$. We see that 
\begin{gather*}
\abs{-\frac{a_1}{8} \frac{a_1^2}{a_0^2} \delta} \leq \abs{a_1 \delta} \abs{ \frac{1}{8}\frac{a_1^2}{a_0^2}} \ll \abs{a_0} \abs{ \frac{1}{8}\frac{a_1^2}{a_0^2}} = \abs{ \frac{a_0}{8}\frac{a_1^2}{a_0^2}} , \\
\abs{-\frac{a_1}{8} \frac{\fnc{U}_1^2}{\fnc{U}_0^2} \delta} \leq \abs{a_1 \delta} \abs{\frac{1}{8} \frac{\fnc{U}_1^2}{\fnc{U}_0^2}} \ll \abs{a_0} \abs{\frac{1}{8} \frac{\fnc{U}_1^2}{\fnc{U}_0^2}} = \abs{\frac{a_0}{8} \frac{\fnc{U}_1^2}{\fnc{U}_0^2}} , \\
\abs{\frac{a_0}{4} \frac{\fnc{U}_1^3}{\fnc{U}_0^3} \delta} \leq \frac{\abs{\fnc{U}_1 \delta}}{\abs{\fnc{U}_0}} \abs{ \frac{a_0}{4} \frac{\fnc{U}_1^2}{\fnc{U}_0^2}} \ll \abs{ \frac{a_0}{4} \frac{\fnc{U}_1^2}{\fnc{U}_0^2}} .
\end{gather*}
Therefore, all three terms of $\beta_1 \delta$ are negligible compared to those of $\beta_0 + \nu_1$, so can be dropped without affecting the dynamics of the modified PDE. Note that we could use the same logic to drop the fourth term of $\beta_0 + \nu_1$. We also examine the terms of $\nu_1 \delta$ and $\nu_0$,
\[
\nu_0 = \frac{a_1}{4} + \frac{a_0}{4} \frac{\fnc{U}_1}{\fnc{U}_0 }  , \quad 
 \nu_1 = \frac{a_1}{4} \frac{ \fnc{U}_1}{\fnc{U}_0 } - \frac{a_0}{4} \frac{\fnc{U}_1^2}{\fnc{U}_0^2 } .
\]
The corresponding first and second terms can be immediately compared through
\[
\abs{\frac{a_1}{4} \frac{ \fnc{U}_1}{\fnc{U}_0 } \delta} \leq \abs{\frac{a_1}{4}} \frac{\abs{\fnc{U}_1 \delta}}{\abs{\fnc{U}_0}} \ll \abs{\frac{a_1}{4}} , \quad 
\abs{\frac{a_0}{4} \frac{\fnc{U}_1^2}{\fnc{U}_0^2 }} \leq \abs{\frac{a_0}{4} \frac{\fnc{U}_1}{\fnc{U}_0 }} \frac{\abs{\fnc{U}_1 \delta}}{\abs{\fnc{U}_0}} \ll \abs{\frac{a_0}{4} \frac{\fnc{U}_1}{\fnc{U}_0 }} ,
\]
allowing us to drop the contributions from $\nu_1 \delta$. We are left with the local modified PDE stated in \eqref{eq_modified_pde_frozen},
\[
\pder[\fnc{V}]{t} + (a_0 - h^2 \beta_0 - h^2 \nu_1) \pder[\fnc{V}]{x} + (a_1 - h^2 \beta_1) \fnc{V}  = h^2 \nu_0 \pder[^2\fnc{V}]{x^2}  .
\]

\section{Catastrophic Errors of Flux-Differencing Schemes Near Vacuum} \label{app_nearvacuum}

In this section we explore the catastrophic errors introduced by the logarithmic and geometric flux-differencing schemes in near-vacuum states. We first show that in these regimes, the truncation error becomes excessively large. We then show that the linearized semidiscrete Jacobian becomes highly stiff and nonnormal, which has implications for time marching and unphysical perturbation growth. 

Consider a smooth nonnegative profile $\fnc{U}(x)$ satisfying $\min_x \fnc{U}(x) = 0$ and define the one-parameter family of solutions 
\[
    \fnc{U}_\epsilon(x) \defn \fnc{U}(x) + \epsilon, \quad \epsilon > 0.
\]
As $\epsilon \rightarrow 0$, this family approaches a near-vacuum state by shifting the same profile downward toward zero. This family is chosen so that the spatial variation of the profile is unchanged. In contrast, a family of the form $\fnc{U}_\epsilon = \epsilon \fnc{U}$ can avoid the numerical difficulties described below by trivially rescaling the variables.

Now consider a flux-differencing volume term \eqref{eq_flux_difference} of degree $p$ for the linear advection equation \eqref{eq_lin_conv}. Let $a(x)=1$, and consider finite-difference refinement where the mesh spacing $h$ is decreased by adding nodes. The exact solution is restricted to the mesh via $\bu = \fnc{U}_\epsilon(\vec{x})$, and for convenience we define $\phi_{i,\epsilon}(x) \defn 2 \fnc{F}^\star(u_i , \fnc{U}_\epsilon(x))$. At node $i$, the volume term \eqref{eq_flux_difference} acting on the exact solution can then be expressed as
\[
2 \sum_{j=1}^N \mat{D}_{ij} \fnc{F}^\star(u_i , u_j) 
= \sum_{j=1}^N \mat{D}_{ij} \phi_{i,\epsilon}(x_j) 
= \left. \pder[\phi_{i,\epsilon}]{x} \right\vert_{x=x_i} + c_i h^p \left. \pder[^{p+1}\phi_{i,\epsilon}]{x^{p+1}} \right\vert_{x=\xi_i}
\]
for some constant $c_i$ and point $\xi_i$ near $x_i$. For the central flux \eqref{eq_central_flux}, consistency is easily established via
\[
\pder[\phi_{i,\epsilon}(x)]{x} = \pder{x} \left(u_i + \fnc{U}_\epsilon(x) \right) = \pder[\fnc{U}_\epsilon(x)]{x} = \pder[\fnc{U}(x)]{x} .
\]
Repeating this argument for higher-order derivatives of $\phi_{i,\epsilon}$, we can easily show that the truncation error of the volume term \eqref{eq_flux_difference} is $c_i h^p \fnc{U}^{(p+1)}(\xi_i)$, where $\fnc{U}^{(k)} \defn \partial^{k} \fnc{U} / \partial x^{k}$. For the logarithmic flux \eqref{eq_log_flux}, more care is required. Letting $\fnc{U}_\epsilon' \defn \partial \fnc{U}_\epsilon / \partial x = \partial \fnc{U} / \partial x$, consistency is established through
\begin{align*}
    \left.\pder[\phi_{i,\epsilon}(x)]{x}\right\vert_{x=x_i} 
    &= \pder{x} \left. \left( \frac{ 2(\fnc{U}_\epsilon(x) - u_i)}{\log(\fnc{U}_\epsilon(x)) - \log(u_i)} \right) \right\vert_{x=x_i} \\
    &= 2 \, \fnc{U}_\epsilon'(x_i) \ \lim_{x \rightarrow x_i} \, \frac{\fnc{U}_\epsilon(x) \left(\beta_\epsilon(x)  - 1 \right) + u_i }{\fnc{U}_\epsilon(x) \beta_\epsilon^2(x) }  , \quad \beta_\epsilon(x) \defn \log(\fnc{U}_\epsilon(x)) - \log(u_i) \\
    &= 2 \, \fnc{U}_\epsilon'(x_i) \ \lim_{x \rightarrow x_i} \, \frac{\fnc{U}_\epsilon'(x) \left(\beta_\epsilon(x)  - 1 \right) + \fnc{U}_\epsilon'(x) }{\fnc{U}_\epsilon'(x) \beta_\epsilon^2(x) + 2 \, \fnc{U}_\epsilon'(x) \beta_\epsilon(x)} , \quad \text{since} \ \beta_\epsilon'(x) = \frac{\fnc{U}_\epsilon'(x)}{\fnc{U}_\epsilon(x)} \\
    &= 2 \, \fnc{U}_\epsilon'(x_i) \ \lim_{\beta_\epsilon(x) \rightarrow 0} \, \frac{1}{\beta_\epsilon(x) + 2 } = \fnc{U}_\epsilon'(x_i) = \fnc{U}'(x_i) .
\end{align*}
Computing higher order derivatives requires tedious algebra, but it is still straightforward to show that
\[
\phi_{i,\epsilon}''(x_i) = \fnc{U}''(x_i) - \frac{(\fnc{U}'(x_i))^2}{3 \, \fnc{U}_\epsilon(x_i)} , \quad \text{and} \quad \phi_{i,\epsilon}'''(x_i) = \fnc{U}'''(x_i) - \frac{\fnc{U}'(x_i) \fnc{U}''(x_i)}{\fnc{U}_\epsilon(x_i)} + \frac{(\fnc{U}'(x_i))^3}{2 \, \fnc{U}_\epsilon^2(x_i)} .
\]
Though a general expression for $\phi_{i,\epsilon}^{(k)}(x)$ can be found, it requires significant effort and patience. For the purposes of this discussion, we emphasize only the important point that $\phi_{i,\epsilon}^{(p+1)}(x)$---and hence the truncation error $c_i h^p \phi_{i,\epsilon}^{(p+1)}(\xi_i)$---contains terms with inverse powers of $\fnc{U}_\epsilon(x) = \fnc{U}(x)+\epsilon$, while the numerators contain derivatives of $\fnc{U}_\epsilon$, which are independent of $\epsilon$. Consequently, in regions near a zero of $\fnc{U}$, as $\epsilon \rightarrow 0$, the numerators remain unchanged while the denominator can be made arbitrarily small, and the truncation error grows rapidly.

Although the singular terms may vanish exactly at a smooth minimum $x_i = x_\text{min}$ if $\fnc{U}^{(k)}(x_\text{min})=0$ for some $k\geq 1$, the inverse-power dependence persists in the surrounding near-vacuum region. For a fixed $\epsilon>0$, the scheme retains its formal order of accuracy under mesh refinement, but the constant multiplying $h^p$ in the truncation error becomes increasingly singular in the regions surrounding $x_\text{min}$ as $\epsilon \rightarrow 0$. The same argument can be repeated for the geometric flux, yielding analogous expressions for $\phi_{i,\epsilon}^{(k)}(x)$, and therefore the same inverse-power dependence on $\fnc{U}_\epsilon(x)$ in near-vacuum regions.

A similar issue can be identified with respect to the nonnormality of the semidiscrete linearized Jacobians $\mat{J}(t)$. Instantaneous perturbation growth (including transient growth) is determined by the symmetric part of $\mat{J}(t)$, $\mat{J}_\mathrm{sym} \defn \tfrac{1}{2}(\mat{J} + \mat{J}^\T)$, which for the logarithmic flux with circulant operators is (see Appendix \ref{app_linearizations_sym})
\[
(\mat{J}_{\mathrm{sym}})_{ij} =
\begin{cases}
- \mat{D}_{ij}\dfrac{2 \log(u_j/u_i) + u_i/u_j - u_j/u_i}{\log^2(u_j/u_i)}, & i\neq j,\\[2.5ex]
- 2\displaystyle\sum_{k\neq i}\mat{D}_{ik}\dfrac{u_k/u_i - \log(u_k/u_i)-1}{\log^2(u_k/u_i)}, & i=j.
\end{cases}
\]
The maximum eigenvalue of $\mat{J}_\mathrm{sym}$ can be bounded from below by considering the restriction of $\mat{J}_\mathrm{sym}$ to the subspace $\text{span}\{\vec{e}_i , \vec{e}_j\}$, where $\vec{e}_k$ is the standard basis vector for node $k$, and $j$ is in the stencil of $i$. Let
\[
\mat P \defn [\mathbf e_i,\mathbf e_j]\in\mathbb R^{N\times 2},
\qquad
\tilde{\mat J}_{\mathrm{sym}} \defn \mat P^T \mat J_{\mathrm{sym}} \mat P
=
\begin{pmatrix}
a & b\\
b & c
\end{pmatrix},
\]
where $a \defn (\mat{J}_\mathrm{sym})_{ii}$, $b \defn (\mat{J}_\mathrm{sym})_{ij}$, and $c \defn (\mat{J}_\mathrm{sym})_{jj}$. If $\vec{y} \in \mathbb{R}^2$ is a non-zero eigenvector of $\tilde{\mat{J}}_\mathrm{sym}$ corresponding to its larger eigenvalue $\lambda_{\max}(\tilde{\mat{J}}_\mathrm{sym}) = \tfrac{1}{2} (a + c) + \sqrt{\tfrac{1}{4} (a - c)^2 + b^2}$, then
\[
\lambda_{\max}(\tilde{\mat{J}}_\mathrm{sym}) = \frac{\vec{y}^\T \tilde{\mat{J}}_\mathrm{sym} \vec{y}}{\vec{y}^\T \vec{y}}
= \frac{\left( \mat{P} \vec{y} \right)^\T \mat{J}_\mathrm{sym} \left( \mat{P} \vec{y} \right)}{\left( \mat{P} \vec{y} \right)^\T \left( \mat{P} \vec{y} \right)} 
\leq \lambda_{\max}(\mat{J}_\mathrm{sym}) \quad \text{since} \ \mat{P}^\T \mat{P} = \mat{I}_{2 \times 2}.
\]
Now fix a mesh through $h>0$, and suppose that node $i$ samples a zero of the base profile, $\fnc{U}(x_i)=0$, while $j$ is a stencil node with $\fnc{U}(x_j)>0$. Under the shifted family $\fnc{U}_\epsilon = \fnc{U}+\epsilon$, we then have
\[
u_i = \epsilon,
\quad
u_j = \fnc{U}(x_j) + \epsilon
\quad \Rightarrow \quad
r_j \defn \frac{u_j}{u_i} = \frac{\fnc{U}(x_j)+\epsilon}{\epsilon} \rightarrow \infty
\quad \text{as} \ \epsilon \rightarrow 0 .
\]
This is precisely the discrete counterpart of shifting the same profile downward toward vacuum: the neighboring differences are inherited from the fixed profile $\fnc{U}$, so the blowup can not be avoided by rescaling the variables. Then
\[
b = - \mat{D}_{ij}  \frac{2 \log(r_j) + 1 / r_j - r_j }{\log^2(r_j)} \sim  \mat{D}_{ij}  \frac{ r_j }{\log^2(r_j)} \rightarrow \sgn{(\mat{D}_{ij})} \cdot \infty \quad \text{as} \ r_j \rightarrow \infty .
\]
We now turn to the asymptotics of $a$. Let $r_k \defn  u_k / u_i$, where $k$ is also in the stencil of $i$. Under the same shifted-family assumption, we have $r_k \rightarrow \infty$ as $\epsilon \rightarrow 0$ whenever $\fnc{U}(x_k)>0$, while $r_k \rightarrow 1$ as $\epsilon \rightarrow 0$ whenever $\fnc{U}(x_k) = 0$ also. Thus each term
\[
a_k \defn - 2 \mat{D}_{ik}  \frac{r_k - \log(r_k) -1 }{\log^2(r_k)}
\]
is either bounded for a bounded $r_k$, or satisfies
\[
a_k \sim - 2 \mat{D}_{ik}  \frac{r_k }{\log^2(r_k)} \quad \text{as} \ r_k \rightarrow \infty.
\]
A similar analysis shows that since $\fnc{U}(x_j) > 0$ by assumption, $u_k / u_j$ is bounded for all $k$ in the stencil of $j$, and therefore each term $c_k$ in the sum of $c$ is bounded.

In the case that each $a_k$ is bounded (which as $\epsilon \rightarrow 0$ can only occur if $\fnc{U}(x_k) = 0$ also), both $a$ and $c$ are bounded, so $\lambda_{\max}(\tilde{\mat{J}}_\mathrm{sym}) \sim \abs{b} \rightarrow \infty$ as $\epsilon \rightarrow 0$. Otherwise, if at least one $a_k$ is unbounded (\ie $\fnc{U}(x_k) > \epsilon$), then $a$ might also grow unboundedly. However, since the stencil is fixed, at worst
\[
a = \alpha \abs{b} + \fnc{O}(1) \quad \text{as} \ \epsilon \rightarrow 0.
\]
for some constant $\alpha$. Then
\[
\lambda_{\max}(\tilde{\mat{J}}_\mathrm{sym}) \sim 
\left(
\frac{\alpha}{2} + \sqrt{\frac{\alpha^2}{4}+1}
\right) |b| + O(1) \rightarrow \infty \quad \text{as} \ \epsilon \rightarrow 0,
\]
since the coefficient in parentheses is strictly positive. Hence
\[
\lambda_{\max}(\mat{J}_\mathrm{sym}) \geq \lambda_{\max}(\tilde{\mat{J}}_\mathrm{sym}) \rightarrow \infty
\quad \text{as} \ \epsilon \rightarrow 0 .
\]
A similar analysis shows the same conclusion for the geometric flux. Note that this does not contradict the perturbation bound \eqref{eq_geom_pert_bound}, since (a) the weighted norm upon which it is based becomes singular, and (b) the bound is proportional to the baseflow through $\max_i{u_i(t)}$, which itself can grow if $u_i \rightarrow 0$ (see Appendix \ref{app_entropy_l2_bound}). It also does not contradict the exponential stability bound of Appendix~\ref{app_linearizations_exp_bound} because in this regime, $\norm{\mat{\Theta}_\text{geom}}_\Hnrm \rightarrow \infty$.

Finally, we also demonstrate that the spectral radius of the semidiscretization can grow unboundedly, \ie $\rho(\mat{J}) \rightarrow \infty$ as $u_i \rightarrow 0$ and $u_j \rightarrow \fnc{U}(x_j) > 0$, or equivalently, $r_j \rightarrow \infty$ as $\epsilon \rightarrow 0$. This has important implications for time marching. Let $\lambda_i$ be the eigenvalues of $\mat{J}$. We consider the trace of the matrix $\mat{J}^2$ since
\[
\abs{\text{tr}(\mat{J}^2)} = \abs{\sum_{i=1}^N \lambda_i^2} \leq \sum_{i=1}^N \abs{\lambda_i}^2 \leq N \rho(\mat{J})^2 \quad \Rightarrow \quad \rho(\mat{J}) \geq \sqrt{\frac{1}{N} \abs{\text{tr}(\mat{J}^2)}} .
\]
We begin by finding an expression for $\text{tr}(\mat{J}^2)$ in terms of the entries of $\mat{J}$,
\[
\text{tr}(\mat{J}^2) = \sum_{i=1}^N \left(\mat{J}^2\right)_{ii} = \sum_{i=1}^N \sum_{j=1}^N \mat{J}_{ij} \mat{J}_{ji} = \sum_{i=1}^N \mat{J}_{ii}^2 + 2 \sum_{i=1}^N \sum_{j<i} \mat{J}_{ij} \mat{J}_{ji} .
\]
Once again, let $r_j \defn u_j / u_i$, such that according to Appendix \ref{app_linearizations}, the logarithmic scheme has Jacobian entries
\[
\mat{J}_{ij} =
\begin{cases}
- 2\mat{D}_{ij} \dfrac{\log(r_j) - 1 + 1/r_j}{\log^2(r_j)}, & i\neq j,\\[2.5ex]
- 2\displaystyle\sum_{k\neq i}\mat{D}_{ik}\dfrac{r_k - \log(r_k)-1}{\log^2(r_k)}, & i=j.
\end{cases}
\]
We begin by examining the cross terms $\mat{J}_{ij} \mat{J}_{ji}$. Using $\mat{D}_{ij} = - \mat{D}_{ji}$, we find
\begin{align*}
\mat{J}_{ij} \mat{J}_{ji} &= - 4 \mat{D}_{ij}^2 \left( \dfrac{\log(r_j) - 1 + 1/r_j}{\log^2(r_j)} \right) \left( \dfrac{\log(1/r_j) - 1 + r_j}{\log^2(1/r_j)} \right) \\
&= - 4 \mat{D}_{ij}^2 \frac{\left( r_j \log(r_j) - r_j + 1 \right)\left( r_j - \log(r_j) - 1 \right)}{r_j \log^4(r_j)} .
\end{align*}
For $r_j \neq 1$, each term in brackets in the numerator is a positive function, meaning $\mat{J}_{ij} \mat{J}_{ji} < 0$. Furthermore, by following the same asymptotic scaling arguments as before, we can show that as $r_j \rightarrow \infty$ or $r_j \rightarrow 0$, $\mat{J}_{ij} \mat{J}_{ji} \rightarrow - \infty$. For any other bounded value of $ 0 < r_j < \infty$, $\mat{J}_{ij} \mat{J}_{ji}$ is bounded. Therefore, if the diagonal components $\mat{J}_{ii}$ are bounded, we immediately establish that $\abs{\text{tr}(\mat{J}^2)}$, and therefore $\rho(\mat{J})$, is unbounded. Unfortunately, $\mat{J}_{ii}$ can indeed be unbounded, as we have already shown that the individual sums of the diagonal entries, $a_k$, can grow unboundedly as $r_k \rightarrow \infty$. Furthermore, cancellations between $a_k$ make it difficult to predict how $\mat{J}_{ii}$ behaves overall, so instead of establishing a result in full generality, we restrict ourselves to a few special cases to only show that, at very least, a blowup of spectral radius is possible.

Assume $\mat{D}$ is a second-order central difference operator, and that ${u_i = \epsilon\rightarrow 0}$ while ${u_{i+1} = \delta_{+} h + \epsilon}$ and ${u_{i-1} = \delta_{-} h + \epsilon}$ for some $\delta_{\pm} > 0$ that depends only on the physical solution profile~$\fnc{U}(x)$. Then $\mat{D}_{i,i+1} = \sfrac{1}{2h}$ and $\mat{D}_{i,i-1} = \sfrac{-1}{2h}$ and $r_{i+1} , r_{i-1} \rightarrow \infty$. Consequently,
\[
\mat{J}_{ii} = \frac{1}{h} \left( \frac{r_{i-1} - \log(r_{i-1}) - 1}{\log^2(r_{i-1})} - \frac{r_{i+1} - \log(r_{i+1}) - 1}{\log^2(r_{i+1})} \right) .
\]
Both terms in brackets approach $\infty$ as $\epsilon \rightarrow 0$ (and $u_i = \epsilon \rightarrow 0$), but they do so at slightly different rates. If $\delta_{+} = \delta_{-}$, then $r_{i+1} = r_{i-1}$, so they cancel and $\mat{J}_{ii} = 0$ (though as shown below, $\mat{J}_{ij} \mat{J}_{ji}$ still blows up). However if $\delta_{+} \neq \delta_{-}$, then an asymptotic expansion about $\epsilon$ in the limit $\epsilon \ll \delta_{\pm}$ shows
\[
\mat{J}_{ii} \sim  \frac{\delta_{-}-\delta_{+}}{\epsilon \log^2(1/\epsilon)} \rightarrow \text{sgn}(\delta_{-}-\delta_{+}) \cdot \infty \quad \text{as} \ \epsilon \rightarrow 0 .
\]
A similar asymptotic expansion of $\mat{J}_{ij} \mat{J}_{ji}$ for $j=i-1$ yields
\[
\mat{J}_{ij} \mat{J}_{ji} \sim - \frac{1}{h} \frac{\delta_{-}}{\epsilon \log^3(1/\epsilon)} \rightarrow - \infty \quad \text{as} \ \epsilon \rightarrow 0 .
\]
Although the signs of $\mat{J}_{ii}^2$ and $\mat{J}_{ij} \mat{J}_{ji}$ can be opposite, note that $\mat{J}_{ii}^2$ grows quicker than $\mat{J}_{ij} \mat{J}_{ji}$, and therefore dominates. As a result, $\abs{\text{tr}(\mat{J}^2)} \rightarrow \infty$ as $\epsilon \rightarrow 0$, and so $\rho(\mat{J}) \rightarrow \infty$ also.

\section{Floquet Theory} \label{app_floquet}

\subsection{The Monodromy Matrix and Floquet Multipliers}

This section summarizes relevant information taken from \cite{Chicone2006}. Consider the time-varying linear system
\begin{equation} \label{eq_time_periodic_system2}
    \der[\bu]{t} = \mat{J}(t) \bu , \quad \bu \in \mathbb{R}^N , \quad \mat{J}(t) = \mat{J}(t + T) \in \mathbb{R}^{N \times N} .
\end{equation}
A \emph{fundamental matrix} $\mat{\Phi}(t) \in \mathbb{C}^{N \times N}$ is a matrix whose columns form a linearly independent basis of solutions to \eqref{eq_time_periodic_system2}. 
The \emph{principal fundamental matrix} is defined by
\begin{equation} \label{eq_principal_fundamental_matrix}
    \mat{\Psi}(t) \defn \mat{\Phi}(t)\mat{\Phi}^{-1}(0) ,
\end{equation}
so that $\mat{\Psi}(0)=\mat{I}$, and every solution $\bu$ of \eqref{eq_time_periodic_system2} satisfies
\[
\bu(t) = \mat{\Psi}(t)\bu(0).
\]
We can view $\mat{\Psi}(t)$ as a matrix of solutions with initial conditions $\vec{e}_i$ in column $i$, or equivalently, the evolution operator acting on an arbitrary initial condition $\bu_0 = \bu(0)$. We can interchangeably call $\mat{\Psi}(t)$ the \emph{state transition matrix}, though this name usually denotes the family of principal fundamental matrices parametrized by different initial times $t_0$ (we assumed $t_0 = 0$ above).
Periodicity of $\mat{J}(t)$ implies that
\[
\mat{\Phi}(t+T) = \mat{\Phi}(t) \mat{\Phi}^{-1}(0) \mat{\Phi}(T) , \quad \text{or} \quad \mat{\Psi}(t+T) = \mat{\Psi}(t)\mat{\Psi}(T).
\]
The matrix $\mat{\Psi}(T)$ is called the \emph{monodromy matrix}. It is invertible, and its eigenvalues $\rho \in \mathbb{C}$ are the \emph{characteristic multipliers} of \eqref{eq_time_periodic_system2}. They are unique, \ie independent of the initial time. 
The Floquet theorem states that there exists a constant matrix $\mat{B} \in \mathbb{C}^{N \times N}$ and a $T$-periodic matrix $\mat{P}(t) \in \mathbb{C}^{N \times N}$ such that
\[
\mat{\Phi}(t) = \mat{P}(t)e^{\mat{B} t},
\qquad 
\mat{P}(t+T)=\mat{P}(t).
\]
In other words, each matrix of solutions---including the principal funcamental matrix $\Psi(t)$---can be expressed as the solutions to a constant-coefficient linear system in a rotated frame via a change of coordinates $\bu \leftrightarrow \tilde{\bu}$,
\[
\bu(t) = \mat{P}(t) \tilde{\bu}(t) , \quad \der[\tilde{\bu}]{t} = \mat{B} \tilde{\bu} , \quad \text{or} \quad \tilde{\bu}(t) = e^{\mat{B} t} \tilde{\bu}_0 .
\]
This representation of $\mat{\Phi}$ is the \emph{Floquet normal form}, and allows us to express the monodromy matrix as $\mat{\Psi}(T) = e^{\mat{B} T}$, since this particular choice of $\mat{\Phi}$ corresponds to $\mat{P}(0) = \mat{I}$. Therefore, the characteristic multipliers are also eigenvalues of $e^{\mat{B} T}$, a constant-coefficient system in the rotated frame.
A complex number $\lambda \in \mathbb{C}$ is called a \emph{Floquet exponent} if $e^{\lambda T} = \rho$ for some characteristic multiplier $\rho$, \ie $\lambda = \tfrac{1}{T} \log(\rho)$. Note that the exponents are not unique, since $\lambda + 2\pi i k/T$, $k \in \mathbb{Z}$ yields the same multiplier $\rho$. Up to this imaginary shift, they correspond to the eigenvalues of $\mat{B}$. Note that if the principal fundamental matrix is known, we can directly express $\mat{P}(t)$ and $\mat{B}$ as
\[
    \mat{P}(t) = \mat{\Psi}(t) e^{-\mat{B} t} , \quad \mat{B} = \tfrac{1}{T} \log(\mat{\Psi}(T)) .
\]

Stability over one period is determined by the moduli of the multipliers:
\[
\abs{\rho} = e^{\Re(\lambda)T}, \quad \abs{\rho} \leq 1 \ \text{for stability}.
\]
Equivalently, $\Re(\lambda)\le 0$ characterizes non-growth over each period. Note that we assume here that all eigenvalues on the unit circle are semisimple; if they are not then $\abs{\rho} \leq 1$ is not sufficient for boundedness. In the special case where $\mat{J}$ is time-invariant, we recover $\mat{P}(t)=\mat{I}$ and $\mat{B}=\mat{J}$, so the Floquet exponents reduce to the eigenvalues of $\mat{J}$. In this sense, Floquet theory generalizes eigenvalue analysis of constant-coefficient systems to periodic systems. Moreover, in the periodic setting, the real parts of the Floquet exponents coincide with the Lyapunov exponents of the system.

To characterize the worst-case scenario growth for long times, we examine the spectral radius of $\mat{\Psi}(T)$, or equivalently, $\max_j \abs{\rho_j}$. The mode with largest asymptotic growth is the eigenvector of $\mat{\Psi}(T)$ corresponding to the characteristic multiplier (eigenvalue of $\mat{\Psi}(T)$) of maximal modulus, called the \emph{dominant Floquet mode}. In general, if $\mat{V}(0)$ is the right eigenvector matrix of $\mat{\Psi}(T)$, and $\mat{M}$ is the diagonal matrix of characteristic multipliers, then by definition, $\mat{V}(T) = \mat{\Psi}(T) \mat{V}(0) = \mat{V}(0) \mat{M}$. It is then possible to show that the Floquet normal form of the monodromy matrix is
\[
\mat{\Psi}(T) = e^{\mat{B} T} , \quad \mat{B} = \mat{V}(0) \mat{\Lambda} \mat{V}^{-1}(0) ,
\]
where $\mat{\Lambda}$ is the diagonal matrix of Floquet exponents. Therefore, $\mat{V}(0)$ is also the eigenvector matrix of $\mat{B}$ with eigenvalue matrix $\mat{\Lambda}$, and since $\mat{B}$ is constant for all $t$, we obtain
\[
\mat{V}(t) = \mat{\Psi}(t) \mat{V}(0) = \mat{P}(t) e^{\mat{B} t} \mat{V}(0) = \mat{P}(t) \mat{V}(0) e^{\mat{\Lambda} t} , \quad \mat{P}(t+T) = \mat{P}(t).
\]
Therefore, Floquet modes $\bv = e^{\lambda t} \vec{p}$ can be separated into a growth term $e^{\lambda t}$ multiplying a periodic part, $\vec{p} (t) \defn \mat{P}(t) \bv_0 = e^{-\lambda t} \mat{\Psi}(t) \bv_0$, where $\vec{p} (t+T) = \vec{p} (t)$.

Within a single period (or even possibly after a few periods), significant growth is possible even if $\max_j \abs{\rho_j} \leq 1$ due to transient dynamics.
Short-time amplification is governed by
\[
\norm{\mat{\Psi}(t)} = \sup_{\norm{\bu}=1} \norm{\mat{\Psi}(t) \bu}
= \sigma_{\max}(\mat{\Psi}(t)) ,
\]
where $\sigma_{\max}$ indicates the maximum singular value. This holds for any weighted Euclidean norm (\eg {$\Hnrm$-norm}), provided that the singular value is computed using a scaling transform (\eg $\sigma_{\max}(\Hnrm^{1/2} \mat{\Psi}(t) \Hnrm^{-1/2})$). If $\sigma_{\max}(\mat{\Psi}(t))>1$, there exists initial data that experiences amplification up to time $t$ due to nonnormality of $\mat{\Psi}(t)$. The mode achieving the largest possible amplification over this time is the right singular vector corresponding to $\sigma_{\max}(\mat{\Psi}(t))$.

\subsection{Numerical Computation of the Monodromy Matrix}

One can numerically compute the monodromy matrix and Floquet multipliers using a variety of different techniques, and for a complete discussion, we refer the reader to the rich literature, \eg \cite{magnus, Iserles_lie_groups, floquet2, floquet1, friedmann_1977}. Below we introduce one simple method that is well suited to our needs: the exponential midpoint method.

We begin by drawing attention to the composition property of time-varying linear systems \eqref{eq_time_periodic_system2},
\[
\mat{\Psi}_{t_0}(t) = \mat{\Psi}_{t_1}(t) \mat{\Psi}_{t_0}(t_1) ,
\]
where $\mat{\Psi}_{t_0}(t) \defn \mat{\Phi}(t) \mat{\Phi}^{-1}(t_0)$ denotes the fundamental matrix (or evolution operator) taking a solution $\bu(t_0)$ defined at some initial time $t_0$ to the solution $\bu(t)$ at later time $t \geq t_0$. In the previous section, we suppressed the subscript $t_0$. The composition property dictates that we can break the interval $[0,T]$ into $K$ subintervals $[t_{j-1}, t_j] \subset [0,T]$ with $0 = t_0 < t_1 < \dots t_K = T$, so that the monodromy matrix can be expressed as 
\begin{equation} \label{eq_monodromy_product}
\mat{\Psi}(T) =  \mat{\Psi}_{t_{K-1}}(t_K) \dots \mat{\Psi}_{t_1}(t_2) \mat{\Psi}_{t_0}(t_1) 
= \prod_{j=1}^K \mat{\Psi}_{t_{j-1}}(t_{j}) ,
\end{equation}
where the matrix product is understood to be in the correct order (lower $t_j$ applied first). This is analogous to taking a step of a time marching method. On each step, however, the principal fundamental matrix $\mat{\Psi}_{t_{j-1}}(t_{j})$ remains to be calculated. One simple method would be to take $N$ linearly independent solutions, propagate them in time through $[t_{j-1},t_j]$ using a numerical method, then assembling $\mat{\Psi}_{t_{j-1}}(t_{j})$ through \eqref{eq_principal_fundamental_matrix}. However, this will be susceptible to errors both from the numerical discretization, and possible ill-conditioning of $\mat{\Phi}$ if the columns are individually propagated to similar solutions. Instead, one can consider series expansions of $\mat{\Psi}_{t_{j-1}}(t_{j})$ through either a Magnus expansion \cite{magnus_orig, magnus}, a Fer expansion \cite{fer}, or a Peano-Baker series (also known as the Dyson series) \cite{dyson1, dyson2}. 

The need for such expansions comes from the fact that only in the very restrictive scenario that the linear operator $\mat{J}$ commutes for all times can we express the solution $\bu$ as a proper exponential (as we do in the scalar case), \ie
\[
\bu(t) = e^{\int_0^t \mat{J}(s) \mathrm{d}s} \bu_0 \quad \text{only if} \quad \mat{J}(t_1) \mat{J}(t_2) = \mat{J}(t_2) \mat{J}(t_1) \ \forall t_1, t_2 \in [0,t] .
\]
In general, the exponential integral must be either interpreted as an infinite sum of time-ordered integrals (Peano-Baker series) or a sum of increasingly unwieldy commutator corrections (Magnus expansion). To the first term in the series, however, the Magnus expansion (as well as the Fer expansion) reduces to the same approximation written above \cite{magnus, Iserles_lie_groups},
\[
\mat{\Psi}_{t_{j-1}}(t_j) \approx e^{\int_{t_{j-1}}^{t_j} \mat{J}(s) \mathrm{d} s} .
\]
Truncating the Magnus expansion in this way yields a second-order accurate exponential integrator, provided that $\mat{J}(t)$ is sufficiently smooth and the above integral is approximated with (at least) a second-order quadrature \cite{Iserles_lie_groups}. We use the midpoint rule, resulting in the ``exponential midpoint'' method \cite{magnus},
\[
\mat{\Psi}_{t_{j-1}}(t_j) \approx e^{\mat{J}(t_{j-1/2}) \Delta t} , \quad t_{j-1/2} \defn \frac{t_{j-1} + t_j}{2} , \quad \Delta t \defn t_j - t_{j-1} .
\]
This method is exact for constant $\mat{J}$, and if $\mat{J}$ is skew-symmetric, it produces orthogonal $\mat{\Psi}_{t_{j-1}}(t_j)$, thus preserving solution norms~\cite{Iserles_lie_groups}. Convergence to the exact monodromy matrix is achieved by increasing the number of subintervals $K$. Occasionally, if the subinterval fundamental matrices are ill-conditioned, the composition of $K$ matrices in \eqref{eq_monodromy_product} can result in an approximation to $\mat{\Psi}(T)$ that yields significant errors when applying an eigenvalue routine to compute the Floquet multipliers. In this case, it may be preferable to apply a periodic Schur decomposition to avoid explicitly computing the matrix product \cite{floquet1}. For the purposes of this paper, however, we have found the simple exponential midpoint method to be sufficiently accurate in constructing well-conditioned and convergent approximations to $\mat{\Psi}(T)$ with relatively small $K$.

\section{Linearizing Burgers to Variable-Coefficient Advection} \label{app_lin_burgers}

The linearization of the inviscid Burgers equation to recover the variable-coefficient advection equation is a standard procedure~\cite{kreiss_lorenz, gassner_stability_2022}, though for completeness, we repeat the process here. Consider the inviscid Burgers equation
\begin{equation} \label{eq_burgers}
    \pder[\fnc{U}]{t} + \pder{x} \left( \frac{1}{2} \fnc{U}^2 \right) = 0
\end{equation}
on the domain $(x,t) \in [0,L] \times [0,T] \subset \mathbb{R} \times \mathbb{R}^+$ with periodic boundary conditions and some initial condition $\fnc{U}(x,t=0)=\fnc{U}_0(x)$. Assuming smoothness, we decompose $\fnc{U}$ into a smooth and periodic \emph{baseflow} $\tilde{\fnc{U}}$ plus a small \emph{perturbation} $\fnc{V}$~\cite{kreiss_lorenz}. Plugging this into \eqref{eq_burgers}, we find
\begin{equation*}
    \pder[\fnc{V}]{t} + \tilde{\fnc{U}} \pder[\fnc{V}]{x} + \pder[\tilde{\fnc{U}}]{x} \fnc{V} + \fnc{V} \pder[\fnc{V}]{x} = - \left( \pder[\tilde{\fnc{U}}]{t} + \tilde{\fnc{U}} \pder[ \tilde{\fnc{U}}]{x} \right)
\end{equation*}
where the initial conditions are $\fnc{V}(x,t=0) = \fnc{V}_0$ and $\tilde{\fnc{U}}(x,t=0)=\tilde{\fnc{U}}_0(x)$, and the boundary conditions are again periodic. Suppose the decomposition was made such that $\tilde{\fnc{U}}$ is an exact solution to \eqref{eq_burgers} with the slightly perturbed initial condition $\tilde{\fnc{U}}(x,t=0)=\tilde{\fnc{U}}_0(x)$. This implies that the right-hand side source term is zero. Furthermore, if the difference between the initial problem and the perturbed problem is small, then $\fnc{V}$ is small, and the right-most term on the left-hand side can be dropped as a higher order perturbation~\cite{kreiss_lorenz}. We obtain the linearized problem
\begin{equation} \label{eq_burgers_lin}
    \pder[\fnc{V}]{t} + \pder{x} \left( \tilde{\fnc{U}} \fnc{V} \right) = 0,
\end{equation}
which governs how small perturbations $\fnc{V}$ to solutions of \eqref{eq_burgers} behave. This is precisely the variable advection equation \eqref{eq_lin_conv} that we have studied in this work, with $a(x,t) = \tilde{\fnc{U}}(x,t)$. Notably, however, the baseflow $\tilde{\fnc{U}}$ is time-varying, and so the frozen baseflow experiments performed in \S\ref{sec:spectra} and \S\ref{sec:time_invariant_accuracy_experiments} with a time-invariant $a(x)$ represent a worst-case scenario. In practice, perturbations do not experience consistent exponential modal growth, but instead alternating cycles of growth and dissipation as the baseflow evolves with the perturbation. Another important distinction to the linear case is that the baseflow $\tilde{\fnc{U}}(x,t)$ develops a shock in short-time, which limits the period of time within which the linearization picture is valid~\cite{kreiss_lorenz}.

In a similar fashion, we can also show that linearizing a split-form semidiscretization of the Burgers equation recovers the split-form product-scheme semidiscretization \eqref{eq_general_linear_split_form} we studied extensively in \S\ref{sec_lin_advec}. Consider the following consistent and conservative split-form semidiscretization for \eqref{eq_burgers} with $\alpha \in [0,1]$,
\begin{equation} \label{burgers_disc}
    \der[\bu]{t} + \alpha \mat{D} \left( \tfrac{1}{2} \mat{U} \bu \right) + (1 - \alpha ) \mat{U} \mat{D} \bu   = \text{SAT} , 
    \quad \mat{U} \defn \diag{\bu} .
\end{equation}
The following simultaneous approximation term (SAT) ensures the scheme is conservative,
\begin{align} \label{burgers_disc_sat}
    \mathrm{SAT} &= \HI \Big( \vec{e}_N \tfrac{1}{2} \left( \tfrac{1}{2} u_N^2 + \tfrac{1}{2} (u_1^R)^2 \right) - \vec{e}_1 \tfrac{1}{2} \left( \tfrac{1}{2} u_1^2 + \tfrac{1}{2} (u_N^L)^2 \right) \notag \\
    &\qquad + \frac{(1 - \alpha)}{2} \vec{e}_N \left( \tfrac{1}{2} u_N^2 + \tfrac{1}{2} (u_1^R)^2 - u_N u_1^R \right) - \frac{(1 - \alpha)}{2} \vec{e}_1 \left( \tfrac{1}{2} u_1^2 + \tfrac{1}{2} (u_N^L)^2 - u_1 u_N^L \right) \\
    &\qquad - \frac{\sigma}{2} \vec{e}_N \abs{ u_N^\star} \left( u_N - u_1^R \right) - \frac{\sigma}{2} \vec{e}_1 \abs{ u_1^\star} \left( u_1 - u_N^L \right) \Big) \notag ,
\end{align}
where $\sigma$ is some upwinding parameter that introduces interface dissipation, and $\abs{ u_i^\star}$ is some approximate intermediate state. $u^R_1$ and $u^L_N$ denote the first and last nodes of the right and left adjacent blocks, respectively. The choice $\alpha=1$ recovers a standard central (or conservative-form) semidiscretization of \eqref{burgers_disc}, while $\alpha=\tfrac{2}{3}$ recovers the well-known energy-stable (or entropy-stable) semidiscretization~\cite{DelReyFernandez2014_review}. Continuous SBP and circulant finite-difference schemes are recovered by setting $\text{SAT} = 0$. Linearizing \eqref{burgers_disc} through $\bu \rightarrow \tilde{\bu} + \bv$ yields a consistent split-form semidiscretization of \eqref{eq_burgers_lin},
\begin{equation} \label{eq_Burgers_disc_lin}
    \der[\bv]{t} + \alpha \mat{D} \tilde{\mat{U}} \bv + (1-\alpha)  \tilde{\mat{U}} \mat{D} \bv + (1-\alpha) \mat{V} \mat{D} \tilde{\bu} = \mathrm{SAT}', \quad \tilde{\mat{U}} \defn \diag{\tilde{\bu}}, \ \mat{V} \defn \diag{\bv} ,
\end{equation}
where the linearization of the SAT term can be complicated if $\sigma \neq 0$ due to the possibly discontinuous nature of $\abs{u_i^\star}$. For simplicity, here we state only the linearization only for $\sigma = 0$,
\begin{equation} \label{eq_Burgers_disc_lin_sat}
\begin{aligned}
    \mathrm{SAT}' &= \tfrac{1}{2} \HI \Big( \vec{e}_N \left( \tilde{u}_N v_N + \tilde{u}_1^R v_1^R \right) - \vec{e}_1 \left( \tilde{u}_1 v_1 + \tilde{u}_N^L v_N^L \right) \\
    & \qquad \qquad + (1 - \alpha) \Big( \vec{e}_N \left( \tilde{u}_N v_N + \tilde{u}_1^R v_1^R - \tilde{u}_N v_1^R - \tilde{u}_1^R v_N \right) \\
    & \qquad \qquad \qquad \qquad - \vec{e}_1 \left( \tilde{u}_1 v_1 + \tilde{u}_N^L v_N^L - \tilde{u}_1 v_N^L - \tilde{u}_N^L v_1 \right) \Big) \Big) .
\end{aligned}
\end{equation}
If $\tilde{\bu}$ is continuous across interfaces, the second line vanishes and we obtain the central-product split-form semidiscretization \eqref{eq_general_linear_split_form} for the variable-coefficient linear advection problem with $a(x_i,t) = \tilde{u}_i$.

\section{A Stability Bound on the Conservative Variables From Entropy Stability} \label{app_entropy_l2_bound}

In this section we consider an arbitrary system of conservation laws $\vecfnc{U} \in \mathbb{R}^n$, and show that any discretely conservative semidiscretization that has a non-increasing entropy function has a bounded numerical solution, provided that certain assumptions are met. These proofs are based on the relative entropy method of \cite{Dafermos2010, svard_2015}, though here we follow \cite{sbpbook} and explicitly consider the semidiscrete case. 

\subsection{Explicit Construction of the Semidiscrete Stability Bound}

Consider any entropy-stable and discretely conservative semi-discretization (\eg the flux-differencing formulations of \cite{fisher_2013, crean}) that can be written as 
\begin{align*}
\der[\bu]{t} = \vecfnc{R}(\bu) , \quad &\text{with} \quad
 \der[\fnc{S}_h(\bu)]{t} = \sum_{i=1}^N \Hnrm_{ii} \vec{w}_i^\T \der[\bu_i]{t} \leq 0 \\
\quad & \ \text{and} \quad 
\sum_{i=1}^N \Hnrm_{ii} \vec{\alpha}^\T \der[\bu_i]{t} = 0 \quad  \forall \ \text{constant} \ \vec{\alpha} \in \mathbb{R}^n,
\end{align*}
where $\bu \in \mathbb{R}^{Nn}$ and $\bu_i \in \mathbb{R}^n$ at each node $i$, $\fnc{S}_h(\bu) \defn \sum_{i=1}^N \mat{H}_{ii} \fnc{S}(\bu_i)$ is a discrete analogue of the total entropy $\int_\Omega \fnc{S}(\vecfnc{U}) \, \mathrm{d}x$, ${\fnc{S} : \fnc{A} \subset \mathbb{R}^n \rightarrow \mathbb{R}}$ is the entropy function of the continuous system, and $\vec{w} \in \mathbb{R}^{Nn}$ are the discrete entropy variables ${\vec{w}_i = \vecfnc{W}(\bu_i) \in \mathbb{R}^n}$, where $\vecfnc{W}^\T \defn \partial \fnc{S} / \partial \vecfnc{U}$. The \emph{admissible set} $\fnc{A} \subset \mathbb{R}^n$ is the set of all possible states $\vecfnc{U} \in \mathbb{R}^n$ for which
\begin{itemize}
    \item the flux $\vecfnc{F}(\vecfnc{U})$ is well defined, and
    \item the entropy $\fnc{S}(\vecfnc{U})$ is well defined and $C^2$, and \item the Hessian of the entropy, $\nabla^2 \fnc{S}(\vecfnc{U})$, is positive definite, \ie $\vecfnc{V}^\T \nabla^2 \fnc{S}(\vecfnc{U}) \vecfnc{V} > 0 \quad \forall \ \vecfnc{V} \in \mathbb{R}^{n}, \ \vecfnc{U} \in \fnc{A}$.
\end{itemize}
The third condition ensures strict convexity of $\fnc{S}$ on convex subsets of $\fnc{A}$. We must now make the following assumption, the importance of which will be discussed further in the following section.
\begin{assumption} \label{assumption_admissible}
At each node $i$, the semidiscrete solution $\bu_i \in \mathbb{R}^n$ remains in a compact convex subset $\fnc{K}$ of the interior of the admissible set, \ie $\fnc{K} \subset \mathrm{int}(\fnc{A})$.
\end{assumption}

To construct a bound on $\norm{\bu}$, we begin by letting $\bar{\vec{u}}$ be some arbitrary constant reference state with $\bar{\vec{u}}_i = \bar{\vecfnc{U}} \in \fnc{K}$, and introducing the discrete relative entropy
\[
\tilde{\fnc{S}_h}(\bu) \defn \fnc{S}_h(\bu)-\fnc{S}_h(\bar{\bu}) - \sum_{i=1}^N \Hnrm_{ii} \bar{\vec{w}}_i^\T (\bu_i-\bar{\bu}_i),
\]
where $\bar{\vec{w}} \defn \vecfnc{W}(\bar{\bu})$. We now establish several intermediate results on $\tilde{\fnc{S}}_h$ before attempting to bound $\norm{\bu}$. First we note that $\tilde{\fnc{S}}_h$ is bounded from above. Since $\bar{\bu}$ is constant, the relative entropy satisfies
\[
\der[\tilde{\fnc{S}}_h(\bu)]{t} = \der[\fnc{S}_h(\bu)]{t} - \sum_{i=1}^N \Hnrm_{ii} \bar{\vec{w}}_i^\T \der[\bu_i]{t} 
\leq - \vecfnc{W}(\bar{\vecfnc{U}})^\T \sum_{i=1}^N \Hnrm_{ii} \der[\bu_i]{t} = 0 ,
\]
where the inequality comes from the entropy stability estimate, and the final equality comes from discrete conservation. Therefore, the relative entropy can be bounded from above via $\tilde{\fnc{S}}_h(\bu(t)) \leq C \defn \tilde{\fnc{S}}_h(\bu_0)$.

Next we note that since $\tilde{\fnc{S}}_h$ differs from $\fnc{S}_h$ only by an affine function of
$\bu$, it has the same Hessian $\nabla^2 \tilde{\fnc{S}}_h(\bu)=\nabla^2 \fnc{S}_h(\bu)$, which is block diagonal with blocks between entries $ni$ and $n(i+1)$ given by 
\[
\mat{H}_{ii} \nabla^2 \fnc{S}(\bu_i) \in \mathbb{R}^{n \times n}.
\]
Admissibility $\bu_i \in \mathrm{int}(\fnc{A})$ ensures that each block is continuous in $\bu_i$ and symmetric positive definite, implying that the full Hessian $\nabla^2 \tilde{\fnc{S}}_h(\bu)$ is also continuous and symmetric positive definite.

The last important property to note is that $\tilde{\fnc{S}}_h$ achieves its minimum $\tilde{\fnc{S}}_h(\bar{\bu}) = 0$ at $\bu = \bar{\bu}$, since 
\[
\left[ \nabla \tilde{\fnc{S}}_h(\bu) \right]_j = 
\sum_{i=1}^N \mat{H}_{ii} \pder[\fnc{S}(\bu_i)]{\bu_j} - \sum_{i=1}^N \Hnrm_{ii} \bar{\vec{w}}_i^\T \pder[\bu_i]{\bu_j} = \mat{H}_{jj} (\vec{w}_j - \vecfnc{W}(\bar{\vecfnc{U}}))^\T = \vec{0}^\T \quad \text{if} \quad \vec{w} = \bar{\vec{w}}.
\]
Expanding about this minimum point $\bar{\bu}$, Taylor's theorem then allows us to express any $\tilde{\fnc{S}}_h(\bu)$ as
\[
\tilde{\fnc{S}}_h(\bu) = \tfrac{1}{2} (\bu - \bar{\bu})^\T \nabla^2 \fnc{S}_h(\vec{z}) (\bu - \bar{\bu}) , \quad \vec{z} \defn \bar{\bu} + \theta (\bu - \bar{\bu}) \ \  \text{for some} \ \theta \in [0,1] .
\]
By Assumption \ref{assumption_admissible}, $\bar{\bu}_i, \bu_i \in \fnc{K}$, so by convexity of $\fnc{K}$, $\vec{z}_i \in \fnc{K}$ also. Furthermore, since $\fnc{K}$ is compact and $\nabla^2 \fnc{S}_h(\bu)$ is both continuous and symmetric positive definite, there exists a constant $\beta > 0$ such that
\[
\bv^\T \nabla^2 \fnc{S}_h(\bu) \bv \geq \beta \norm{\bv}_\Hnrm^2 \quad \forall \ \bv \in \mathbb{R}^{Nn}, \ \bu_i \in \fnc{K}.
\]
This follows from continuous functions always attaining their minimums on compact sets. For example, a sufficient uniform positive definiteness constant is 
\begin{equation} \label{eq_coercivity_const}
\beta = \min_{\vecfnc{U} \in \fnc{K}} \lambda_{\min} \! \left( \nabla^2 \fnc{S}(\vecfnc{U}) \right)  .
\end{equation}
Returning to our Taylor expansion of $\tilde{\fnc{S}}_h(\bu)$, we find
\[
\tilde{\fnc{S}}_h(\bu_0) = C \geq \tilde{\fnc{S}}_h(\bu) = \tfrac{1}{2} (\bu - \bar{\bu})^\T \nabla^2 \fnc{S}_h(\vec{z}) (\bu - \bar{\bu}) \geq \tfrac{1}{2} \beta \norm{\bu - \bar{\bu}}_\Hnrm^2 .
\]
Finally, using the norm relations $\norm{\bu}^2 \leq \left( \norm{\bu - \bar{\bu}} + \norm{\bar{\bu}} \right)^2 \leq 2 \norm{\bu - \bar{\bu}}^2 + 2 \norm{\bar{\bu}}^2$, we obtain 
\begin{equation} \label{eq_entropy_l2_bound}
\norm{\bu}_\Hnrm^2 \leq \frac{4 C}{\beta} + 2 \norm{\bar{\bu}}_\Hnrm^2 .
\end{equation}

\subsection{A Discussion on the Validity of the Admissibility Assumptions}

In order to construct the stability bound \eqref{eq_entropy_l2_bound}, we required Assumption \ref{assumption_admissible}. Below we summarize where each part of this assumption was used:
\begin{itemize}
    \item $\fnc{K}$ lies in the interior of the admissible set $\fnc{A}$ because we need an open set to ensure the Hessian $\nabla^2 \fnc{S}(\vecfnc{U})$ is well defined, and that Taylor's theorem applies.
    \item Admissibility $\bu_i \in \mathrm{int}(\fnc{A})$ ensures the Hessian is everywhere continuous and symmetric positive definite.
    \item Convexity of $\fnc{K}$ allows us to express $\tilde{\fnc{S}}_h(\bu)$ in a quadratic form evaluated at a state that remains in $\fnc{K}$.
    \item Compactness of $\fnc{K}$ and $\fnc{S}$ being $C^2$ allows the pointwise positive definiteness of the Hessian to be upgraded to uniform positive definiteness, ensuring the existence of some $\beta > 0$.
\end{itemize}

Given a specific system of conservation laws, the admissible set $\fnc{A}$ is often straightforward to identify. For example, for the Euler and Navier--Stokes equations with $\fnc{S} = -\rho s/(\gamma-1)$ and $s \defn \log(p/\rho^\gamma)$, smoothness and positive definiteness of $\nabla^2 \fnc{S}(\vecfnc{U})$ is guaranteed under positivity of density $\rho$ and pressure $p$ \cite{svard_2015},
\[
\fnc{A}_{\text{Euler}} \defn \left\lbrace ( \rho, \rho \vec{v}, e) : \rho > 0 , \ p \defn ( \gamma-1)\left(e - \tfrac{1}{2} \rho \norm{\vec{v}}^2\right) > 0 \right\rbrace ,
\]
which forms a convex set \cite{zhang_shu}. For the geometric and logarithmic flux-differencing semidiscretizations of the variable-coefficient linear advection equation discussed in \S\ref{sec:flux_differencing}, where $\fnc{S}_\text{geom}(\fnc{U}) \defn -2 \sqrt{\fnc{U} / a}$ and $\fnc{S}_\text{log}(\fnc{U}) \defn \fnc{U} \log(a\fnc{U}) - \fnc{U}$, respectively, the admissible sets are even simpler,
\[
\fnc{A}_{\text{geom}} = \fnc{A}_{\text{log}} \defn \left\lbrace \fnc{U} : 0 < \fnc{U} < \infty \right\rbrace .
\]
since $\fnc{S}_\text{geom}''(\fnc{U}) = 1 / (2 \sqrt{a} \, \fnc{U}^{3/2})$ and $\fnc{S}_\text{log}''(\fnc{U}) = 1 / \fnc{U}$. Intervals on $\mathbb{R}$ trivially form convex sets.

Showing that Assumption \ref{assumption_admissible} is valid remains a difficult task, though in some cases, it is possible to establish partial guarantees. For example, we recognize that the discrete entropy stability bound naturally restricts the solution $\bu$ to a sublevel subset of $\mathbb{R}^{Nn}$, given by
\[
\fnc{K}_{\bu} \defn \left\lbrace \bu : \fnc{S}_h(\bu) \defn \sum_{i=1}^N \mat{H}_{ii} \fnc{S}(\bu_i) \leq \fnc{S}_h(\bu_0)  \right\rbrace \subset \mathbb{R}^{Nn} .
\]
This is still insufficient to guarantee that $\bu_i$ remain in compact subsets of $\mathrm{int}(\fnc{A})$, however. For this we need additional conditions to be met. For example, if we can show that $\fnc{S}_h(\bu)$ is coercive, \ie $\fnc{S}_h(\bu) \rightarrow \infty$ as $\norm{\bu} \rightarrow \infty$, then one can deduce that not only $\fnc{K}_{\bu}$ is bounded, but also that any pointwise solution $\bu_i$ is restricted to some bounded subset $\fnc{K} \subset \mathbb{R}^n$, since a pointwise blowup would also imply violation of the discrete entropy stability bound.

For the logarithmic scheme entropy, proving coercivity is straightforward. Temporarily assume $u_i>0$. Each nodal entropy $\fnc{S}_\text{log}(u_i)$ has a unique minimum at $u_i = 1/a_i$ with $\fnc{S}_\text{log}(1/a_i) = - 1/a_i$. Furthermore, they are strictly convex with limits $\fnc{S}_\text{log}(u_i) \rightarrow 0$ as $u_i \rightarrow 0^+$ and $\fnc{S}_\text{log}(u_i) \rightarrow \infty$ as $u_i \rightarrow \infty$. Positivity of $\mat{H}_{ii}$ then ensures that no dramatic cancellations can occur in the summation, establishing coercivity of $\fnc{S}_h(\bu)$, and by extension, pointwise boundedness $u_i \leq u_\text{max} < \infty$. From this, we could for example set the uniform positive definiteness constant \eqref{eq_coercivity_const} as $\beta = 1/u_{\max}$.

On the other hand, the geometric scheme is notably not coercive, since $\fnc{S}_\text{geom}(u_i) \rightarrow - \infty$ as $u_i \rightarrow \infty$. It is still possible to prove the existence of a maximum state $u_\text{max}$, but we must rely on other arguments. For example, for scalar equations, discrete conservation plus positivity yields the mesh-dependent bound
\[
\sum_{i=1}^N \Hnrm_{ii} u_i = M > 0 \ \text{constant} , \quad \Rightarrow \quad u_i \leq \frac{M}{\min_{j \in [1,N]} \Hnrm_{jj}} .
\]
One could justify, therefore, that since the simple mesh-dependent bound presented above is sufficient to establish $L^2$ stability, the entropy-based bound \eqref{eq_entropy_l2_bound} redundant. Although this is technically true, we highlight first that \eqref{eq_entropy_l2_bound} results in a sharper bound if a mesh-independent $\beta$ can be found, and second that this simple approach does not always extend to general systems of conservation laws.

The notable issue with Assumption \ref{assumption_admissible} for the geometric and logarithmic schemes is of course enforcing the left boundary of the admissible set, \ie $u_i > 0$. There is no natural mechanism in place to prevent the numerical solution from approaching vacuum states, or from leaving the admissible set $\fnc{A}$ altogether. This exact issue was shown to occur in \S\ref{sec:nearvacuum}. These schemes must therefore be augmented by some positivity-preservation mechanism that ensures $0 < u_\text{min} \leq u_i$. Once this is true, we can say
\[
u_i \in \fnc{K} \defn \lbrace u : u_\text{min} \leq u \leq u_\text{max} \rbrace \subset \mathrm{int}(\fnc{A}) .
\]
In $\mathbb{R}^n$, compactness is equivalent to closed and bounded, and as previously mentioned, convexity is trivial in $\mathbb{R}$. For these scalar schemes, therefore, the key to satisfying Assumption \ref{assumption_admissible} is the inclusion of a numerical positivity-preservation mechanism.

Unfortunately, for the Euler equations, establishing that Assumption \ref{assumption_admissible} is valid remains a challenge even after positivity preservation is included. This is because even though the numerical solution can be constrained to a convex subset $\fnc{K}$ of the admissible set $\fnc{A}$, it is difficult to guarantee that $\fnc{K}$ is compact. Indeed, the entropy is non-coercive in exactly the same qualitative sense as the geometric scalar entropy. For example, if $\rho$ is bounded but $p \rightarrow \infty$, then $s \rightarrow \infty$, and hence $\fnc{S} \rightarrow - \infty$. In order to establish compactness of $\fnc{K}$, therefore, we must appeal to some other arguments. This challenge has been noted before, with the existence of a uniform positive definiteness constant $\beta$ being left as a conjecture \cite{merriam}.

Here, we provide a simple argument using discrete conservation arguments---similar to those used in the scalar case---to construct pointwise (albeit mesh-dependent) upper bounds on the density, pressure, and internal energy. Let the discrete conservative variables be denoted $\bu_i \defn ( \rho_i , \vec{m}_i, e_i )$, and $\Hnrm \defn \min_{i \in [1,N]} \Hnrm_{ii}$. Suppose that the positivity-preservation mechanism provides strict nodal lower bounds
\[
0 < \rho_{\min} \le \rho_i,
\qquad
0 < p_{\min} \le p_i \defn (\gamma-1)\left( e_i - \frac{1}{2 \rho_i} \norm{\vec{m}_i}^2 \right)
\qquad
\forall i \in [1,N].
\]
By discrete conservation of mass and total energy, the finite quantities
\[
M_\rho \defn \sum_{i=1}^N \Hnrm_{ii} \rho_i > 0,
\qquad
M_e \defn \sum_{i=1}^N \Hnrm_{ii} e_i > 0
\]
remain constant in time, where positivity of total energy follows from
\[
e_i = \frac{p_i}{\gamma-1} + \frac{\norm{\vec{m}_i}^2}{2\rho_i} \geq \frac{p_{\min}}{\gamma-1} > 0 .
\]
As in the scalar case, this immediately yields the mesh-dependent pointwise upper bounds
\[
\rho_i \leq \frac{M_\rho}{\Hnrm_{\min}}, \qquad
e_i \leq \frac{M_e}{\Hnrm_{\min}} .
\]
These pointwise bounds also immediately establish boundedness of pressure and momentum via
\[
p_i \leq (\gamma-1) e_i \leq (\gamma-1) \frac{M_e}{\Hnrm_{\min}} , \quad
\norm{\vec{m}_i} \leq
\sqrt{2 \rho_i \left( e_i - \frac{p_{\min}}{\gamma-1} \right)} \leq
\sqrt{2 \frac{M_\rho}{\Hnrm_{\min}} \left( \frac{M_e}{\Hnrm_{\min}} - \frac{p_{\min}}{\gamma - 1} \right)}.
\]
Therefore every nodal state $\bu_i$ remains in the set
\[
\mathcal{K} \defn \left\lbrace \vecfnc{U} \defn (\rho,\vec{m},e) \in \mathbb{R}^{d+2} :
\rho_{\min} \leq \rho \le \frac{M_\rho}{\Hnrm_{\min}}, \ \ 
\frac{p_{\min}}{\gamma-1} \leq e \leq \frac{M_e}{H_{\min}} , \ \
p(\vecfnc{U}) \geq p_{\min} ,
\right\}.
\]
This set is closed and bounded in $\mathbb{R}^{d+2}$, hence compact. Moreover, it is convex, since the bounds on $\rho$ and $e$ define a convex box $(\rho,e) \subset \mathbb{R}^2$, while the superlevel set $\{\vecfnc{U} : p(\vecfnc{U}) \geq p_{\min}\}$ is also convex because $p(\vecfnc{U})$ is concave, and the intersection of convex sets is itself convex. Therefore, Assumption~\ref{assumption_admissible} is satisfied on any fixed mesh. Unfortunately, this argument is suboptimal, as the resulting uniform positive definiteness constant $\beta$ need not remain bounded away from zero under mesh refinement. Indeed, these bounds are already sufficient to establish $L^2$ boundedness of the conservative variables on a fixed mesh without appealing to the entropy-based estimate \eqref{eq_entropy_l2_bound}. However, because the bounds are mesh-dependent, they are insufficient for establishing a meaningful stability estimate in the limit $h \to 0$. Nevertheless, they do provide some justification for the validity of Assumption~\ref{assumption_admissible} on fixed meshes, and offer some intuition for the conjecture that entropy-stable schemes admit mesh-independent uniform positive definiteness constants $\beta > 0$.

\bibliographystyle{spmpsci}      
\bibliography{references.bib}

\end{document}